\documentclass{amsart}
\usepackage{newclude}
\usepackage{tikz}
\usepackage{tikz-cd}
\usepackage[status=draft]{fixme}
\usepackage{standalone}
\usepackage{amsfonts,amsmath,stmaryrd,txfonts,mathrsfs,amsthm}
\usepackage[all,2cell]{xy}
\usepackage[neveradjust]{paralist}
\usepackage{hyperref}
\usepackage{mathtools}
\usepackage{tikz}
\usetikzlibrary{trees}
\usetikzlibrary{topaths}
\usetikzlibrary{decorations}
\usetikzlibrary{decorations.pathreplacing}
\usetikzlibrary{decorations.pathmorphing}
\usetikzlibrary{decorations.markings}
\usetikzlibrary{matrix,backgrounds,folding}
\usetikzlibrary{chains,scopes,positioning,fit}
\usetikzlibrary{arrows,shadows}
\usetikzlibrary{calc} 
\usetikzlibrary{chains}
\usetikzlibrary{shapes,shapes.geometric,shapes.misc}
\usepackage{smbicat} 
\usepackage{mathpartir}

\makeatletter
\let\ea\expandafter

\def\mdef#1#2{\ea\ea\ea\gdef\ea\ea\noexpand#1\ea{\ea\ensuremath\ea{#2}}}
\def\alwaysmath#1{\ea\ea\ea\global\ea\ea\ea\let\ea\ea\csname your@#1\endcsname\csname #1\endcsname
  \ea\def\csname #1\endcsname{\ensuremath{\csname your@#1\endcsname}}}

\newcommand{\cA}{\ensuremath{\mathcal{A}}}
\newcommand{\cB}{\ensuremath{\mathcal{B}}}
\newcommand{\cC}{\ensuremath{\mathcal{C}}}
\newcommand{\cD}{\ensuremath{\mathcal{D}}}
\newcommand{\cE}{\ensuremath{\mathcal{E}}}

\newcommand{\cH}{\ensuremath{\mathcal{H}}}

\newcommand{\cK}{\ensuremath{\mathcal{K}}}
\newcommand{\cL}{\ensuremath{\mathcal{L}}}
\newcommand{\cM}{\ensuremath{\mathcal{M}}}

\newcommand{\cS}{\ensuremath{\mathcal{S}}}

\newcommand{\cV}{\ensuremath{\mathcal{V}}}

\newcommand{\lC}{\ensuremath{\mathbb{C}}}
\newcommand{\lD}{\ensuremath{\mathbb{D}}}
\newcommand{\lE}{\ensuremath{\mathbb{E}}}

\newcommand{\fa}{\ensuremath{\mathfrak{a}}}

\newcommand{\fl}{\ensuremath{\mathfrak{l}}}

\newcommand{\fr}{\ensuremath{\mathfrak{r}}}
\newcommand{\fs}{\ensuremath{\mathfrak{s}}}

\newcommand{\fu}{\ensuremath{\mathfrak{u}}}

\newcommand{\fx}{\ensuremath{\mathfrak{x}}}

\newcommand{\fB}{\ensuremath{\mathfrak{B}}}
\newcommand{\fC}{\ensuremath{\mathfrak{C}}}

\mdef\fahat{\hat{\fa}}

\newcommand{\fhat}{\ensuremath{\hat{f}}}
\newcommand{\ghat}{\ensuremath{\hat{g}}}

\DeclareSymbolFont{bbold}{U}{bbold}{m}{n}
\DeclareSymbolFontAlphabet{\mathbbb}{bbold}

\newcommand{\alhat}{\ensuremath{\hat{\alpha}}}
\newcommand{\behat}{\ensuremath{\hat{\beta}}}
\newcommand{\gmhat}{\ensuremath{\hat{\gamma}}}
\newcommand{\dehat}{\ensuremath{\hat{\delta}}}

\mdef\del{\partial}
\mdef\delbar{\overline{\partial}}

\newcommand{\inv}{^{-1}}

\mdef\hf{\textstyle\frac{1}{2}}
\mdef\thrd{\textstyle\frac{1}{3}}
\mdef\qtr{\textstyle\frac{1}{4}}

\newcommand{\op}{^{\mathrm{op}}}
\newcommand{\co}{^{\mathrm{co}}}

\let\adj\dashv
\SelectTips{cm}{}
\newdir{ >}{{}*!/-10pt/@{>}}    

\let\iso\cong
\let\eqv\simeq

\mdef\Id{\mathrm{Id}}
\mdef\id{\mathrm{id}}
\alwaysmath{ell}
\alwaysmath{infty}
\alwaysmath{odot}
\def\frc#1/#2.{\frac{#1}{#2}}   
\mdef\ten{\mathrel{\otimes}}
\mdef\bigten{\bigotimes}
\mdef\sqten{\mathrel{\boxtimes}}
\def\pow(#1,#2){\mathop{\pitchfork}(#1,#2)} 

\newcommand{\cat}[1]{\ensuremath{\mathbf{#1}}}
\newcommand{\doub}[1]{\ensuremath{\mathbb{#1}}}

\makeatletter
\def\slashedarrowfill@#1#2#3#4#5{%
  $\m@th\thickmuskip0mu\medmuskip\thickmuskip\thinmuskip\thickmuskip
   \relax#5#1\mkern-7mu%
   \cleaders\hbox{$#5\mkern-2mu#2\mkern-2mu$}\hfill
   \mathclap{#3}\mathclap{#2}%
   \cleaders\hbox{$#5\mkern-2mu#2\mkern-2mu$}\hfill
   \mkern-7mu#4$%
}

\def\Rightslashedarrowfill@{%
  \slashedarrowfill@\Relbar\Relbar\Mapstochar\Rightarrow}
\newcommand\xslashedRightarrow[2][]{%
  \ext@arrow 0055{\Rightslashedarrowfill@}{#1}{#2}}

\newcommand{\Rightthreecell}{\RRightarrow}

\tikzstyle{doubletick}=[-implies, double equal sign distance, postaction={decorate},decoration={markings,mark=at position .5 with {\node at (0,0) {$\bullet$};}}]

\tikzstyle{doubletick}=[-implies, double equal sign distance, postaction={decorate},decoration={markings,mark=at position .5 with {\node at (0,0) {$\bullet$};}}]

\tikzstyle{doubleloose}=[-implies, double equal sign distance, postaction={decorate},decoration={markings,mark=at position .5 with {\node at (0,-0.0) {$\bullet$};}}]

\tikzstyle{vertarrow}=[-implies, ->, postaction={decorate},decoration={markings,mark=at position .5 with {\node at (0,-0.0) {$\bullet$};}}]

\tikzstyle{threecell}=[-implies, double equal sign distance, postaction={decorate},decoration={markings,mark=at position .5 with {\node at (0,-0.0) {$\bullet$};}}]

\tikzstyle{darrow}=[-implies, double equal sign distance, postaction={decorate},decoration={markings,mark=at position .5 with {\node at (0,0) {$\circ$};}}]

\tikzstyle{horarrow}=[-implies, ->, postaction={decorate},decoration={markings,mark=at position .5 with {\node at (0,0) {$\circ$};}}]

\tikzstyle{doubleloose2}=[-implies, double equal sign distance, postaction={decorate},decoration={markings,mark=at position .5 with {\node at (0,-0.01) {$\bullet$};}}] 

\tikzstyle{onecell}=[-implies, postaction={decorate},decoration={markings,mark=at position .5 with {\node at (0,-0.01) {$\circ$};}, mark=at position .90 with {\node at (0,-0.01) {$>$};}}]

\tikzstyle{darrow}=[-implies, double equal sign distance, postaction={decorate},decoration={markings,mark=at position .5 with {\node at (0,0) {$\circ$};}}]

\tikzstyle{doubletight}=[-implies, double equal sign distance]

\tikzstyle{doubleeq}=[-implies, double equal sign distance, postaction={decorate},decoration={markings,mark=at position .5 with {\node at (0,0) {\large $\circ$};}}]

\tikzstyle{doubletighteq}=[double equal sign distance]

\tikzstyle{doubleeq}=[double equal sign distance]

\newcommand{\too}[1][]{\ensuremath{\overset{#1}{\longrightarrow}}}
\newcommand{\ot}{\ensuremath{\leftarrow}}

\mdef\we{\overset{\sim}{\longrightarrow}}
\mdef\leftwe{\overset{\sim}{\longleftarrow}}

\let\maps\colon

\newsavebox{\DDownarrowbox}
\savebox{\DDownarrowbox}{\tikz[scale=1.5]{\draw[-implies,double equal
sign distance] (0,.1) -- (0,-.1); \draw (0,.1) -- (0,-.1);}}
\newcommand{\DDownarrow}{\mathrel{\raisebox{-.2em}{\usebox{\DDownarrowbox}}}}

\newsavebox{\Diagarrowbox}
\savebox{\Diagarrowbox}{\tikz[scale=1.5]{\draw[-implies,double equal sign distance] (.1,.1) -- (-.05,-.05)}}

\newsavebox{\Diagloosearrowbox}
\savebox{\Diagloosearrowbox}{\tikz[scale=1.5]{\draw[-implies,doubleloose] (.1,.1) -- (-.05,-.05)}}

\newsavebox{\RRightarrowbox}
\savebox{\RRightarrowbox}{\tikz[scale=1.5]{\draw[-implies,double equal
sign distance] (-.1,0) -- (.1,0); \draw (-.1,0) -- (.1,0);}}
\newcommand{\RRightarrow}{\mathrel{\raisebox{.2em}{\usebox{\RRightarrowbox}}}}

\newsavebox{\looseRightarrowbox}
\savebox{\looseRightarrowbox}{\tikz[scale=1.5]{\draw[doubleloose2] (-.15,0) -- (.15,0);}}
\newcommand{\looseRightarrow}{\mathrel{\raisebox{-.25em}{\usebox{\looseRightarrowbox}}}}
\newcommand{\xlooseRightarrow}[1]{\mathrel{\overset{#1}{\raisebox{-.25em}{\usebox{\looseRightarrowbox}}}}}

\let\xto\xrightarrow
\let\xot\xleftarrow
\def\rightarrowtailfill@{\arrowfill@{\Yright\joinrel\relbar}\relbar\rightarrow}
\newcommand\xrightarrowtail[2][]{\ext@arrow 0055{\rightarrowtailfill@}{#1}{#2}}

\def\twoheadrightarrowfill@{\arrowfill@{\relbar\joinrel\relbar}\relbar\twoheadrightarrow}
\newcommand\xtwoheadrightarrow[2][]{\ext@arrow 0055{\twoheadrightarrowfill@}{#1}{#2}}

\def\slashedarrowfill@#1#2#3#4#5{%
  $\m@th\thickmuskip0mu\medmuskip\thickmuskip\thinmuskip\thickmuskip
   \relax#5#1\mkern-7mu%
   \cleaders\hbox{$#5\mkern-2mu#2\mkern-2mu$}\hfill
   \mathclap{#3}\mathclap{#2}%
   \cleaders\hbox{$#5\mkern-2mu#2\mkern-2mu$}\hfill
   \mkern-7mu#4$%
}
\def\rightslashedarrowfill@{%
  \slashedarrowfill@\relbar\relbar\mapstochar\rightarrow}
\newcommand\xslashedrightarrow[2][]{%
  \ext@arrow 0055{\rightslashedarrowfill@}{#1}{#2}}
\mdef\hto{\xslashedrightarrow{}}
\mdef\htoo{\xslashedrightarrow{\quad}}
\let\xhto\xslashedrightarrow

\def\defthm#1#2{%
  \newtheorem{#1}{#2}[section]%
  \expandafter\def\csname #1autorefname\endcsname{#2}%
  \expandafter\let\csname c@#1\endcsname\c@thm}
\newtheorem{thm}{Theorem}[section]

\defthm{cor}{Corollary}
\defthm{prop}{Proposition}
\defthm{lem}{Lemma}
\defthm{sch}{Scholium}
\defthm{assume}{Assumption}
\defthm{claim}{Claim}
\defthm{conj}{Conjecture}
\defthm{hyp}{Hypothesis}
\theoremstyle{definition}
\defthm{defn}{Definition}
\defthm{notn}{Notation}
\theoremstyle{remark}
\defthm{rmk}{Remark}
\defthm{eg}{Example}
\defthm{egs}{Examples}
\defthm{ex}{Exercise}
\defthm{ceg}{Counterexample}

\def\thmqedhere{\expandafter\csname\csname @currenvir\endcsname @qed\endcsname}

\let\c@equation\c@thm
\numberwithin{equation}{section}

\mathtoolsset{showonlyrefs,showmanualtags}

\newlength\oldleftmargini       
\newlength\oldleftmarginii
\newlength\oldleftmarginiii
\newlength\oldleftmarginiv
\newlength\oldleftmarginv
\newlength\oldleftmarginvi
\newcount\maxenum
\newif\ifkillspacing
\def\@adjust@enum@labelwidth{%
  \advance\@listdepth by 1\relax
  \ifkillspacing                
    \csname c@\@enumctr\endcsname\maxenum
    \settowidth{\@tempdima}{%
      \csname label\@enumctr\endcsname\hspace{\labelsep}}%
    \csname leftmargin\romannumeral\@listdepth\endcsname
      \@tempdima
  \else                         
    \csname fixspacing\romannumeral\@listdepth\endcsname
  \fi
  \advance\@listdepth by -1\relax}
\def\fixspacingi{\ifnum\oldleftmargini=0\setlength\oldleftmargini\leftmargini\else\setlength\leftmargini\oldleftmargini\fi}
\def\fixspacingii{\ifnum\oldleftmarginii=0\setlength\oldleftmarginii\leftmarginii\else\setlength\leftmarginii\oldleftmarginii\fi}
\def\fixspacingiii{\ifnum\oldleftmarginiii=0\setlength\oldleftmarginiii\leftmarginiii\else\setlength\leftmarginiii\oldleftmarginiii\fi}
\def\fixspacingiv{\ifnum\oldleftmarginiv=0\setlength\oldleftmarginiv\leftmarginiv\else\setlength\leftmarginiv\oldleftmarginiv\fi}
\def\fixspacingv{\ifnum\oldleftmarginv=0\setlength\oldleftmarginv\leftmarginv\else\setlength\leftmarginv\oldleftmarginv\fi}
\def\fixspacingvi{\ifnum\oldleftmarginvi=0\setlength\oldleftmarginvi\leftmarginvi\else\setlength\leftmarginvi\oldleftmarginvi\fi}

\def\pl@label#1#2{%
  \edef\pl@the{\noexpand#1{\@enumctr}}%
  \pl@lab\expandafter{\the\pl@lab\csname yourthe\@enumctr\endcsname}%
  \advance\@tempcnta1
  \pl@loop}
\def\@enumlabel@#1[#2]{%
  \@plmylabeltrue
  \@tempcnta0
  \pl@lab{}%
  \let\pl@the\pl@qmark
  \expandafter\pl@loop\@gobble#2\@@@
  \ifnum\@tempcnta=1\else
    \PackageWarning{paralist}{Incorrect label; no or multiple
      counters.\MessageBreak The label is: \@gobble#2}%
  \fi
  \expandafter\edef\csname label\@enumctr\endcsname{\the\pl@lab}%
  \expandafter\edef\csname the\@enumctr\endcsname{\the\pl@lab}%
  \expandafter\let\csname yourthe\@enumctr\endcsname\pl@the
  #1}

\alwaysmath{alpha}
\alwaysmath{beta}
\alwaysmath{gamma}
\alwaysmath{Gamma}
\alwaysmath{delta}
\alwaysmath{Delta}
\alwaysmath{epsilon}
\mdef\ep{\varepsilon}
\alwaysmath{zeta}
\alwaysmath{eta}
\alwaysmath{theta}
\alwaysmath{Theta}
\alwaysmath{iota}
\alwaysmath{kappa}
\alwaysmath{lambda}
\alwaysmath{Lambda}
\alwaysmath{mu}
\alwaysmath{nu}
\alwaysmath{xi}
\alwaysmath{pi}
\alwaysmath{rho}
\alwaysmath{sigma}
\alwaysmath{Sigma}
\alwaysmath{tau}
\alwaysmath{upsilon}
\alwaysmath{Upsilon}
\alwaysmath{phi}
\alwaysmath{Pi}
\alwaysmath{Phi}
\mdef\ph{\varphi}
\alwaysmath{chi}
\alwaysmath{psi}
\alwaysmath{Psi}
\alwaysmath{omega}
\alwaysmath{Omega}
\let\al\alpha
\let\be\beta
\let\gm\gamma

\let\de\delta

\let\th\theta

\makeatother

\tikzstyle{tickarrow}=[->,postaction={decorate},decoration={markings,mark=at position .5 with {\draw[-] (0,-0.1) -- (0,0.1);}},line width=0.50]

\let\cref\autoref
\UseAllTwocells
\title{Constructing symmetric monoidal bicategories functorially}
\author{Linde Wester Hansen}
\author{Michael Shulman}
\thanks{This material is based on research sponsored by The United States Air Force Research Laboratory under agreement number FA9550-15-1-0053.  The U.S.~Government is authorized to reproduce and distribute reprints for Governmental purposes notwithstanding any copyright notation thereon.  The views and conclusions contained herein are those of the authors and should not be interpreted as necessarily representing the official policies or endorsements, either expressed or implied, of the United States Air Force Research Laboratory, the U.S.~Government, or Carnegie Mellon University.}
\mdef\cMod{\mathcal{M}\mathit{od}}
\mdef\cCat{\mathcal{C}\mathit{at}}
\mdef\cTwocat{2\text{-}\mathcal{C}\mathit{at}}
\mdef\cBicat{\mathcal{B}\mathit{icat}}
\mdef\cCat{\mathcal{C}\mathit{at}}
\mdef\cVect{\mathcal{V}\mathit{ect}}
\mdef\cHilb{\mathcal{H}\mathit{ilb}}
\mdef\cProf{\mathcal{P}\mathit{rof}}
\mdef\cSpan{\mathcal{S}\mathit{pan}}
\mdef\cMon{\mathcal{M}\mathit{on}}
\mdef\cBr{\mathcal{B}\mathit{r}}
\mdef\cSyl{\mathcal{S}\mathit{yl}}
\mdef\cSym{\mathcal{S}\mathit{ym}}
\mdef\cH{\mathcal{L}}
\mdef\fBicat{\mathfrak{Bicat}}
\mdef\fDbl{\mathfrak{Dbl}}
\mdef\fH{\mathfrak{H}}
\mdef\fC{\mathfrak{T}}
\mdef\lMod{\mathbb{M}\mathsf{od}}
\mdef\lnCob{n\mathbb{C}\mathsf{ob}}
\let\ltwo\bbtwo
\mdef\lProf{\mathbb{P}\mathsf{rof}}
\mdef\lSpan{\mathbb{S}\mathsf{pan}}
\mdef\lMat{\mathbb{M}\mathsf{at}}
\mdef\cDbl{\mathcal{D}\mathit{bl}}
\mdef\cDblcf{\cDbl_{c\mathbf{f}}}
\mdef\cDbllf{\cDbl_{l\mathbf{f}}}
\mdef\cDblf{\cDbl_{\mathbf{f}}}
\mdef\fDblf{\mathfrak{Dbl}_{\mathbf{f}}}
\mdef\cMonDbll{\mathcal{M}\mathit{on}\mathcal{D}\mathit{bl}_l}
\mdef\cMonDblc{\mathcal{M}\mathit{on}\mathcal{D}\mathit{bl}_c}
\mdef\cMonDblp{\mathcal{M}\mathit{on}\mathcal{D}\mathit{bl}_p}
\mdef\fchk{\check{f}}
\mdef\conj{\Yleft}
\mdef\Conj{\mathcal{C}\mathit{onj}}
\mdef\Icon{\mathcal{I}\mathit{con}}
\mdef\id{\mathsf{id}}
\mdef\C{\mathbb{C}}
\mdef\D{\mathbb{D}}
\mdef\E{\mathbb{E}}
\mdef\F{\mathbb{F}}
\newcommand{\bAlg}{\mathbb{A}\mathsf{lg}}
\let\lAlg\bAlg
\newcommand{\cAlg}{\mathcal{A}\mathit{lg}}

\newcommand{\tens}{\otimes}
\newcommand{\onecell}{\rightarrow}

\newcommand{\verc}{\cdot}
\newcommand{\horc}{\bullet}
\newcommand{\comp}{\circ}
\newcommand{\looseid}{\Id}
\newcommand{\tightid}{1}
\newcommand{\transid}{\id}
\newcommand{\mult}{\tinymult[gray dot]}
\newcommand{\unit}{\tinyunit[gray dot]}

\newcommand{\lop}{^{l\cdot\mathrm{op}}}
\newcommand{\ttop}{^{t\cdot\mathrm{op}}}
\newcommand{\tlop}{^{tl\cdot\mathrm{op}}}
\newcommand{\lco}{^{l\cdot\mathrm{co}}}
\newcommand{\tco}{^{t\cdot\mathrm{co}}}
\newcommand{\tlco}{^{tl\cdot\mathrm{co}}}
\newcommand{\lcoop}{^{l\cdot\mathrm{coop}}}

\newcommand{\hora}{a^{\horc}}
\newcommand{\horr}{r^{\horc}}
\newcommand{\horl}{l^{\horc}}
\newcommand{\compa}{a^{\comp}}
\newcommand{\compr}{r^{\comp}}
\newcommand{\compl}{l^{\comp}}

\newcounter{mondbl}             

\begin{document}

\begin{abstract}
  We present a method of constructing monoidal, braided monoidal, and symmetric monoidal bicategories from corresponding types of monoidal double categories that satisfy a lifting condition.
  Many important monoidal bicategories arise naturally in this way, and applying our general method is much easier than explicitly verifying the coherence laws of a monoidal bicategory for each example.
  Abstracting from earlier work in this direction, we express the construction as a functor between locally cubical bicategories that preserves monoid objects; this ensures that it also preserves monoidal functors, transformations, adjunctions, and so on.
  Examples include the monoidal bicategories of algebras and bimodules, categories and profunctors, sets and spans, open Markov processes, parametrized spectra, and various functors relating them.
\end{abstract}

\maketitle
\setcounter{tocdepth}{1}
\tableofcontents

\section{Introduction}
\label{sec:introduction}

Symmetric monoidal bicategories are important in many contexts.
However, the definition of even a monoidal bicategory
(see~\cite{gps:tricats,nick:tricats}), let alone a symmetric monoidal
one
(see~\cite{kv:2cat-zam,kv:bm2cat,bn:hda-i,ds:monbi-hopfagbd,crans:centers,mccrudden:bal-coalgb,gurski:brmonbicat}),
a monoidal functor between such (see~\cite{nick:tricatsbook,mccrudden:bal-coalgb}),
or a monoidal transformation or modification (see~\cite{sp:thesis})
is quite imposing, and time-consuming to verify in any example.

In this paper we describe a method for constructing (symmetric) monoidal
bicategories, as well as functors and transformations between them, which is hardly more difficult than constructing a pair
of ordinary (symmetric) monoidal categories.
While not universally applicable, this method applies in many cases of interest.
The underlying idea has often been implicitly used in particular cases, such as
bicategories of enriched profunctors, but to our knowledge the first
general statement was claimed in~\cite[Appendix B]{shulman:frbi}.
In the unpublished~\cite{shulman:smbicat}, the first author worked out the details for the construction of monoidal bicategories themselves.
Here we include that work and build on it further to construct monoidal functors, transformations, and so on between monoidal bicategories as well, making the entire construction into a functor.%
\footnote{See also~\cite[\S5]{gg:ldstr-tricat} which generalizes the construction in a different direction, showing that every sufficiently nice locally cubical bicategory has an underlying tricategory.}

The method relies on the fact that in many bicategories, the 1-cells
are not the most fundamental notion of `morphism' between the objects.
For instance, in the bicategory \cMod\ of rings, bimodules, and
bimodule maps, the more fundamental notion of morphism between objects
is a ring homomorphism. The addition of these extra morphisms promotes
a bicategory to a \emph{double category}, or a category internal to
\cCat.  The extra morphisms are usually stricter than the 1-cells in
the bicategory and easier to deal with for coherence questions; in
many cases it is quite easy to show that we have a \emph{symmetric
  monoidal double category}.  The central observation is that in most
cases (when the natural transformations have `loosely strong companions') we can then `lift' this
symmetric monoidal structure to the original bicategory.  That is, we
prove the following theorem:

\begin{thm}\label{thm:mondbl-monbi-intro}
  If \lD\ is a monoidal double category, of which the monoidal constraints have loosely strong companions, then its underlying bicategory $\cH(\lD)$ is a monoidal bicategory.  If \lD\ is braided
  or symmetric, so is $\cH(\lD)$.
\end{thm}

In~\cite{shulman:smbicat} this theorem was proven by explicitly constructing liftings of all the coherence data, but in the present paper we take a more functorial viewpoint:
we extend the operation $\cH$ that takes a double category to its underlying bicategory to a suitable sort of ``functor'', and show that this functor is product-preserving.
Thus, just as a product-preserving functor between ordinary categories automatically preserves not just internal monoids but also monoid homomorphisms, the functor $\cH$ preserves monoidal objects as well as functors, transformations, and so on between them.
In fact, we will show that $\cH$ induces another functor from a ``category'' of monoidal double categories to one of monoidal bicategories, thereby preserving all kinds of composition as well.

The tricky part is deciding into what kind of categorical structure we should assemble our double categories and bicategories, and thus what kind of functor $\cH$ should be.
To start with, double categories most naturally assemble into a strict 2-category, while bicategories most naturally form a tricategory; and since a 2-category can be considered a degenerate tricategory, we could work with tricategories all the way through.
However, tricategories are really too weak for our purposes.
On the one hand, manipulating all the coherences in a tricategory, let alone a functor between tricategories, is exceedingly difficult.
On the other hand, even the tricategory of bicategories is considerably stricter than an arbitrary tricategory.
In addition to suggesting that stricter alternatives are available, this also means that if we treated bicategories as forming a fully weak tricategory, then an internal notion of ``monoid'' in that tricategory would not coincide exactly with a monoidal bicategory as usually defined, but would have extra unnecessary coherences added, making for yet more work in relating such a definition to the now-accepted one.

The first alternative to tricategories one might consider is what was called in~\cite{shulman:psalg} an \emph{iconic tricategory}.
This is a tricategory-like structure whose coherences for composition along 0-cells are \emph{icons}~\cite{lack:icons} rather than fully general pseudonatural transformations; more precisely it is a bicategory enriched over the monoidal 2-category of bicategories, pseudofunctors, and icons.
Informally, an iconic tricategory is one where composition of 1-cells along 0-cells is strictly associative and unital (though composition of 2-cells along 0-cells need not be).
The tricategory of bicategories is indeed iconic, as of course is the strict 2-category of double categories regarded as a tricategory, and $\cH$ can be made into a product-preserving ``iconic functor'' between them.
However, while the notion of iconic tricategory suffices for our \emph{input} data, it is insufficient for our \emph{output} data: the tricategory of monoidal bicategories is not iconic.\footnote{The same problem occurs for the ``bicategory-enriched categories'' of~\cite{verity:base-change}.}

There is, however, a stricter structure than tricategories that does encompass monoidal bicategories: a bicategory enriched over \emph{double} categories, introduced in~\cite{gg:ldstr-tricat} under the name \emph{locally cubical bicategory}.
In addition to 0-cells and 1-cells, a locally cubical bicategory has \emph{two} kinds of 2-cells, as well as 3-cells inhabiting a square boundary of 2-cells.
Just as a bicategory can be regarded as a double category that is trivial in one direction, an iconic tricategory can be regarded as a locally cubical bicategory --- although in the case of bicategories, it is more natural to take the additional kind of 2-cells to be icons.
Similarly, it is shown in~\cite{gg:ldstr-tricat} that monoidal bicategories form a locally cubical bicategory, with an appropriate notion of ``monoidal icon'' as the additional 2-cells.

Thus, we can hope to show that in general, internal monoids in a locally cubical bicategory with products form another locally cubical bicategory; that this reproduces the standard notions of monoidal double category and monoidal bicategory; and that any product-preserving functor between locally cubical bicategories with products induces another functor between their locally cubical bicategories of internal monoids (of all sorts).
But unfortunately, here the extra generality of locally cubical bicategories becomes a problem: the usual definition of monoidal bicategory in fact \emph{relies} on the additional strictness of bicategories, and doesn't make sense in the generality of any locally cubical bicategory.
We expect this could be dealt with using ``local fibrancy'' conditions, but for simplicity we will instead assume that our \emph{input} locally cubical bicategories and functors are ``1-strict'', i.e.\ they satisfy roughly the same strictness property of an iconic tricategory (composition of 1-cells is strictly associative), even though our output ones may not be.
Thus, we finally obtain our desired statement:

\begin{thm}\label{thm:functor-intro}
  The assignment $\cH$ extends to a functor between the locally cubical bicategories of monoidal, braided, and symmetric double categories and bicategories.
  In particular, it preserves monoidal functors, monoidal transformations, and composites thereof.
\end{thm}

In fact, we actually prove several theorems of this sort, depending on whether the monoidal functors and transformations in question are chosen to be \emph{lax}, \emph{colax}, or \emph{strong}, i.e.\ whether they preserve the monoidal structure up to a transformation in one direction, the other direction, or an invertible transformation.
This distinction for functors is already known for ordinary monoidal categories; for monoidal bicategories such a threefold choice is also available for transformations.
Note that this laxity is only relative to the monoidal structure: on the underlying bicategories, all our functors and transformations will be strong/pseudo, preserving composition up to invertible transformations.

One might hope to incorporate both lax and colax functoriality in a single theorem.
For instance, as noted in~\cite{gp:double-adjoints,shulman:dblderived}, lax and colax morphisms \emph{themselves} form the two kinds of morphisms in a double category!
A functoriality theorem at this level would have the advantage of also preserving ``mates'' in this double category, including for instance doctrinal adjunctions~\cite{kelly:doc-adjn} between lax and colax monoidal functors.
However, this would require a more complicated 3-dimensional structure such as that of~\cite{gp:intercategories-i}; we do not pursue it here.

\begin{rmk}
We also expect similar theorems to be true in higher dimensions.  For
instance, Chris Douglas~\cite{douglas:tfttalk} has suggested that many
apparent tricategories are more naturally bicategories internal to
\cCat\ or categories internal to \cTwocat; and in most such cases
arising in practice, we can again `lift' the coherence to give a
tricategory.
We propose the term \textbf{$(n\times k)$-category}
(pronounced ``$n$-by-$k$-category'') for an $n$-category internal to
$k$-categories, which has $(n+1)(k+1)$ different types of
cells in an $(n+1)$ by $(k+1)$ grid.  Thus
double categories are \textbf{1x1-categories}, while in
place of tricategories we may consider 2x1-categories and
1x2-categories --- or even 1x1x1-categories, i.e.\ triple categories,
as in~\cite{gp:intercategories-i,gp:intercategories-ii}.%
\footnote{If we also generalize the traditional terminology to say that an $(n+1)$-category is a category \emph{enriched} in $n$-categories even when $n$ is a symbol more general than a natural number, we could say that locally cubical bicategories are ``$((1\times 1)+1)$-categories''.}
Any
$(n\times k)$-category with a suitable lifting property
should have an underlying $(n+k)$-category, but this discards an increasing amount of structure as $n$ and $k$ grow.


There is a case to be made that often the extra cells should
\emph{not} be discarded.  But sometimes
it really is the underlying $(n+k)$-category one cares about; for
instance, the Baez-Dolan cobordism hypothesis is about the $(n+1)$-category of cobordisms, not
the $(n\times 1)$-category from which it is constructed
(see~\cite{lurie:tft}).  Thus we believe there is an indisputable value to
results such as \autoref{thm:mondbl-monbi-intro} and \autoref{thm:functor-intro}.
\end{rmk}

Proceeding to the contents of this paper, in
\S\ref{sec:symm-mono-double} we review the definition of symmetric
monoidal double categories, and in \S\ref{sec:comp-conj} we recall the
notions of `companion' and `conjoint' whose presence supplies the
necessary lifting property and prove their essential uniqueness.
(Double categories with companions and conjoints have also been called ``framed bicategories''~\cite{shulman:frbi}, and are roughly equivalent to ``proarrow equipments''~\cite{wood:proarrows-i}).
Then in \S\ref{sec:1x1-to-bicat} we show that taking the underlying bicategory of a double category defines a functor from the locally cubical bicategory of double categories with ``loosely strong companions'' to the locally cubical bicategory of bicategories.

In order to prove that this functor preserves the monoidal structure, in \S\ref{sec:mono-objects} we define monoidal structures abstractly on objects of a locally cubical bicategory, and we show that the monoidal objects and cells form a new locally cubical bicategory. Furthermore, we prove that any 1-strict product-preserving functor between 1-strict locally cubical bicategories preserves monoidal objects and cells of all sorts, and indeed induces another functor of locally cubical bicategories.
There is a lot to check here, but the hardest part is writing down all the definitions in the appropriate generality! Then we specialize this to the functor from double categories to bicategories, yielding functorial constructions of monoidal, braided, and symmetric monoidal bicategories.

Finally, in \S\ref{sec:Alg} we illustrate our method with a number of examples, including algebras and bimodules, categories and profunctors, dagger Frobenius algebras and bimodules, black-boxing of open Markov processes, cartesian bicategories and cartesian double categories, and parametrized spectra.

We would like to thank Peter May, Tom Fiore, Stephan Stolz, Chris
Douglas, Nick Gurski, Jamie Vicary, Julian Hedges, and Bob Par\'{e} for helpful discussions and comments.


\section{Symmetric monoidal double categories}
\label{sec:symm-mono-double}

In this section, we recall basic notions of double categories to fix our terminology and notation, and define monoidal double categories and functors between them in an explicit way.
Double categories go back originally to Ehresmann
in~\cite{ehresmann:cat-str}; a brief introduction can be found
in~\cite{ks:r2cats}.  Other references
include~\cite{multi_funct_i,gp:double-limits,gp:double-adjoints,aleiferi2018cartesian}.

\begin{defn}\label{def:dblcat}
  A \textbf{(pseudo) double category} \lD\ consists of a `category of
  objects' $\lD_0$ and a `category of arrows' $\lD_1$, with structure
  functors
  \begin{align*}
    U&\maps \lD_0\to \lD_1\\
    S,T&\maps \lD_1\rightrightarrows \lD_0\\
    \odot&\maps \lD_1\times_{\lD_0}\lD_1\to \lD_1
  \end{align*}
  (where the pullback is over
  $\lD_1\too[T]\lD_0\overset{S}{\longleftarrow} \lD_1$) such that
  \begin{alignat*}{2}
    S(U_A) &= A &\qquad
    S(M\odot N) &= SN\\
    T(U_A) &= A &\qquad
    T(M\odot N) &= TM
  \end{alignat*}
  naturally, and equipped with natural isomorphisms
  \begin{align*}
    \fa &: (M\odot N) \odot P \too[\iso] M \odot (N \odot P)\\
    \fl &: U_B \odot M \too[\iso] M\\
    \fr &: M \odot U_A \too[\iso] M
  \end{align*}
  such that $S(\fa)$, $T(\fa)$, $S(\fl)$, $T(\fl)$, $S(\fr)$, and
  $T(\fr)$ are all identities, and such that the standard coherence
  axioms for a monoidal category or bicategory (such as Mac Lane's
  pentagon; see~\cite{maclane}) are satisfied.
\end{defn}

Just as a bicategory can be thought of as a category weakly
\emph{enriched} over \cCat, a pseudo double category can be thought of
as a category weakly \emph{internal} to \cCat.  We will usually omit
the adjective ``pseudo'' for conciseness.

We call the objects of $\lD_0$ \textbf{objects} or \textbf{0-cells},
and we call the morphisms of $\lD_0$ \textbf{tight 1-cells}
and write them as $f\maps A\to B$.  We call the objects of $\lD_1$
\textbf{loose 1-cells}; if $M$ is a 1-cell with $S(M)=A$ and
$T(M)=B$, we write $M\maps A\hto B$.  We call a morphism $\alpha\maps
M\to N$ of $\lD_1$ with $S(\alpha)=f$ and $T(\alpha)=g$ a
\textbf{2-cell} and draw it as follows:
\begin{equation}\label{eq:square}
  \xymatrix@-.5pc{
    A \ar[r]|{|}^{M}  \ar[d]_f \ar@{}[dr]|{\Downarrow\alpha}&
    B\ar[d]^g\\
    C \ar[r]|{|}_N & D
  }.
\end{equation}
We will sometimes say ``$k$-morphism'' as a synonym for ``$k$-cell'' (tight or loose).
Note that composition of tight 1-cells is strictly associative and unital, while that of loose 1-cells is only weakly so.
This is the case in the majority of examples, and is part of what makes monoidal double categories so much easier to construct than monoidal bicategories.
(It is, however, possible to define double categories that are weak in both directions~\cite{verity:base-change}.)

The words ``tight'' and ``loose'', borrowed from~\cite{ls:limlax}, are intended to suggest one kind of morphism that is ``stricter'' than another.
Traditionally the two kinds of 1-arrow in a double category are called ``vertical'' and ``horizontal'' with reference to how they are drawn, but this creates confusion because some authors draw the tight morphisms (those with strictly associative composition) vertically and others draw them horizontally.
We usually draw the tight 1-cells vertically and the loose ones horizontally, but the terminology ``tight'' and ``loose'' unambiguously identifies which arrows we are talking about, independently of our conventions about how to draw pictures.%
\footnote{In~\cite{shulman:smbicat} a similar effect was intended by simply distinguishing between ``1-morphisms'' (the tight ones) and ``1-cells'' (the loose ones), but since ``morphism'' and ``cell'' (and ``arrow'') are generally used interchangeably in higher category theory this was not very successful.}

We write the composition of tight 1-morphisms $A\too[f] B\too[g] C$
and the tight composition of 2-morphisms $M\too[\alpha]
N\too[\beta] P$ as $g\circ f$ and $\beta\circ\alpha$, or sometimes
just $gf$ and $\beta\alpha$.  We write the loose composition of
1-cells $A\xhto{M} B \xhto{N} C$ as $A\xhto{N\odot M} C$ and that of
2-morphisms
\[\vcenter{\xymatrix{ \ar[r]|-@{|}^-{} \ar[d] \ar@{}[dr]|{\Downarrow\alpha} &
     \ar[r]|-@{|}^-{} \ar[d] \ar@{}[dr]|{\Downarrow\beta} &\ar[d]\\
  \ar[r]|-@{|}_-{} & \ar[r]|-@{|}_-{} & }}
\qquad\text{as}\qquad
\vcenter{\xymatrix@C=4pc{ \ar[r]|-@{|}^-{} \ar[d] \ar@{}[dr]|{\Downarrow\;\be\odot\al} &  \ar[d]\\
  \ar[r]|-@{|}_-{} & }}.\]

The two different compositions of 2-morphisms obey an interchange law,
by the functoriality of $\odot$:
\[(M_1\odot M_2) \circ (N_1\odot N_2) = (M_1\circ N_1)\odot (M_2\circ N_2).
\]
Every object $A$ has a tight identity $1_A$ and a loose unit
$U_A$, every loose 1-cell $M$ has an identity 2-morphism $1_M$, every
tight 1-morphism $f$ has a loose unit 2-morphism $U_f$, and we
have $1_{U_A} = U_{1_A}$ (by the functoriality of $U$).



A 2-morphism~\eqref{eq:square} where $f$ and $g$ are identities (such
as the constraint isomorphisms $\fa,\fl,\fr$) is called
\textbf{globular}.  Every double category \lD\ has a
\textbf{loose bicategory} $\cH(\lD)$ consisting of the objects,
loose 1-cells, and globular 2-morphisms.  In the literature, this is often called the ``horizontal'' or ``vertical'' bicategory of a double category (depending on conventions). Conversely, many naturally
occurring bicategories are actually the loose bicategory of some
naturally ocurring double category.  Here are just a few examples.

\begin{eg}
  The double category \lnCob\ has as objects closed $n$-manifolds, as
  tight 1-cells diffeomorphisms, as loose 1-cells cobordisms, and as
  2-cells diffeomorphisms between cobordisms.  Its loose bicategory $\cH(\lnCob)$
  is the usual bicategory of cobordisms.
\end{eg}

\begin{eg}
  The double category \lMod\ has as objects rings, as tight 1-cells ring
  homomorphisms, as loose 1-cells bimodules, and as 2-cells equivariant
  bimodule maps.  Its loose bicategory $\cMod = \cH(\lMod)$ is
  the usual bicategory of rings and bimodules. 
\end{eg}

\begin{eg}
  The double category \lProf\ has as objects categories, as
  tight 1-cells functors, as loose 1-cells \emph{profunctors} (a profunctor
  $A\hto B$ is a functor $B\op\times A\to \mathbf{Set}$), and as
  2-cells natural transformations.  Bicategories such as
  $\cH(\lProf)$ are commonly encountered in category theory,
  especially the enriched versions.
\end{eg}

\begin{eg}
  The double category \lSpan\ has as objects sets, as tight 1-cells functions, as loose 1-cells spans, and as 2-cells maps of spans over pairs of functions.  Its loose bicategory $\cH(\lSpan)$ is the usual bicategory of spans.
\end{eg}

Further examples will be found in Section~\ref{sec:Alg}.


As opposed to bicategories, which naturally form a tricategory, double
categories naturally form a \emph{2-category}, a much simpler object.

\begin{defn}
  Let \lD\ and \lE\ be double categories.  A \textbf{(pseudo double)
    functor} $F \maps \lD\to \lE$ consists of the following.
  \begin{itemize}
  \item Functors $F_0\maps \lD_0 \to \lE_0$ and $F_1\maps \lD_1 \to
    \lE_1$ such that $S\circ F_1 = F_0\circ S$ and $T\circ F_1 =
    F_0\circ T$.
  \item Natural transformations $F_\odot\maps F_1M \odot F_1N \to
    F_1(M\odot N)$ and $F_U\maps U_{F_0 A} \to F_1(U_A)$, whose
    components are globular isomorphisms, and which satisfy the usual
    coherence axioms for a monoidal functor or pseudofunctor
    (see~\cite[\S{}XI.2]{maclane}).
  \end{itemize}
\end{defn}

\begin{defn}\label{thm:dbl-transf}
  A \textbf{(tight) transformation} between two functors $\alpha:
  F\to G:\lD\to\lE$ consists of natural transformations $\alpha_0\maps
  F_0\to G_0$ and $\alpha_1\maps F_1\to G_1$ (both usually written as
  $\alpha$), such that $S(\alpha_{M}) = \alpha_{SM}$ and
  $T(\alpha_{M}) = \alpha_{TM}$, and such that
  \[\vcenter{\xymatrix@-.5pc{
      FA \ar@{=}[d] \ar[r]|{|}^{FM}
      \ar@{}[drr]|{\Downarrow F_\odot} &
      FB \ar[r]|{|}^{FN} &
      FC \ar@{=}[d]\\
      FA \ar[rr]|{F(N\odot M)} \ar[d]_{\alpha_A}
      \ar@{}[drr]|{\Downarrow \alpha_{N\odot M}} &&
      FC \ar[d]^{\alpha_C}\\
      GA \ar[rr]|{|}_{G(N\odot M)} && GC
    }} =
  \vcenter{\xymatrix@-.5pc{
      FA \ar[d]_{\alpha_A} \ar@{}[dr]|{\Downarrow \alpha_M} \ar[r]|{|}^{FM} &
      FB \ar[d]|{\alpha_B} \ar@{}[dr]|{\Downarrow \alpha_N} \ar[r]|{|}^{FN} &
      FC \ar[d]^{\alpha_C}\\
      GA \ar@{=}[d] \ar[r]|{|}_{GM} \ar@{}[drr]|{\Downarrow G_\odot} &
      GB \ar[r]|{|}_{GN} &
      GC \ar@{=}[d]\\
      GA \ar[rr]|{|}_{G(N\odot M)} && GC
    }}\]
  for all 1-cells $M\colon A\hto B$ and $N\colon B\hto C$, and
  \[\vcenter{\xymatrix@-.5pc{
      FA \ar[rr]|{|}^{U_{FA}} \ar@{=}[d]
      \ar@{}[drr]|{\Downarrow F_U} &&
      FA \ar@{=}[d]\\
      FA \ar[rr]|{F(U_A)} \ar[d]_{\alpha_A}
      \ar@{}[drr]|{\Downarrow \alpha_{U_A}} &&
      FA \ar[d]^{\alpha_A}\\
      GA \ar[rr]|{|}_{G(U_A)} && GA
    }} =
  \vcenter{\xymatrix@-.5pc{
      FA \ar[rr]|{|}^{U_{FA}} \ar[d]_{\alpha_A}
      \ar@{}[drr]|{\Downarrow U_{\alpha_A}} &&
      FA \ar[d]^{\alpha_A}\\
      GA \ar[rr]|{U_{GA}} \ar@{=}[d]
      \ar@{}[drr]|{\Downarrow G_U} &&
      GA \ar@{=}[d]\\
      GA \ar[rr]|{|}_{G(U_A)} && GA.
    }}\]
  for all objects $A$.
\end{defn}

We write \cDbl\ for the strict 2-category of double categories, functors, and transformations, and $\mathbf{Dbl}$ for its underlying 1-category. As pseudofunctors compose strictly associatively, this is well-defined. Note that a 2-cell $\al$ in \cDbl\ is an isomorphism just when each
$\al_A$, \emph{and} each $\al_M$, is invertible.

\begin{rmk}\label{rmk:dblop}
  A double category has three different opposites:
  in the \emph{loose opposite} $\lD\lop$ we reverse the loose 1-cells (and reverse the 2-cells in the loose direction) but not the tight ones, in the \emph{tight opposite} $\lD\ttop$ we reverse the tight ones but not the loose ones, and in the \emph{double opposite} $\lD\tlop$ we reverse both.
  These define 2-functors $(-)\lop : \cDbl\to\cDbl$, $(-)\ttop : \cDbl\co \to \cDbl$, and $(-)\tlop : \cDbl\co \to \cDbl$, where $(-)\co$ denotes reversal of 2-cells but not 1-cells in a 2-category.
\end{rmk}

The 2-category \cDbl\ gives us an easy way to define what we mean by a
\emph{symmetric monoidal double category}.
We say that a strict 2-category $\cK$ has \emph{finite products} if it has a strictly terminal object, i.e.\ an object $1$ such that each category $\cK(D,1)$ is \emph{isomorphic} to the terminal category, and any two objects $D,E$ have a strict product, i.e.\ a span $D \ot D\times E \to E$ such that each induced functor $\cK(X,D\times E) \to \cK(X,D) \times \cK(X,E)$ is an \emph{isomorphism} of categories.
In any such 2-category with finite products there is a notion of a \emph{pseudomonoid} (perhaps braided or symmetric), which generalizes the notion of monoidal category in \cCat.
We omit the general definition, since we will give a more general one in \cref{sec:mono-objects};
for now we just specialize it to \cDbl\ and obtain the following.

\begin{defn}\label{def:symmondoub}
  A \textbf{monoidal double category} is a double category equipped
  with functors $\ten\maps \lD\times\lD\to\lD$ and $I\maps * \to\lD$,
  and invertible transformations
  \begin{align*}
    \mathord{\otimes} \circ (\Id\times \mathord{\otimes})
    &\iso \mathord{\otimes} \circ (\mathord{\otimes} \times \Id)\\
    \mathord{\otimes} \circ (\Id\times I) &\iso \Id\\
    \mathord{\otimes} \circ (I\times \Id) &\iso \Id
  \end{align*}
  satisfying the usual axioms.  If it additionally has a braiding
  isomorphism
  \begin{align*}
    \mathord{\otimes} &\iso \mathord{\otimes} \circ \tau
  \end{align*}
  (where $\tau\maps \lD\times\lD\iso \lD\times\lD$ is the twist)
  satisfying the usual axioms, then it is \textbf{braided} or
  \textbf{symmetric}, according to whether or not the braiding is
  self-inverse.
\end{defn}

Unpacking this definition more explicitly, we see that a monoidal
double category is a double category together with the following
structure.
\begin{enumerate}
\item $\lD_0$ and $\lD_1$ are both monoidal categories.
\item If $I$ is the monoidal unit of $\lD_0$, then $U_I$ is the
  monoidal unit of $\lD_1$.\footnote{Actually, all the above
    definition requires is that $U_I$ is coherently \emph{isomorphic
      to} the monoidal unit of $\lD_1$, but we can always choose them
    to be equal without changing the rest of the structure.
    More precisely, we may choose the unit isomorphism
    $I_U : U_{I_{\lD_0}} \to I_{\lD_1}$ of the functor $I:* \to \lD$
    to be an identity.}
\item The functors $S$ and $T$ are strict monoidal, i.e.\ $S(M\ten N)
  = SM\ten SN$ and $T(M\ten N)=TM\ten TN$ and $S$ and $T$ also
  preserve the associativity and unit constraints.
\item \label{eq:mon1} We have globular isomorphisms derived from $\otimes_{\odot}$ and $\otimes_U$
  \[\fx\maps (M_1\ten N_1)\odot (M_2\ten N_2)\too[\iso] (M_1\odot M_2)\ten (N_1\odot N_2)\]
  and
  \[\fu\maps U_{A\ten B} \too[\iso] (U_A \ten U_B)\]
  such that the following diagrams commute, expressing that $\ten$ is a pseudo double functor:
\begin{equation}\label{eq:mondoub1}
\begin{aligned}
\begin{tikzpicture}[xscale=1.8, yscale=1.5]
\node (tl) at (0,2) {$((M_1 \tens N_1)\odot (M_2 \tens N_2)) \odot (M_3 \tens N_3)$};
\node (tr) at (4,2) {$((M_1 \odot M_2) \tens (N_1 \odot N_2)) \odot (M_3 \tens N_3)$};
\node (ml) at (0,1) {$(M_1 \tens N_1) \odot ((M_2 \tens N_2) \odot (M_3 \tens N_3))$};
\node (mr) at (4,1) {$((M_1 \odot M_2) \odot M_3) \tens ((N_1 \odot N_2) \odot N_3)$};
\node (bl) at (0,0) {$(M_1 \tens N_1) \odot ((M_2 \odot M_3) \tens (N_2 \odot N_3))$};
\node (br) at (4,0) {$(M_1 \odot (M_2 \odot M_3)) \tens (N_1 \odot (N_2 \odot N_3))$};
\draw[->] (tl) to node[above] {$\fx \odot \id$} (tr);
\draw[->] (tl) to node[left]{$\fa$} (ml);
\draw[->] (ml) to node[left]{$\id\odot \fx$} (bl);
\draw[->] (tr) to node[left]{$\fx$} (mr);
\draw[->] (mr) to node[left]{$\fa \odot \fa$} (br);
\draw[->] (bl) to node[above] {$\fx$} (br);
\end{tikzpicture}
    \end{aligned}
\end{equation}


\begin{equation}\label{eq:mondoub2}
\begin{aligned}
\begin{tikzpicture}[xscale=2, yscale=1.5]
\node (tl) at (0,2) {$(M \tens N) \odot U_{C \tens D}$};
\node (tr) at (4,2) {$(M \tens N) \odot ( U_C \tens U_D)$};
\node (ml) at (0,1) {$M \tens N$};
\node (mr) at (4,1) {$(M \odot U_C) \tens (N \odot U_D)$};
\draw[->] (tl) to node[above] {$\id \odot \fu$} (tr);
\draw[->] (tl) to node[left]{$\fr$} (ml);
\draw[->] (tr) to node[left]{$\fx$} (mr);
\draw[->] (mr) to node[above] {$\fr \tens \fr$} (ml);
\end{tikzpicture}
    \end{aligned}
\end{equation}

    
    \begin{equation}\label{eq:mondoub3}
\begin{aligned}
\begin{tikzpicture}[xscale=2, yscale=1.5]
\node (tl) at (0,2) {$U_{A \tens B} \odot (M \tens N) $};
\node (tr) at (4,2) {$(U_A \tens U_B) \odot (M \tens N) $};
\node (ml) at (0,1) {$M \tens N$};
\node (mr) at (4,1) {$(U_A \odot M ) \tens (U_B \odot N)$};
\draw[->] (tl) to node[above] {$\fu \odot \id$} (tr);
\draw[->] (tl) to node[left]{$\fl$} (ml);
\draw[->] (tr) to node[left]{$\fx$} (mr);
\draw[->] (mr) to node[above] {$\fl \tens \fl$} (ml);
\end{tikzpicture}
    \end{aligned}
\end{equation}

\item \label{eq:mon2}The following diagrams commute, expressing that the
  associativity isomorphism for $\ten$ is a transformation of double
  categories.
\begin{equation}\label{eq:mondoub4}
\begin{aligned}
\begin{tikzpicture}[xscale=1.8, yscale=1.5]
\node (tl) at (0,2) {$((M_1 \tens N_1) \tens P_1) \odot ((M_2 \tens N_2) \tens P_2) $};
\node (tr) at (4,2) {$(M_1 \tens (N_1 \tens P_1)) \odot (M_2 \tens (N_2 \tens P_2))$};
\node (ml) at (0,1) {$((M_1 \tens N_1) \odot (M_2 \tens N_2)) \tens (P_1 \odot P_2)$};
\node (mr) at (4,1) {$(M_1 \odot M_2)  \tens ((N_1 \tens P_1) \odot (N_2 \tens P_2))$};
\node (bl) at (0,0) {$((M_1 \odot M_2) \tens (N_1 \odot N_2)) \tens (P_1 \odot P_2)$};
\node (br) at (4,0) {$(M_1 \odot M_2) \tens ((N_1 \odot N_2) \tens (P_1 \odot P_2))$};
\draw[->] (tl) to node[above] {$\fa \odot \fa$} (tr);
\draw[->] (tl) to node[left]{$\fx$} (ml);
\draw[->] (ml) to node[left]{$ \fx \tens \id$} (bl);
\draw[->] (tr) to node[left]{$\fx$} (mr);
\draw[->] (mr) to node[left]{$\id \tens \fx$} (br);
\draw[->] (bl) to node[above] {$\fa$} (br);
\end{tikzpicture}
    \end{aligned}
\end{equation}
  
        
        \begin{equation}\label{eq:mondoub5}
\begin{aligned}
\begin{tikzpicture}[xscale=1.8, yscale=1.5]
\node (tl) at (0,2) {$U_{(A \tens B) \tens C}$};
\node (tr) at (4,2) {$U_{A \tens (B \tens C)}$};
\node (ml) at (0,1) {$U_{A \tens B} \tens U_C$};
\node (mr) at (4,1) {$U_A \tens U_{B \tens C}$};
\node (bl) at (0,0) {$(U_A \tens U_B) \tens U_C$};
\node (br) at (4,0) {$U_A \tens (U_B \tens U_C)$};
\draw[->] (tl) to node[above] {$U_{\alpha_{A,B,C}}$} (tr);
\draw[->] (tl) to node[left]{$\fu$} (ml);
\draw[->] (ml) to node[left]{$\fu \tens \id$} (bl);
\draw[->] (tr) to node[left]{$\fu$} (mr);
\draw[->] (mr) to node[left]{$\id \tens \fu$} (br);
\draw[->] (bl) to node[above] {$\alpha_{U_A, U_B, U_C}$} (br);
\end{tikzpicture}
    \end{aligned}
\end{equation}

\item \label{eq:mon3}The following diagrams commute, expressing that the unit
  isomorphisms for $\ten$ are transformations of double categories.
 \begin{equation}\label{eq:mondoub6}
\begin{aligned}
\begin{tikzpicture}[xscale=2, yscale=1.5]
\node (tl) at (0,2) {$(M \tens U_I) \odot (N \tens U_I) $};
\node (tr) at (4,2) {$(M \odot N) \tens ( U_I \odot U_I)$};
\node (ml) at (0,1) {$M \odot N$};
\node (mr) at (4,1) {$(M \odot N) \tens U_I$};
\draw[->] (tl) to node[above] {$\fx$} (tr);
\draw[->] (tl) to node[left]{$\rho_M \odot \rho_M$} (ml);
\draw[->] (tr) to node[left]{$\id \tens \rho_{U_I}$} (mr);
\draw[->] (mr) to node[above] {$\rho_{M \odot N}$} (ml);
\end{tikzpicture}
    \end{aligned}
\end{equation}
      
      \begin{equation}\label{eq:mondoub7}
\begin{aligned}
\begin{tikzpicture}[yscale=1.5]
\node (tl) at (0,2) {$U_{A\tens I}$};
\node (tr) at (4,2) {$U_A \tens  U_I$};
\node (mr) at (4,1) {$U_A$};
\draw[->] (tl) to node[above]{$\fu$} (tr);
\draw[->] (tl) to node[left]{$U_{\rho_A}$} (mr);
\draw[->] (tr) to node[right]{$\rho_{U_A}$} (mr);
\end{tikzpicture}
    \end{aligned}
\end{equation}

 \begin{equation}\label{eq:mondoub8}
\begin{aligned}
\begin{tikzpicture}[xscale=2, yscale=1.5]
\node (tl) at (0,2) {$(U_I \tens M) \odot ( U_I \tens N) $};
\node (tr) at (4,2) {$( U_I \odot U_I) \tens (M \odot N) $};
\node (ml) at (0,1) {$M \odot N$};
\node (mr) at (4,1) {$U_I \tens (M \odot N) $};
\draw[->] (tl) to node[above] {$\fx$} (tr);
\draw[->] (tl) to node[left]{$\lambda_M \odot \lambda_N$} (ml);
\draw[->] (tr) to node[left]{$\lambda_{U_I} \tens \id$} (mr);
\draw[->] (mr) to node[above] {$\lambda_{M \odot N}$} (ml);
\end{tikzpicture}
    \end{aligned}
\end{equation}
      
      \begin{equation}\label{eq:mondoub9}
\begin{aligned}
\begin{tikzpicture}[yscale=1.5]
\node (tl) at (0,2) {$U_{I\tens A}$};
\node (tr) at (4,2) {$U_I \tens  U_A$};
\node (mr) at (4,1) {$U_A$};
\draw[->] (tl) to node[above]{$\fu$} (tr);
\draw[->] (tl) to node[left]{$U_{\lambda_A}$} (mr);
\draw[->] (tr) to node[right]{$\lambda_{U_A}$} (mr);
\end{tikzpicture}
    \end{aligned}
\end{equation}

  \setcounter{mondbl}{\value{enumi}}
\end{enumerate}
Similarly, a \textbf{braided monoidal double category} is a monoidal double
category with the following additional structure.
\begin{enumerate}\setcounter{enumi}{\value{mondbl}}
\item $\lD_0$ and $\lD_1$ are braided monoidal categories.
\item The functors $S$ and $T$ are strict braided monoidal (i.e.\ they
  preserve the braidings).
\item \label{eq:braid1} The following diagrams commute, expressing that the braiding is
  a transformation of double categories. 
  \begin{equation}\label{eq:brmondoub1}
\begin{aligned}
\begin{tikzpicture}[xscale=2, yscale=1.5]
\node (tl) at (0,2) {$(M_1 \odot M_2) \tens (N_1 \odot N_2)$};
\node (tr) at (4,2) {$(N_1 \odot N_2) \tens (M_1 \odot M_2)$};
\node (ml) at (0,1) {$(M_1 \tens N_1) \odot (M_2 \tens N_2)$};
\node (mr) at (4,1) {$(N_1 \tens M_1) \odot (N_2 \tens M_2)$};
\draw[->] (tl) to node[above] {$\fs$} (tr);
\draw[->] (tl) to node[left]{$\fx$} (ml);
\draw[->] (tr) to node[left]{$\fx$} (mr);
\draw[->] (ml) to node[above] {$\fs \odot \fs$} (mr);
\end{tikzpicture}
    \end{aligned}
\end{equation}
 \begin{equation}\label{eq:brmondoub2}
\begin{aligned}
\begin{tikzpicture}[xscale=2, yscale=1.5]
\node (tl) at (0,2) {$U_A \tens U_B$};
\node (tr) at (4,2) {$U_{A \tens B}$};
\node (ml) at (0,1) {$U_B \tens U_A$};
\node (mr) at (4,1) {$U_{B \tens A}$};
\draw[->] (tl) to node[above] {$\fu$} (tr);
\draw[->] (tl) to node[left]{$\fs$} (ml);
\draw[->] (tr) to node[left]{$U_{\fs}$} (mr);
\draw[->] (ml) to node[above] {$\fu$} (mr);
\end{tikzpicture}
    \end{aligned}
\end{equation}

  \setcounter{mondbl}{\value{enumi}}
\end{enumerate}
Finally, a \textbf{symmetric monoidal double category} is a braided one such that
\begin{enumerate}\setcounter{enumi}{\value{mondbl}}
\item $\lD_0$ and $\lD_1$ are in fact symmetric monoidal.
\end{enumerate}
While there are a fair number of coherence diagrams in this definition, most of
them are fairly small, and in any given case most or all of them are
fairly obvious.  Thus, verifying that a given double category is
(braided or symmetric) monoidal is not a great deal of work.

\begin{eg}
  The examples \lnCob, \lMod, \lProf, and \lSpan\ are all easily seen to be
  symmetric monoidal under the tensor product of rings, disjoint union
  of manifolds, and cartesian product of categories, respectively.
\end{eg}

\begin{rmk}
  In a 2-category with finite products there is additionally the
  notion of a \emph{cartesian object}: one such that the diagonal
  $D\to D\times D$ and projection $D\to 1$ have right adjoints.  Any
  cartesian object is a symmetric pseudomonoid in a canonical way,
  just as any category with finite products is a monoidal category
  with its cartesian product.  Many of the ``cartesian bicategories''
  considered in~\cite{cw:cart-bicats-i,ckww:cartbicats-ii} are in
  fact the loose bicategory of some cartesian object in \cDbl,
  and inherit their monoidal structure in this way.
  Cartesian double categories have recently been further studied by~\cite{aleiferi2018cartesian}.
\end{rmk}

Two further general methods for constructing symmetric monoidal double
categories can be found in~\cite{shulman:frbi}; we will return to them in \S\ref{sec:Alg}.

\begin{rmk}
  The general yoga of internalization says that an $X$ internal to
  $Y$s internal to $Z$s is equivalent to a $Y$ internal to $X$s
  internal to $Z$s, but this is only strictly true when the
  internalizations are all strict.  We have defined a symmetric
  monoidal double category to be a (pseudo) symmetric monoid internal
  to (pseudo) categories internal to categories, but one could also
  consider a (pseudo) category internal to (pseudo) symmetric monoids
  internal to categories, i.e.\ a pseudo internal category in the
  2-category
  $\mathcal{S}\mathit{ym}\mathcal{M}\mathit{on}\mathcal{C}\mathit{at}$
  of symmetric monoidal categories and strong symmetric monoidal
  functors.  This would give \emph{almost} the same definition, except
  that $S$ and $T$ would only be strong monoidal (preserving $\ten$ up
  to isomorphism) rather than strict monoidal.  We prefer our
  definition, since $S$ and $T$ are strict monoidal in almost all
  examples, and keeping track of their constraints would be tedious.
\end{rmk}

Just as every bicategory is equivalent to a strict 2-category, it is
proven in~\cite{gp:double-limits} that every pseudo double category is
equivalent to a strict double category (one in which the associativity
and unit constraints for $\odot$ are identities).  Thus, from now on
we will usually omit to write these constraint isomorphisms (or
equivalently, implicitly strictify our double categories).  We
\emph{will} continue to write the constraint isomorphisms for the
monoidal structure $\ten$, since these are where the whole question
lies.

We now move on to define functors and transformations of monoidal double categories.
Like monoidal double categories themselves, these are also special cases of a notion that makes sense internal to any 2-category with products.

\begin{defn}\label{def:monfunc}
  Let $\D$, $\E$ be (braided/symmetric) monoidal double categories.  A {\bf (braided or symmetric) lax monoidal double functor} $F: \D \rightarrow \E$ is a pseudo double functor $F$, together with transformations $\phi : \otimes \circ (F,F) \rightarrow F \circ \otimes$ and $\phi_u:I_{\E}\rightarrow F \circ I_{\D}$ satisfying the usual axioms for (braided/symmetric) monoidal functors with respect to $\otimes$.
\end{defn}

Unfolding the definitions gives us:

\begin{enumerate}
\item $F_0$ and $F_1$ are (braided/symmetric) monoidal functors.
\item The equalities $F_0 \circ S_\D = S_\E \circ F_1$ and $F_0 \circ T_\D = T_\E \circ F_1$ are strict equalities of monoidal functors.
\item The following diagrams commute, expressing that $\phi$ is a transformation of double categories:

\begin{tikzpicture}[xscale=2, yscale=2]
\node (tl) at (0,2) {$(FN \otimes FL) \odot (FM \otimes FK)$};
\node (tr) at (4,2) {$F(N\otimes L) \odot F(M \otimes K)$};
\draw[->] (tl) to node[above] {$\phi \odot \phi$} (tr);
\node (ml) at (0,1) {$(FN \odot FM) \otimes (FL \odot FK)$};
\node (mr) at (4,1) {$F((N \otimes L) \odot (M \otimes K))$};
\node (bl) at (0,0) {$F(N \odot M) \otimes F(L \odot K)$};
\node (br) at (4,0) {$F((N \odot M)\otimes(L\odot K))$};
\draw[->] (tl) to node[left]{$\xi$} (ml);
\draw[->] (ml) to node[left]{$F_{\odot} \otimes F_{\odot}$} (bl);
\draw[->] (tr) to node[left]{$F_{\odot}$} (mr);
\draw[->] (mr) to node[left]{$F(\xi)$} (br);
\draw[->] (bl) to node[above] {$\phi$} (br);
\end{tikzpicture}

\begin{tikzpicture}[yscale=2]
\node (tl) at (0,2) {$U_{FA \otimes FB}$};
\node (tr) at (4,2) {$U_{F(A \otimes B)} $};
\draw[->] (tl) to node[above] {$U_{\phi}$} (tr);
\node (ml) at (0,1) {$U_{FA} \otimes U_{FB}$};
\node (mr) at (4,1) {$F(U_{A\otimes B}) $};
\node (bl) at (0,0) {$F(U_A) \otimes F(U_B)$};
\node (br) at (4,0) {$F(U_A \otimes U_B)$};
\draw[->] (tl) to node[left]{$u$} (ml);
\draw[->] (ml) to node[left]{$F_u \otimes F_u$} (bl);
\draw[->] (tr) to node[left]{$F_U$} (mr);
\draw[->] (mr) to node[left]{$F \circ u$} (br);
\draw[->] (bl) to node[above] {$\phi$} (br);
\end{tikzpicture}

\end{enumerate}

When the natural transformations are in the opposite direction, the functor is {\bf colax monoidal}, and when they are isomorphisms, the functor is {\bf strong monoidal}.

\begin{defn}\label{Def:monverttrans}
  Let $\D$, $\E$ be monoidal double categories and let $(F, \phi) ,(G,\psi): \D \rightarrow \E$ be monoidal double functors. A \textbf{monoidal transformation} $\alpha: F \rightarrow G$ is a tight transformation such that $\alpha_0$ and $\alpha_1$ are monoidal natural transformations.
  Explicitly, this means (in the lax case) that the following equalities hold:

\begin{equation}
\begin{aligned}
\begin{tikzpicture}
\node (tl) at (0,4) {$FA \otimes FB$};
\node (tr) at (4,4) {$FC \otimes FD$};
\node (ml) at (0,2) {$F(A\otimes B)$};
\node (mr) at (4,2) {$F(C \otimes D)$};
\node (bl) at (0,0) {$G(A \otimes B)$};
\node (br) at (4,0) {$G(C \otimes D)$};
\draw[style=tickarrow] (tl) to node [above] {$FM \otimes FN$} (tr);
\draw[style=tickarrow] (ml) to node [above] {$F(M\otimes N)$} (mr);
\draw[->] (tl) to node [left] {$\phi_{A,B}$} (ml);
\draw[->] (tr) to node [right] {$\phi_{C,D}$}(mr);
\draw[->] (ml) to node [left] {$\alpha_{A\otimes B}$} (bl);
\draw[->] (mr) to node [right] {$\alpha_{C \otimes D}$} (br);
\draw[style=tickarrow] (bl) to node [above] {$G(M \otimes N)$} (br);
\node at (2,3) {$\Downarrow \phi_{M,N}$};
\node at (2,1) {$\Downarrow \alpha_{M \otimes N}$};
\end{tikzpicture}
\end{aligned}
=
\begin{aligned}
\begin{tikzpicture}
\node (tl) at (0,4) {$FA \otimes FB$};
\node (tr) at (4,4) {$FC \otimes FD$};
\node (ml) at (0,2) {$GA \otimes GB$};
\node (mr) at (4,2) {$GC \otimes GD$};
\node (bl) at (0,0) {$G(A \otimes B)$};
\node (br) at (4,0) {$G(C \otimes D)$};
\draw[style=tickarrow] (tl) to node [above] {$FM \otimes FN$} (tr);
\draw[style=tickarrow] (ml) to node [above] {$GM \otimes GN$} (mr);
\draw[->] (tl) to node [left] {$\alpha_A \otimes \alpha_B$} (ml);
\draw[->] (tr) to node [right] {$\alpha_C \otimes \alpha_D$} (mr);
\draw[->] (ml) to node [left] {$\psi_{A,B}$} (bl);
\draw[->] (mr) to node [right] {$\psi_{C,D}$} (br);
\draw[style=tickarrow] (bl) to node [above] {$G(M \otimes N)$} (br);
\node at (2,3) {$\Downarrow \alpha_M \otimes \alpha_N$};
\node at (2,1) {$\Downarrow \psi_{M,N}$};
\end{tikzpicture}
\end{aligned}
\end{equation}

\begin{equation}
\begin{aligned}
\begin{tikzpicture}
\node (tl) at (0,4) {$I_{\mathbb{E}}$};
\node (tr) at (4,4) {$I_{\mathbb{E}}$};
\node (ml) at (0,2) {$F(I_{\mathbb{D}})$};
\node (mr) at (4,2) {$F(I_{\mathbb{D}})$};
\node (bl) at (0,0) {$G(I_{\mathbb{D}})$};
\node (br) at (4,0) {$G(I_{\mathbb{D}})$};
\draw[style=tickarrow] (tl) to node [above] {$U_{I_{\mathbb{E}}}$} (tr);
\draw[style=tickarrow] (ml) to node [above] {$F(U_{I_{\mathbb{D}}})$} (mr);
\draw[->] (tl) to node [left] {$\phi_{u_0}$} (ml);
\draw[->] (tr) to node [right] {$\phi_{u_1}$}(mr);
\draw[->] (ml) to node [left] {$\alpha_{I_{\mathbb{D}}}$} (bl);
\draw[->] (mr) to node [right] {$\alpha_{I_{\mathbb{D}}}$} (br);
\draw[style=tickarrow] (bl) to node [below] {$G(U_{I_{\mathbb{D}}})$} (br);
\node at (2,3) {$\Downarrow \phi_{u_1}$};
\node at (2,1) {$\Downarrow \alpha_{U_{I_{\mathbb{D}}}}$};
\end{tikzpicture}
\end{aligned}
=
\begin{aligned}
\begin{tikzpicture}
\node (tl) at (0,4) {$I_{\mathbb{E}}$};
\node (tr) at (4,4) {$I_{\mathbb{E}}$};
\node (ml) at (0,2) {$G(I_{\mathbb{D}})$};
\node (mr) at (4,2) {$G(I_{\mathbb{D}})$};
\draw[style=tickarrow] (tl) to node [above] {$U_{I_{\mathbb{E}}}$} (tr);
\draw[style=tickarrow] (ml) to node [below] {$G(U_{I_{\mathbb{D}}})$} (mr);
\draw[->] (tl) to node [left] {$\psi_{u}$} (ml);
\draw[->] (tr) to node [right] {$\psi_{u}$}(mr);
\node at (2,3) {$\Downarrow \psi_{u_1}$};
\end{tikzpicture}
\end{aligned}
\end{equation}

A {\bf braided or symmetric monoidal tight transformation} is a monoidal transformation between braided/symmetric monoidal functors.
\end{defn}

We have three strict 2-categories $\cMonDbll, \cMonDblc,\cMonDblp$ of monoidal double categories and lax, colax, or pseudo monoidal functors, respectively.
(More generally, we have three 2-categories of pseudomonoids in any 2-category with finite products.)


\section{Companions and conjoints}
\label{sec:comp-conj}

Suppose that \lD\ is a monoidal double category; when does
$\cH(\lD)$ become a monoidal bicategory?  It clearly has a
unit object $I$, and the pseudo double functor $\ten\maps
\lD\times\lD\to\lD$ clearly induces a functor $\ten\maps
\cH(\lD)\times\cH(\lD)\to\cH(\lD)$.  However, the problem is that the
constraint isomorphisms such as $A\ten (B\ten C)\iso (A\ten B)\ten C$
are \emph{tight} 1-cells, which get discarded when we pass to
$\cH(\lD)$.  Thus, in order for $\cH(\lD)$ to inherit a symmetric
monoidal structure, we must have a way to make tight 1-cells
into loose ones.  Thus is the purpose of the following
definition.

\begin{defn}\label{def:companion}
  Let \lD\ be a double category and $f\maps A\to B$ a tight
  1-cell.  A \textbf{companion} of $f$ is a loose 1-cell
  $\fhat\maps A\hto B$ together with 2-morphisms
  \begin{equation*}
    \begin{array}{c}
      \xymatrix@-.5pc{
        \ar[r]|-@{|}^-{\fhat} \ar[d]_f \ar@{}[dr]|{\Downarrow \epsilon_{\hat{f}} }
        & \ar@{=}[d]\\
        \ar[r]|-@{|}_-{U_B} & }
    \end{array}\quad\text{and}\quad
    \begin{array}{c}
      \xymatrix@-.5pc{
        \ar[r]|-@{|}^-{U_A} \ar@{=}[d] \ar@{}[dr]|{\Downarrow \eta_{\hat{f}}}
        & \ar[d]^f\\
        \ar[r]|-@{|}_-{\fhat} & }
    \end{array}
  \end{equation*}
  such that the following equations hold.
  \begin{align}\label{eq:compeqn}
    \begin{array}{c}
      \xymatrix@-.5pc{
        \ar[r]|-@{|}^-{U_A} \ar@{=}[d] \ar@{}[dr]|{\Downarrow \eta_{\hat{f}}}
        & \ar[d]^f\\
        \ar[r]|-{\fhat} \ar[d]_f \ar@{}[dr]|
        {\Downarrow  \epsilon_{\hat{f}} }
        & \ar@{=}[d]\\
        \ar[r]|-@{|}_-{U_B} & }
    \end{array} &= 
    \begin{array}{c}
      \xymatrix@-.5pc{ \ar[r]|-@{|}^-{U_A} \ar[d]_f
        \ar@{}[dr]|{\Downarrow U_f} &  \ar[d]^f\\
        \ar[r]|-@{|}_-{U_B} & }
    \end{array}
    &
    \begin{array}{c}
      \xymatrix@-.5pc{
        \ar[r]|-@{|}^-{U_A} \ar@{=}[d] \ar@{}[dr]|{ \Downarrow \eta_{\hat{f}}}&
        \ar[r]|-@{|}^{\fhat}\ar[d]|f \ar@{}[dr]|{\Downarrow  \epsilon_{\hat{f}} }
        & \ar@{=}[d]\\
        \ar[r]|-@{|}_-{\fhat} &
        \ar[r]|-@{|}_-{U_B} &}
    \end{array} &=
    \begin{array}{c}
      \xymatrix@-.5pc{
        \ar[r]|-@{|}^-{\fhat} \ar@{=}[d] \ar@{}[dr]|{\Downarrow 1_{\fhat}}
        & \ar@{=}[d]\\
        \ar[r]|-@{|}_-{\fhat} & }
    \end{array}
  \end{align}
  A \textbf{conjoint} of $f$, denoted $\fchk\maps B\hto A$, is a
  companion of $f$ in the ``loose opposite'' double category $\lD\lop$.
\end{defn}

\begin{rmk}
  We momentarily suspend our convention of pretending that our double
  categories are strict to mention that the second
  equation in~\eqref{eq:compeqn} actually requires an insertion of unit
  isomorphisms to make sense.
\end{rmk}

The form of this definition is due
to~\cite{gp:double-adjoints,dpp:spans}, but the ideas date back
to~\cite{bs:dblgpd-xedmod}; see
also~\cite{bm:dbl-thin-conn,fiore:pscat}.  In the terminology of these
references, a \emph{connection} on a double category is equivalent to
a strictly functorial choice of a companion for each tight arrow.


\begin{egs}
  \lMod, \lnCob, and \lProf\ have companions and conjoints for all tight morphisms.  In \lMod, the companion
  of a ring homomorphism $f\maps A\to B$ is $B$ regarded as an
  $A$-$B$-bimodule via $f$ on the left, and dually for its conjoint.
  In \lnCob, companions and conjoints are obtained by regarding a
  diffeomorphism as a cobordism.  And in \lProf, companions and
  conjoints are obtained by regarding a functor $f\maps A\to B$ as a
  `representable' profunctor $B(f-,-)$ or $B(-,f-)$.
\end{egs}


The existence of companions and conjoints gives us a way to `lift'
tight 1-cells to loose ones.  What is even more crucial
for our application, however, is that these liftings are unique up to
isomorphism, and that these isomorphisms are canonical and coherent.
This is the content of the following lemmas.  We state most of them
only for companions, but all have dual versions for conjoints.


\begin{lem}\label{thm:theta}
  Let $\fhat\maps A\hto B$ and $\fhat'\maps A\hto B$ be companions of
  $f$ (that is, each comes \emph{equipped with} 2-morphisms as in
  \autoref{def:companion}).  Then there is a unique globular isomorphism
  $\theta_{\fhat,\fhat'}\maps \fhat\too[\iso]\fhat'$ such that
  \begin{equation}\label{eq:comp-iso}
    \vcenter{\xymatrix@R=1.5pc{
        \ar[r]|-@{|}^-{U_A} \ar@{=}[d] \ar@{}[dr]|{\Downarrow \eta_{\hat{f}}} &  \ar[d]^f\\
        \ar[r]|-{\fhat} \ar@{=}[d] \ar@{}[dr]|{\Downarrow \theta_{\fhat,\fhat'}} &  \ar@{=}[d]\\
        \ar[r]|-{\fhat'} \ar[d]_f \ar@{}[dr]|{\Downarrow \epsilon_{\hat{f}'}} & \ar@{=}[d]\\
        \ar[r]|-@{|}_-{U_B} & }} \quad = \quad
    \vcenter{\xymatrix@-.5pc{ \ar[r]|-@{|}^-{U_A} \ar[d]_f
        \ar@{}[dr]|{\Downarrow U_f} &  \ar[d]^f\\
        \ar[r]|-@{|}_-{U_B} & .}}
  \end{equation}
\end{lem}
\begin{proof}
  Composing~\eqref{eq:comp-iso} on the left with
  $\vcenter{\xymatrix@-.5pc{ \ar[r]|-@{|}^-{U_A} \ar@{=}[d]
      \ar@{}[dr]|{\Downarrow \eta_{\hat{f}}} & \ar[d]^f\\ \ar[r]|-@{|}_-{\fhat'} & }}$
  and on the right with $\vcenter{\xymatrix@-.5pc{
      \ar[r]|-@{|}^-{\fhat} \ar[d]_f \ar@{}[dr]|{\Downarrow \epsilon_{\hat{f}}}&
      \ar@{=}[d]\\ \ar[r]|-@{|}_-{U_B} & }}$, and using the second
  equation~\eqref{eq:compeqn}, we see that if~\eqref{eq:comp-iso} is
  satisfied then $\theta_{\fhat,\fhat'}$ must be the composite
  \begin{equation}
    \vcenter{\xymatrix@-.5pc{
        \ar[r]|-@{|}^-{U_A} \ar@{=}[d] \ar@{}[dr]|{\Downarrow \eta_{\hat{f}'}}&
        \ar[r]|-@{|}^-{\fhat} \ar[d]|f \ar@{}[dr]|{\Downarrow \epsilon_{\hat{f}}}
        & \ar@{=}[d]\\
        \ar[r]|-@{|}_-{\fhat'} &
        \ar[r]|-@{|}_-{U_B} &}}\label{eq:theta}
  \end{equation}
  Two applications of the first equation~\eqref{eq:compeqn} shows that
  this indeed satisfies~\eqref{eq:comp-iso}.  As for its being an
  isomorphism, we have the dual composite $\theta_{\fhat',\fhat}$:
  \[\vcenter{\xymatrix@-.5pc{
      \ar[r]|-@{|}^-{U_A} \ar@{=}[d] \ar@{}[dr]|{\Downarrow \eta_{\hat{f}}} &
      \ar[r]|-@{|}^{\fhat'}\ar[d]|f \ar@{}[dr]|
{\Downarrow \epsilon_{\hat{f}'}}
      & \ar@{=}[d]\\
      \ar[r]|-@{|}_-{\fhat} &
      \ar[r]|-@{|}_-{U_B} &}}\]
  which we verify is an inverse using~\eqref{eq:compeqn}:
  \[\vcenter{\xymatrix@-.5pc{
      \ar[r]|-@{|}^{U_A}\ar@{=}[d] \ar@{}[dr]|{=} &
      \ar[r]|-@{|}^{U_A}\ar@{=}[d] \ar@{}[dr]|{\Downarrow \eta_{\hat{f}'}} &
      \ar[r]|-@{|}^{\fhat}\ar[d]|f \ar@{}[dr]|{\Downarrow \epsilon_{\hat{f}}} &
      \ar@{=}[d]\\
      \ar[r]|{U_A}\ar@{=}[d] \ar@{}[dr]|{\Downarrow \eta_{\hat{f}}} &
      \ar[r]|{\fhat'}\ar[d]|f \ar@{}[dr]|{\Downarrow \epsilon_{\hat{f}'}} &
      \ar[r]|{U_B}\ar@{=}[d] \ar@{}[dr]|{=} &
      \ar@{=}[d]\\
      \ar[r]|-@{|}_{\fhat} &
      \ar[r]|-@{|}_{U_B} &
      \ar[r]|-@{|}_{U_B} &
    }} \;=\;
  \vcenter{\xymatrix@-.5pc{
      \ar[r]|-@{|}^-{U_A} \ar@{=}[d] \ar@{}[dr]|{\Downarrow  \eta_{\hat{f}}}&
      \ar[r]|-@{|}^{\fhat}\ar[d]|f \ar@{}[dr]|{\Downarrow \epsilon_{\hat{f}}}
      & \ar@{=}[d]\\
      \ar[r]|-@{|}_-{\fhat} &
      \ar[r]|-@{|}_-{U_B} &}} \;=\;
  \vcenter{\xymatrix@-.5pc{
      \ar[r]|-@{|}^-{\fhat} \ar@{=}[d] \ar@{}[dr]|{\Downarrow 1_{\fhat}}
      & \ar@{=}[d]\\
      \ar[r]|-@{|}_-{\fhat} & }}\]
  (and dually).
\end{proof}

\begin{lem}\label{thm:theta-id}
  For any companion \fhat\ of $f$ we have $\theta_{\fhat,\fhat}=1_{\fhat}$.
\end{lem}
\begin{proof}
  This is the second equation~\eqref{eq:compeqn}.
\end{proof}

\begin{lem}\label{thm:theta-compose-vert}
  Suppose that $f$ has three companions $\fhat$, $\fhat'$, and
  $\fhat''$.  Then $\theta_{\fhat,\fhat''} = \theta_{\fhat',\fhat''}
  \circ\theta_{\fhat,\fhat'}$.
\end{lem}
\begin{proof}
  By definition, we have
  \[\theta_{\fhat',\fhat''} \circ\theta_{\fhat,\fhat'} =\;
  \vcenter{\xymatrix@-.5pc{
      \ar[r]|-@{|}^{U_A}\ar@{=}[d] \ar@{}[dr]|{=} &
      \ar[r]|-@{|}^{U_A}\ar@{=}[d] \ar@{}[dr]|{\Downarrow \eta_{\hat{f}'}} &
      \ar[r]|-@{|}^{\fhat}\ar[d]|f \ar@{}[dr]|{\Downarrow \epsilon_{\hat{f}}} &
      \ar@{=}[d]\\
      \ar[r]|{U_A}\ar@{=}[d] \ar@{}[dr]|{\Downarrow \eta_{\hat{f}''}} &
      \ar[r]|{\fhat'}\ar[d]|f \ar@{}[dr]|{\Downarrow \epsilon_{\hat{f}'}} &
      \ar[r]|{U_B}\ar@{=}[d] \ar@{}[dr]|{=} &
      \ar@{=}[d]\\
      \ar[r]|-@{|}_{\fhat''} &
      \ar[r]|-@{|}_{U_B} &
      \ar[r]|-@{|}_{U_B} &
    }} \;=\;
  \vcenter{\xymatrix@-.5pc{
      \ar[r]|-@{|}^-{U_A} \ar@{=}[d] \ar@{}[dr]|{\Downarrow \eta_{\hat{f}''}}&
      \ar[r]|-@{|}^{\fhat}\ar[d]|f \ar@{}[dr]|{\Downarrow \epsilon_{\hat{f}}}
      & \ar@{=}[d]\\
      \ar[r]|-@{|}_-{\fhat''} &
      \ar[r]|-@{|}_-{U_B} &}} \;=
  \theta_{\fhat,\fhat''}\]
  as desired.
\end{proof}

\begin{lem}\label{thm:comp-unit}
  $U_A\maps A\hto A$ is always a companion of $1_A\maps A\to A$ in a
  canonical way.
\end{lem}
\begin{proof}
  We take both defining 2-morphisms to be
  $1_{U_A}$; the truth of~\eqref{eq:compeqn} is evident.
\end{proof}

\begin{lem}\label{thm:comp-compose}
  Suppose that $f\maps A\to B$ has a companion \fhat\ and $g\maps B\to
  C$ has a companion \ghat.  Then $\ghat\odot\fhat$ is a companion of
  $gf$.
\end{lem}
\begin{proof}
  We take the defining 2-morphisms to be the composites
  \[\vcenter{\xymatrix@-.5pc{
      \ar[r]|-@{|}^-{\fhat} \ar[d]_f \ar@{}[dr]|{\Downarrow \epsilon_{\hat{f}}}&
      \ar[r]|-@{|}^-{\ghat} \ar@{=}[d] \ar@{}[dr]|{1_{\ghat}} &
      \ar@{=}[d]\\
      \ar[r]|-{U_B} \ar[d]_g \ar@{}[dr]|{U_g} &
      \ar[r]|-{\ghat} \ar[d]|g \ar@{}[dr]|{\Downarrow \epsilon_{\hat{g}}}&
      \ar@{=}[d]\\
      \ar[r]|-@{|}_-{U_C} &
      \ar[r]|-@{|}_-{U_C} &
    }}\quad\text{and}\quad
  \vcenter{\xymatrix@-.5pc{
      \ar[r]|-@{|}^-{U_A} \ar@{=}[d] \ar@{}[dr]|{\Downarrow \eta_{\hat{f}}} &
      \ar[r]|-@{|}^-{U_A} \ar[d]|f \ar@{}[dr]|{U_f} &
      \ar[d]^f\\
      \ar[r]|-{\fhat} \ar@{=}[d] \ar@{}[dr]|{1_{\fhat}} &
      \ar[r]|-{U_B} \ar@{=}[d] \ar@{}[dr]|{\Downarrow \eta_{\hat{g}}}&
      \ar[d]^g\\
      \ar[r]|-@{|}_-{\fhat} &
      \ar[r]|-@{|}_-{\ghat} &
    }}
  \]
  It is easy to verify that these satisfy~\eqref{eq:compeqn}, using
  the interchange law for $\odot$ and $\circ$ in a double category.
\end{proof}

\begin{lem}\label{thm:theta-compose-horiz}
  Suppose that $f\maps A\to B$ has companions $\fhat$ and $\fhat'$,
  and that $g\maps B\to C$ has companions $\ghat$ and $\ghat'$.  Then
  $\theta_{\ghat,\ghat'}\odot \theta_{\fhat,\fhat'}  =
    \theta_{\ghat\odot\fhat, \ghat'\odot\fhat'}$.
\end{lem}
\begin{proof}
  Using the interchange law for $\odot$ and $\circ$, we have:
  \begin{align}
    \theta_{\ghat\odot\fhat, \ghat'\odot\fhat'} &=\;
    \vcenter{\xymatrix@-.5pc{
        \ar[r]|-@{|}^-{U_A} \ar@{=}[d] \ar@{}[dr]|{\Downarrow \eta_{\hat{f}'}}&
        \ar[r]|-@{|}^-{U_A} \ar[d]|f \ar@{}[dr]|{U_f} &
        \ar[r]|-@{|}^-{\fhat} \ar[d]|f \ar@{}[dr]|{\Downarrow \epsilon_{\hat{f}}}&
        \ar[r]|-@{|}^-{\ghat} \ar@{=}[d] \ar@{}[dr]|{1_{\fhat}} &
        \ar@{=}[d]\\
        \ar[r]|-{\fhat'} \ar@{=}[d] \ar@{}[dr]|{1_{\ghat}} &
        \ar[r]|-{U_B} \ar@{=}[d] \ar@{}[dr]|{\Downarrow \eta_{\hat{g}'}} &
        \ar[r]|-{U_B} \ar[d]|g \ar@{}[dr]|{U_g} &
        \ar[r]|-{\ghat} \ar[d]|g \ar@{}[dr]|{\Downarrow \epsilon_{\hat{g}}} &
        \ar@{=}[d]\\
        \ar[r]|-@{|}_-{\fhat'} &
        \ar[r]|-@{|}_-{\ghat'} &
        \ar[r]|-@{|}_-{U_C} &
        \ar[r]|-@{|}_-{U_C} &
      }}
    \;=\;
    \vcenter{\xymatrix@-.5pc{
        \ar[r]|-@{|}^-{U_A} \ar@{=}[d] \ar@{}[dr]|{\Downarrow \eta_{\hat{f}'}}&
        \ar[r]|-@{|}^-{\fhat} \ar[d]|f \ar@{}[dr]|{\Downarrow \epsilon_{\hat{f}}}&
        \ar[r]|-@{|}^-{\ghat} \ar@{=}[d] \ar@{}[dr]|{1_{\fhat}} &
        \ar@{=}[d]\\
        \ar[r]|-{\fhat'} \ar@{=}[d] \ar@{}[dr]|{1_{\ghat}} &
        \ar[r]|-{U_B} \ar@{=}[d] \ar@{}[dr]|{\Downarrow \eta_{\hat{g}'}} &
        \ar[r]|-{\ghat} \ar[d]|g \ar@{}[dr]|{\Downarrow \epsilon_{\hat{g}}} &
        \ar@{=}[d]\\
        \ar[r]|-@{|}_-{\fhat'} &
        \ar[r]|-@{|}_-{\ghat'} &
        \ar[r]|-@{|}_-{U_C} &
      }}\\
    &=\;
    \vcenter{\xymatrix@-.5pc{
        \ar[r]|-@{|}^-{U_A} \ar@{=}[d] \ar@{}[dr]|{\Downarrow \eta_{\hat{f}'}}&
        \ar[r]|-@{|}^-{\fhat} \ar[d]|f \ar@{}[dr]|{\Downarrow \epsilon_{\hat{f}}}&
        \ar[r]|-@{|}^-{U_B} \ar@{=}[d] \ar@{}[dr]|{1_{U_B}} &
        \ar[r]|-@{|}^-{\ghat} \ar@{=}[d] \ar@{}[dr]|{1_{\fhat}} &
        \ar@{=}[d]\\
        \ar[r]|-{\fhat'} \ar@{=}[d] \ar@{}[dr]|{1_{\ghat}} &
        \ar[r]|-{U_B} \ar@{=}[d] \ar@{}[dr]|{1_{U_B}} &
        \ar[r]|-{U_B} \ar@{=}[d] \ar@{}[dr]|{\Downarrow \eta_{\hat{g}'}} &
        \ar[r]|-{\ghat} \ar[d]|g \ar@{}[dr]|{\Downarrow \epsilon_{\hat{g}}} &
        \ar@{=}[d]\\
        \ar[r]|-@{|}_-{\fhat'} &
        \ar[r]|-@{|}_-{U_B} &
        \ar[r]|-@{|}_-{\ghat'} &
        \ar[r]|-@{|}_-{U_C} &
      }}\;=\;
    \vcenter{\xymatrix@-.5pc{
        \ar[r]|-@{|}^-{U_A} \ar@{=}[d] \ar@{}[dr]|{\Downarrow \eta_{\hat{f}'}}&
        \ar[r]|-@{|}^-{\fhat} \ar[d]|f \ar@{}[dr]|{\Downarrow \epsilon_{\hat{f}}}&
        \ar[r]|-@{|}^-{U_B} \ar@{=}[d] \ar@{}[dr]|{\Downarrow \eta_{\hat{g}'}}&
        \ar[r]|-@{|}^-{\ghat} \ar[d]|g \ar@{}[dr]|{\Downarrow \epsilon_{\hat{g}}}& \ar@{=}[d]\\
        \ar[r]|-@{|}_-{\fhat'} &
        \ar[r]|-@{|}_-{U_B} &
        \ar[r]|-@{|}_-{\ghat'} &
        \ar[r]|-@{|}_-{U_C} &
      }}\\
    &=\;
    \theta_{\ghat,\ghat'}\odot \theta_{\fhat,\fhat'} 
  \end{align}
  as desired.
\end{proof}

\begin{lem}\label{thm:theta-unit}
  If $f\maps A\to B$ has a companion \fhat, then
  $\theta_{\fhat,\fhat\odot U_A}$ and $\theta_{\fhat,U_B\odot \fhat}$
  are equal to the unit constraints $\fhat \iso \fhat\odot U_A$ and
  $\fhat\iso U_B\odot \fhat$.
\end{lem}
\begin{proof}
  By definition, we have
  \[\theta_{\fhat,\fhat\odot U_A} =\;
  \vcenter{\xymatrix@-.5pc{
      \ar[r]|-@{|}^-{U_A} \ar@{=}[d] \ar@{}[dr]|{\Downarrow 1_{U_A}} &
      \ar[r]|-@{|}^-{U_A} \ar@{=}[d] \ar@{}[dr]|{1_{U_A}} &
      \ar@{=}[d] \ar[rr]|-@{|}^-{\fhat} \ar@{}[ddrr]|{\Downarrow \epsilon_{\hat{f}}}&& \ar@{=}[dd]\\
      \ar[r]|-{U_A} \ar@{=}[d] \ar@{}[dr]|{1_{U_A}} &
      \ar[r]|-{U_A} \ar@{=}[d] \ar@{}[dr]|{\Downarrow \eta_{\hat{f}}}&
      \ar[d]^f\\
      \ar[r]|-@{|}_-{U_A} &
      \ar[r]|-@{|}_-{\fhat} & \ar[rr]|-@{|}^-{U_B} &&
    }}\;=\;
  \vcenter{\xymatrix{ \ar[r]|-@{|}^-{U_A} \ar@{=}[d]
      \ar@{}[dr]|{\Downarrow 1_{U_A}} &  \ar@{=}[d]\\
      \ar[r]|-@{|}_-{U_A} & }}
  \]
  which, bearing in mind our suppression of unit and associativity
  constraints, means that in actuality it is the unit constraint
  $\fhat \iso \fhat\odot U_A$.  The other case is dual.
\end{proof}

\begin{lem}\label{thm:comp-func}
  Let $F\maps \lD\to\lE$ be a functor between double categories and
  let $f\maps A\to B$ have a companion \fhat\ in \lD.  Then $F(\fhat)$
  is a companion of $F(f)$ in \lE.
\end{lem}
\begin{proof}
  We take the defining 2-morphisms to be
  \[\vcenter{\xymatrix@R=1.5pc@C=3pc{
      \ar[r]|-@{|}^-{F(\fhat)} \ar[d]_{F(f)}
      \ar@{}[dr]|{F(\Downarrow \epsilon_{\hat{f}})} &  \ar@{=}[d]\\
      \ar[r]|-{F(U_B)} \ar@{=}[d] \ar@{}[dr]|\iso &  \ar@{=}[d]\\
      \ar[r]|-@{|}_-{U_{F(B)}} & }}
  \quad\text{and}\quad
  \vcenter{\xymatrix@R=1.5pc@C=3pc{
      \ar[r]|-@{|}^-{U_{FA}} \ar@{=}[d] \ar@{}[dr]|\iso & \ar@{=}[d]\\
      \ar[r]|-{F(U_{A})} \ar@{=}[d] \ar@{}[dr]|{F(\Downarrow \eta_{\hat{f}})} & 
      \ar[d]^{F(f)}\\
      \ar[r]|-@{|}_-{F(\fhat)} & .}}\]
  The axioms~\eqref{eq:compeqn} follow directly from those for \fhat.
\end{proof}


\begin{lem}\label{thm:theta-func}
  Suppose that $F\maps \lD\to\lE$ is a functor and that $f\maps A\to
  B$ has companions \fhat\ and $\fhat'$ in \lD.  Then
  $\theta_{F(\fhat),F(\fhat')} = F(\theta_{\fhat,\fhat'})$.
\end{lem}
\begin{proof}
  Using the axioms of a pseudo double functor and the definition of
  the 2-morphisms in \autoref{thm:comp-func}, we have
  \begin{equation}
    F(\theta_{\fhat,\fhat'})
    =\;
    \vcenter{\xymatrix@C=3pc{
        \ar[r]|-@{|}^-{F(\fhat)}
        \ar[d] \ar@{}[dr]|{\Downarrow F(\eta_{\hat{f}'} \odot\epsilon_{\hat{f}})} &  \ar[d]\\
        \ar[r]|-@{|}_-{F(\fhat')} &}}
    \;=\;
    \vcenter{\xymatrix@C=3pc{
        \ar[rr]|-@{|}^-{F(\fhat)}
        \ar@{=}[d] \ar@{}[drr]|\iso &&  \ar@{=}[d]\\
        \ar[r]|-@{|}^-{F(U_{A})} \ar@{=}[d]
        \ar@{}[dr]|{\Downarrow F(\eta_{\hat{f}'})} &
        \ar[r]|-@{|}^-{F(\fhat)} \ar[d]|{F(f)}
        \ar@{}[dr]|{\Downarrow F(\epsilon_{\hat{f}})}
        & \ar@{=}[d]\\
        \ar[r]|-@{|}_-{F(\fhat')} \ar@{}[drr]|\iso\ar@{=}[d] &
        \ar[r]|-@{|}_-{U_{F(B)}} & \ar@{=}[d]\\
        \ar[rr]|-@{|}_-{F(\fhat')} && }}
    \;=\;
    \vcenter{\xymatrix@R=1.5pc@C=3pc{
        \ar[r]|-@{|}^-{U_{F(A)}} \ar@{=}[d] \ar@{}[dr]|\iso &
        \ar[r]|-@{|}^-{F(\fhat)} \ar@{=}[d] \ar@{}[dr]|=
        & \ar@{=}[d]\\
        \ar[r]|-{F(U_{A})} \ar@{=}[d] \ar@{}[dr]|{\Downarrow F(\eta_{\hat{f}'})} &
        \ar[r]|-{F(\fhat)} \ar[d]|{F(f)} \ar@{}[dr]|{\Downarrow F(\epsilon_{\hat{f}})}
        & \ar@{=}[d]\\
        \ar[r]|-{F(\fhat')}  \ar@{=}[d] \ar@{}[dr]|= &
        \ar[r]|-{F(U_{B})} \ar@{}[dr]|\iso  \ar@{=}[d] & \ar@{=}[d]\\
        \ar[r]|-@{|}_-{F(\fhat')} &
        \ar[r]|-@{|}_-{U_{F(B)}} &}}
    \;=
    \theta_{F(\fhat),\,F(\fhat')}
  \end{equation}
  as desired.
\end{proof}

\begin{lem}\label{thm:comp-iso}
  If $f\maps A\to B$ is a tight isomorphism with a companion \fhat,
  then \fhat\ is a conjoint of its inverse $f\inv$.
\end{lem}
\begin{proof}
  The composites
  \[\vcenter{\xymatrix@-.5pc{
      \ar[r]|-@{|}^{\fhat}\ar[d]_f \ar@{}[dr]|{\Downarrow} &
      \ar@{=}[d]\\
      \ar[r]|{U_B}\ar[d]_{f\inv} \ar@{}[dr]|{\Downarrow U_{f\inv}} &
      \ar[d]^{f\inv}\\
      \ar[r]|-@{|}_{U_A} &
    }}\quad\text{and}\quad
  \vcenter{\xymatrix@-.5pc{
      \ar[r]|-@{|}^{U_B}\ar[d]_{f\inv} \ar@{}[dr]|{\Downarrow U_{f\inv}} &
      \ar[d]^{f\inv}\\
      \ar[r]|{U_A}\ar@{=}[d] \ar@{}[dr]|{\Downarrow} &
      \ar[d]^f\\
      \ar[r]|-@{|}_{\fhat} &
    }}
  \]
  exhibit \fhat\ as a conjoint of $f\inv$.
\end{proof}

\begin{lem}\label{thm:compconj-adj}
  If $f\maps A\to B$ has both a companion \fhat\ and a conjoint \fchk,
  then we have an adjunction $\fhat\adj\fchk$ in $\cH\lD$.  If $f$ is
  an isomorphism, then this is an adjoint equivalence.
\end{lem}
\begin{proof}
  The unit and counit of the adjunction $\fhat\adj\fchk$ are the
  composites
  \[\vcenter{\xymatrix@-.5pc{
      \ar[r]|-@{|}^{U_A}\ar@{=}[d] \ar@{}[dr]|{\Downarrow \eta_{\hat{f}}} &
      \ar[r]|-@{|}^{U_A}\ar[d]|{f} \ar@{}[dr]|{\Downarrow \eta_{\check{f}}} &
      \ar@{=}[d]\\
      \ar[r]|-@{|}_{\fhat} &
      \ar[r]|-@{|}_{\fchk} &
    }}\quad\text{and}\quad
  \vcenter{\xymatrix@-.5pc{
      \ar[r]|-@{|}^{\fchk}\ar@{=}[d] \ar@{}[dr]|{\Downarrow \epsilon_{\check{f}}} &
      \ar[r]|-@{|}^{\fhat}\ar[d]|{f} \ar@{}[dr]|{\Downarrow \epsilon_{\hat{f}}} &
      \ar@{=}[d]\\
      \ar[r]|-@{|}_{U_B} &
      \ar[r]|-@{|}_{U_B} &
    }}
  \]
  The triangle identities follow from~\eqref{eq:compeqn}.  If $f$ is
  an isomorphism, then by the dual of \autoref{thm:comp-iso}, \fchk\
  is a companion of $f\inv$.  But then by \autoref{thm:comp-compose}
  $\fchk\odot \fhat$ is a companion of $1_A=f\inv \circ f$ and
  $\fhat\odot\fchk$ is a companion of $1_B = f\circ f\inv$, and hence
  \fhat\ and \fchk\ are equivalences.  We can then check that in this
  case the above unit and counit actually are the isomorphisms
  $\theta$, or appeal to the general fact that any adjunction
  involving an equivalence is an adjoint equivalence.
\end{proof}


\begin{lem}\label{lem:FUtheta}
Suppose $F:\lD \rightarrow \lE$ is a functor of double categories. The 2-cell $F_U$ is equal to $\theta_{\id_{FA}, F\id_A}$.
\end{lem}

\begin{proof}
We show that equation~\eqref{eq:comp-iso} holds when we substitute $\theta_{\id_{FA}, F\id_A}$ by $\hat{F}_U$.  Unfolding the definitions of $\eta_{U_{FA}}$, $\epsilon_{FU_A}$, and $\hat{F}_U$, and applying functoriality of $F$, we obtain an expression that can be rewritten to $U_{\id_{FA}}$. It follows that $F_U$ is a $\theta$-isomorphism, by the uniqueness of $\theta s$ in this expression.
\end{proof}

\begin{rmk}
  It is tempting to want to state a general coherence theorem along the lines of ``any two composites of $\theta$-isomorphisms having the same source and target are equal.''
  However, like statements such as ``any two composites of constraints in a monoidal category are equal'', this statement is actually literally false, because to determine a $\theta$-isomorphism requires not only a source and target but also the choice of companion data.
  If a given 1-cell is a companion of the same 1-cell in more than one way (which is the case as soon as it has any nontrivial automorphisms), then there will be different $\theta$-isomorphisms with the same source and target.
  This is analogous to how in a particular monoidal category there can be ``accidental'' composites of constraints that are not covered by the coherence theorem.
  It is probably possible to state a general coherence theorem for $\theta$-isomorphisms that is sufficiently careful to be true, but we will not need this.
\end{rmk}

\begin{rmk}
  Since all the tight constraints of a monoidal double category are invertible, to construct its underlying monoidal bicategory we only need it to have companions (and hence, by \cref{thm:comp-iso}, conjoints) for all tight \emph{isomorphisms}.
  In~\cite{gg:ldstr-tricat} double categories of this sort were called \textbf{fibrant}; one might also say \textbf{isofibrant} for emphasis.
  Note that this condition is equivalent to asking that the (source, target) functor $\lD_1 \to \lD_0\times\lD_0$ is an isofibration (i.e.\ has the isomorphism-lifting property).
  
  To lift lax or colax monoidal \emph{functors}, and noninvertible transformations between monoidal functors, to the bicategorical level, we require our double categories to have companions (or conjoints, depending on the directions) for noninvertible tight morphisms as well.
  In~\cite{shulman:frbi} double categories with companions and conjoints for \emph{all} tight morphisms were called \emph{framed bicategories}; this is equivalent to asking that the (source, target) functor $\lD_1 \to \lD_0\times\lD_0$ be a Grothendieck fibration or opfibration (either assumption implies the other).
  In \cref{sec:1x1-to-bicat} we will need a further condition to ensure that the ``componentwise'' companion of a tight transformation is pseudo, rather than lax or colax, natural.
\end{rmk}


\section{From double categories to bicategories}
\label{sec:1x1-to-bicat}

We are now equipped to lift structures on double categories to
their loose bicategories.  In this section we show that passage
from double categories to bicategories is given by a functor of locally cubical bicategories. In order to prove this, we first give an intermediate result that $\cH$ lifts to a functor of hom-bicategories
\begin{align}
    \cDbl(\lD,\lE) &\too \cBicat(\cH(\lD),\cH(\lE))
\end{align}

As a point of notation, we write $\odot$ for the composition of
1-cells in a bicategory, since our bicategories are generally of the
form $\cH(\lD)$.  As advocated by Max Kelly, we say \textbf{functor}
to mean a morphism between bicategories that preserves composition up
to coherent isomorphism; equivalent terms include \emph{weak 2-functor},
\emph{pseudofunctor}, and \emph{homomorphism}.
We will not discuss lax functors (a.k.a. ``morphisms'') between bicategories at all in this paper.

Recall that the assignment $\cH$ sends each double category $\lC$ to the loose bicategory  $\cH(\lD)$ of objects, 1-cells, and globular 2-morphisms of $\lD$.  Note that functors of double categories and bicategories compose strictly associatively; hence, we can talk about the 1-categories of double categories and bicategories, which we denote ${\bf Dbl}$ and ${\bf Bicat}$ respectively.

\begin{thm}\label{thm:1-func}
 If \lD\ is a double category, then $\cH(\lD)$ is a bicategory, and
  any functor $F\maps \lD\to\lE$ induces a functor $\cH(F)\maps
  \cH(\lD)\to\cH(\lE)$.  In this way $\cH$ defines a functor of
  1-categories $\mathbf{Dbl}\to \mathbf{Bicat}$.
\end{thm}
\begin{proof}
 The constraints of $F$ are all globular, hence give constraints for
  $\cH(F)$.  Functoriality is evident.
\end{proof}

The action of \cH\ on transformations is less obvious. It
requires the presence of companions or conjoints to lift the part of the data given by tight morphisms to loose 1-cells. Before we discuss how this works, we briefly recall some definitions regarding transformations between functors of bicategories.

If $F,G\maps \cA\to\cB$ are functors between bicategories, then an
\textbf{oplax transformation} $\al\maps F\to G$ consists of 1-cells
$\al_A\maps FA\to GA$ and 2-cells
\[\vcenter{\xymatrix{ \ar[r]^{Ff}\ar[d]_{\al_A} \drtwocell\omit{\al_f} &  \ar[d]^{\al_B}\\
  \ar[r]_{Gf} & }}\]
such that for any 2-cell $\xymatrix{A \rtwocell^f_g{x} & B}$ in \cA,
\begin{equation}
  \label{eq:laxtransf-nat}
  \vcenter{\xymatrix@R=1pc@C=3pc{
      \rtwocell^{Ff}_{Fg}{Fx}\ar[dd]_{\al_A} 
      &  \ar[dd]^{\al_B}\\
      \drtwocell\omit{\al_g} & \\
      \ar[r]_{Gg} & }}\;=\;
  \vcenter{\xymatrix@R=1pc@C=3pc{
      \ar[r]^{Ff}\ar[dd]_{\al_A} \drtwocell\omit{\al_f} &
      \ar[dd]^{\al_B}\\ & \\
      \rtwocell^{Gf}_{Gg}{Gx} & }}
\end{equation}
and moreover for any $A$ and any $f,g$ in \cA,
\begin{equation}
  \vcenter{\xymatrix@R=5pc{
      \rtwocell^{1_{FA}}_{F(1_A)}{\iso} \ar[d]_{\al_A} \drtwocell\omit{\al_{1_A}} &  \ar[d]^{\al_A}\\
      \rtwocell^{G(1_A)}_{1_{GA}}{\iso} & }} \;=\;
  \vcenter{\xymatrix{ \ar[r]^{1_{FA}}\ar[d]_{\al_A} \drtwocell\omit{\iso}&  \ar[d]^{\al_A}\\
      \ar[r]_{1_{GA}} &
    }}
  \quad\text{and}\quad
  \vcenter{\xymatrix{
      \ar[r]|{Ff}\ar[d]_{\al_A} \drtwocell\omit{\al_f}
      \rruppertwocell^{F(gf)}{\iso}
      &
      \ar[r]|{Fg}\ar[d]|{\al_B} \drtwocell\omit{\al_g} &
      \ar[d]^{\al_C}\\
      \ar[r]|{Gf} \rrlowertwocell_{G(gf)}{\iso} & \ar[r]|{Gg} & }}
  \;=\;
  \vcenter{\xymatrix{ \ar[r]^{F(gf)}\ar[d]_{\al_A} \drtwocell\omit{\al_{gf}} &  \ar[d]^{\al_C}\\
      \ar[r]_{G(gf)} & }}\label{eq:laxtransf-ax}
\end{equation}
It is a \textbf{lax transformation} if the 2-cells $\al_f$ go the
other direction, and a \textbf{pseudo transformation} if they are
isomorphisms.

When two functors of bicategories agree on objects, there is a simpler notion of transformation between them, called an \emph{icon}.
An icon is, morally speaking, an oplax transformation whose 1-cell components are all identities; but as noted by~\cite{lack:icons} this can be reexpressed without referring to these identity morphisms at all, yielding a definition that is easier to work with (because identity 1-cells in a bicategory are not strict).

\begin{defn}
Let $\cD, \cE$ be bicategories, and let $F,G: \cD \rightarrow \cE$ be functors that agree on objects. An \textbf{icon} $\alpha: F \Rightarrow G$ is given by a family of 2-cells $\alpha_f : Ff \Rightarrow Gf$ indexed by the 1-cells of $\cD$, which are natural in $f$ and such that for all objects $A, B, C$ and 1-cells $A \xrightarrow{f} B \xrightarrow{g} C$ the following equations hold:

\begin{equation}\label{eq:iconeq}
\begin{aligned}
\begin{tikzpicture}[yscale=2,xscale=3]
\node (tl) at (0,1) {$I_{FA}$};
\node (tr) at (1,1) {$F I_A$};
\node (bl) at (0,0) {$I_{GA}$};
\node (br) at (1,0) {$GI_A$};
\draw[->] (tl) to node[above]{$\iso$} (tr);
\draw[->] (bl) to node[below]{$\iso$} (br);
\draw[->] (tl) to node[left]{$\iso$} (bl);
\draw[->] (tr) to node[right]{$\alpha_{I_A}$} (br);
\end{tikzpicture}
\end{aligned}
\qquad
\begin{aligned}\begin{tikzpicture}[yscale=2, xscale=3]
\node (tl) at (0,1) {$F(g) \odot F(f)$};
\node (tr) at (1,1) {$F(g \odot f)$};
\node (bl) at (0,0) {$G(g) \odot G(f)$};
\node (br) at (1,0) {$G(g \odot f)$};
\draw[->] (tl) to node[above]{$\iso$} (tr);
\draw[->] (bl) to node[below]{$\iso$} (br);
\draw[->] (tl) to node[left]{$\alpha_g \odot \alpha_f$}(bl);
\draw[->] (tr) to node[right]{$\alpha_{g \odot f}$} (br);
\end{tikzpicture}
\end{aligned}
\end{equation}
\end{defn}

Recall also that if $\al,\al'\maps F\to G$ are oplax transformations,
a \textbf{modification} $\mu\maps \al\to\al'$ consists of 2-cells
$\mu_A\maps \al_A\to\al'_A$ such that
\begin{equation}
  \vcenter{\xymatrix@C=1pc@R=2.5pc{ \ar[rr]^{Ff}\dtwocell_{\al'_A}^{\al_A}{\mu_A}  &
      \drtwocell\omit{\al_f} &  \ar[d]^{\al_B}\\
      \ar[rr]_{Gf} && }} \quad=\quad
  \vcenter{\xymatrix@C=1pc@R=2.5pc{ \ar[rr]^{Ff}\ar[d]_{\al'_A} \drtwocell\omit{\al'_f} && 
      \dtwocell^{\al_B}_{\al'_B}{\mu_B}\\
      \ar[rr]_{Gf} && }}\label{eq:modif-ax}
\end{equation}
There is an evident notion of modification between lax transformations
as well.
We have three bicategories
\[ \cBicat_c(\cA,\cB) \qquad \cBicat_l(\cA,\cB) \qquad \cBicat_p(\cA,\cB) \]
whose objects are the functors $\cA\to\cB$, whose morphisms are colax, lax, and pseudo transformations respectively, and whose 2-morphisms are modifications.

A \textbf{pseudo natural adjoint equivalence} is, by definition, an internal adjoint equivalence in $\cBicat_p(\cA,\cB)$.
However, by doctrinal adjunction~\cite{kelly:doc-adjn}, an internal adjoint equivalence in $\cBicat_c(\cA,\cB)$ or $\cBicat_l(\cA,\cB)$ is automatically pseudo natural as well.

Let $\cDblcf$ denote the sub-2-category of $\cDbl$ containing all double categories and all functors between them, but only the tight transformations $\al:F\to G: \lD\to\lE$ such that each tight component $\al_A$ has a loose companion $\widehat{\al_A}$.
Note that if \lE\ is (iso)fibrant, every invertible $\al$ has this property, and if \lE\ has companions for all tight 1-morphisms then every transformation has this property.

\begin{thm}\label{thm:h-locfr}
  We have a functor of bicategories
  \begin{align}
    \cDblcf(\lD,\lE) &\too \cBicat_c(\cH(\lD),\cH(\lE))\\
    F &\mapsto \cH(F)\\
    \al &\mapsto \alhat.
  \end{align}
\end{thm}

Note that we are here regarding the 1-category $\cDblcf(\lD,\lE)$ as a bicategory with only identity 2-cells.
Since any functor of bicategories preserves adjoint equivalences, and an adjoint equivalence in a mere category is simply an isomorphism, it follows that if $\al$ is an isomorphism then $\alhat$ is (equipped with the structure of) a pseudo natural adjoint equivalence.

\begin{proof}
  We define the 1-cell components of $\alhat$ by choosing companions $\alhat_A = \widehat{\al_A}$ of each component of $\al$.
  The 2-cell component $\alhat_f$ is the composite
  \begin{equation}
    \vcenter{\xymatrix@R=1.5pc@C=2.5pc{
        \ar[r]|-@{|}^{U_{FA}}\ar@{=}[d] \ar@{}[dr]|{\Downarrow \eta_{\hat{\alpha}_A}} &
        \ar[r]^{Ff}\ar[d]|{\al_A} \ar@{}[dr]|{\Downarrow \al_f} &
        \ar[r]|-@{|}^{\alhat_B}\ar[d]|{\al_B} \ar@{}[dr]|{\Downarrow \epsilon_{\hat{\alpha}_B}} &
        \ar@{=}[d]\\
        \ar[r]|-@{|}_{\alhat_A} &
        \ar[r]_{Gf} &
        \ar[r]|-@{|}_{U_{GB}} & 
      }}\label{eq:oplax-2cell}
  \end{equation}
  Equations~\eqref{eq:laxtransf-nat} and~\eqref{eq:laxtransf-ax}
  follow directly from \autoref{thm:dbl-transf}.

  It is left to construct the constraints and check the axioms for functors of bicategories. Suppose we are given $\al\maps F\to G$ and $\be\maps G\to H$.  Then by
  \autoref{thm:comp-compose}, $\behat_A\odot\alhat_A$ is a companion
  of $\be_A\circ \al_A$, so we have a canonical isomorphism given by the icon
  \[\theta_{\widehat{\be\al}_A, \,\behat_A\odot\alhat_A}\maps
  \widehat{\be\al}_A \too[\iso] \behat_A\odot\alhat_A.
  \]
  Of course, we also have $\theta_{\widehat{1_A},U_A}\maps
  \widehat{1_A} \too[\iso] U_A$ by \autoref{thm:comp-unit}.  These
  constraints are automatically natural, since $\cDbl(\lD,\lE)$ has no
  nonidentity 2-cells.  The axiom for the composition constraint says
  that two constructed isomorphisms of the form
  \[\widehat{\gm\be\al}_A \too[\iso] (\gmhat_A \odot \behat_A)\odot \alhat_A\]
  are equal.  However, both $\widehat{\gm\be\al}_A$ and $(\gmhat_A
  \odot \behat_A)\odot \alhat_A$ are companions of $\gm_A\be_A\al_A$,
  and both of these isomorphisms are constructed from composites of $\theta$s;
  hence they are both equal to
  \[\theta_{\widehat{\gm\be\al}_A,\, (\gmhat_A \odot \behat_A)\odot
    \alhat_A}\] and thus equal to each other.  The same argument
  applies to the axioms for the unit constraint; thus we have a functor of bicategories.
\end{proof}

Of course, the functor constructed in \cref{thm:h-locfr} depends on the choices of companions made in the proof.
However up to equivalence it does not depend on these choices:

\begin{lem}\label{thm:h-locfr-uniq}
  Suppose we make two different sets of choices of companions for each component of a tight transformation in $\cDblcf(\lD,\lE)$, yielding by the proof of \cref{thm:h-locfr} two different functors
  \[\cH,\cH'\maps \cDblcf(\lD,\lE)\too \cBicat_c(\cH(\lD),\cH(\lE)).\]
  Then the isomorphisms $\theta$ from \autoref{thm:theta} fit together
  into an invertible icon $\cH\cong \cH'$.
\end{lem}
\begin{proof}
  We must first show that for a given transformation $\al\maps F\to
  G\maps \lD\to\lE$ in \cDbl, the isomorphisms \th\ correspond to 2-cells $\alhat \iso \alhat'$ of $\cBicat_c(\cH(\lD),\cH(\lE))$; that is, they form invertible
  modifications.
  Substituting~\eqref{eq:oplax-2cell} and the definition of \th\
  into~\eqref{eq:modif-ax}, this becomes the assertion that
  \begin{equation}
    \vcenter{\xymatrix@R=1.5pc@C=2pc{
        &
        \ar[r]|-@{|}^{U_{FA}}\ar@{=}[d] \ar@{}[dr]|{\Downarrow \eta_{\hat{\alpha}_A}} &
        \ar[r]^{Ff}\ar[d]|{\al_A} \ar@{}[dr]|{\Downarrow \al_f} &
        \ar[r]|-@{|}^{\alhat_B}\ar[d]|{\al_B} \ar@{}[dr]|{\Downarrow \epsilon_{\hat{\alpha}_B}} &
        \ar@{=}[d]\\
        \ar[r]|-@{|}^{U_{FA}} \ar@{=}[d] \ar@{}[dr]|{\Downarrow \eta_{\hat{\alpha}_A'}} &
        \ar[r]|{\alhat_A} \ar[d]|{\al_A} \ar@{}[dr]|{\Downarrow \epsilon_{\hat{\alpha}_A}}&
        \ar[r]_{Gf}  \ar@{=}[d] &
        \ar[r]|-@{|}_{U_{GB}} & \\
        \ar[r]|-@{|}_{\alhat_A'} & \ar[r]|-@{|}_{U_{GB}}&&
      }} \;=\;
    \vcenter{\xymatrix@R=1.5pc@C=2pc{
        && \ar@{=}[d] \ar[r]|-@{|}^{U_{FA}} \ar@{}[dr]|{\Downarrow \eta_{\hat{\alpha}_B'}} &
        \ar[d]|{\al_B} \ar[r]|-@{|}^{\alhat_B} \ar@{}[dr]|{\Downarrow \epsilon_{\hat{\alpha}_B}}
        &
        \ar@{=}[d] &\\
        \ar[r]|-@{|}^{U_{FA}}\ar@{=}[d] \ar@{}[dr]|{\Downarrow \eta_{\hat{\alpha}_A'}} &
        \ar[r]^{Ff}\ar[d]|{\al_A} \ar@{}[dr]|{\Downarrow \al_f} &
        \ar[r]|{\alhat_B'}\ar[d]|{\al_B} \ar@{}[dr]|{\Downarrow \epsilon_{\hat{\alpha}_B'}} &
        \ar@{=}[d] \ar[r]|-@{|}_{U_{GB}}&\\
        \ar[r]|-@{|}_{\alhat_A'} &
        \ar[r]_{Gf} &
        \ar[r]|-@{|}_{U_{GB}} & .
      }}
  \end{equation}
  This follows from two applications of~\eqref{eq:compeqn}, one for
  $\alhat_A$ and one for $\alhat_B'$.
  Now we show that these form an invertible icon. The compatibility axiom with 2-cells is vacuous since $\cDblcf(\lD,\lE)$ has
  no nonidentity 2-cells, while the functoriality requirement follows from uniqueness of $\theta$s,
  since all the constraints involved are also instances of \th.
\end{proof}

\begin{rmk}
  There is a similar result for conjoints instead of companions, but hedged about with duality.
  Recall that a conjoint in $\lE$ is the same as a companion in the ``loose opposite'' $\lE\lop$, and more generally the 2-functor $(-)\lop:\cDbl \to \cDbl$ takes transformations with conjoints to transformations with companions.
  Thus if we denote by $\cDbllf(\lD,\lE)$ the category of functors and transformations having componentwise conjoints, then $(-)\lop$ defines a functor $\cDbllf(\lD,lE) \to \cDblcf(\lD\lop,\lE\lop)$.

  We also have $\cH(\lE\lop) = (\cH\lE)\op$, where $(-)\op$ denots the bicategory obtained by reversing 1-cells but not 2-cells.
  Since $(-)\op$ reverses the direction of transformations, and interchanges lax and colax transformations, we have a functor $\cBicat_c(\cA,\cB) \to \cBicat_l(\cA\op,\cB\op)\op$.
  Thus, composing all of this up, we obtain a functor
  \begin{multline*}
    \cDbllf(\lD,\lE) \to \cDblcf(\lD\lop,\lE\lop) \to \cBicat_c(\cH(\lD\lop),\cH(\lE\lop))\\ = \cBicat_c((\cH\lD)\op, (\cH\lE)\op) \to \cBicat_l(\cH\lD, \cH\lE)\op.
  \end{multline*}
  In other words, a transformation with componentwise conjoints induces a \emph{lax} transformation in the \emph{opposite} direction between loose bicategories.
  If a transformation has both componentwise companions and conjoints (such as if $\lE$ has all companions and conjoints), then the resulting colax and lax transformations are doctrinally adjoint; in~\cite{shulman:smbicat} we called such a pair a ``conjunctional transformation''.
\end{rmk}

So far we have been able to deal with arbitrary tight transformations and their resulting colax (or lax) transformations.
But when we come to talk about tricategories or locally cubical bicategories, we have to restrict to pseudo natural transformations on the side of bicategories, since there is no tricategory (or locally cubical bicategory) containing arbitrary colax transformations: the interchange law only holds laxly.
(One could write down a notion of ``colax tricategory'', but we will not attempt this.)
This restriction means we need a corresponding restriction on the double-categorical side.

\begin{defn}
  A transformation $\alpha$ in $\cDbl$ has \textbf{loosely strong companions} if each component $\alpha_A$ has a loose companion, and the resulting colax transformation $\hat\alpha$ is actually pseudo natural.
  We write $\cDblf$ for the 2-category of double categories, functors, and transformations with loosely strong companions.
\end{defn}

\begin{rmk}\label{thm:iso-strong}
  If $\alpha$ is an isomorphism and its codomain has companions, then $\alpha$ automatically has loosely strong companions.
  For the companions of $\alpha^{-1}$ are inverse adjoint equivalences of those of $\alpha$, making $\alpha$ into a pseudo natural adjoint equivalence.
\end{rmk}

We are now ready to prove that $\cH$ lifts to a functor of locally cubical bicategories. The definition of a locally cubical bicategory can be found in~\cite{gg:ldstr-tricat}, where it is called a locally cubical bicategory. It can also be derived as a bicategory enriched in the monoidal 2-category $\cDbl$, in the following sense.

\begin{defn}\label{def:lcbc}
Let $\mathcal{V}$ be a monoidal 2-category. A \textbf{$\mathcal{V}$-enriched bicategory} $\fB$ consists of the following data.
\begin{enumerate}
\item A collection of objects $\fB=\{A,B,C,...\} $
\item For each two objects $A,B$ of $\fB$, an object $\fB(A,B)$ of $\mathcal{V}$. 
\item For each three objects $A, B, C$ of $\fB$, a 1-cell  $\comp: \fB(B,C) \otimes \fB(A,B) \rightarrow \fB(A,C)$ of $\mathcal{V}$. 
\item For each object $A$ of $\fB$, a 1-cell \newline $I_A: 1 \rightarrow \fB(A,A)$ of $\mathcal{V}$, where $1$ is the initial object of $\mathcal{V}$. 
\item For each four objects $A,B,C,D$ of $\fB$, a 2-cell of $\mathcal{V}$, depicted below. 
  \begin{center}
    \begin{tikzpicture}[yscale=1.5, xscale=2.5]
\node (A) at (0,1){$\fB(C,D) \times \fB(B, C) \times \fB(A,B)$};
\node (C) at (2,1){$\fB(C,D) \times \fB(A,C)$};
\node (D) at (0,0) {$\fB(B,D)  \times \fB(A,B)$};
\node (E) at (2,0) {$\fB(A, D) $};
\node (B) at (1,.5) {$\Downarrow \alpha$};
\draw[->] (A) to node[above]{$\otimes \times \id$} (C);
\draw[->] (A) to node[left]{$\id \times \otimes$} (D);
\draw[->] (C) to node[right]{$\otimes$} (E);
\draw[->] (D) to node[above]{$\otimes$} (E);
\end{tikzpicture}
    \end{center}
\item For each two objects $A,B$ of $\fB$, a natural isomorphisms 
  \begin{center}
    \begin{equation*}
\begin{aligned}
\begin{tikzpicture}[yscale=2.5, xscale=3]
\node (A) at (0,1){$\fB(A,B) \times 1$};
\node (B) at (1,0) {$\fB(A,B)$};
\node (D) at (0,0) {$\fB(B,B) \times \fB(A,B)$};
\node (C) at (.3,.3){$\Downarrow \lambda$};
\draw[->] (A) to node[above]{$\iso$} (B);
\draw[->] (D) to node[below]{$\otimes$} (B);
\draw[->] (A) to node[left]{$  I_B \times \id$} (D);
\end{tikzpicture}
\hspace{1cm}
\begin{tikzpicture}[yscale=2.5, xscale=3]
\node (A) at (0,1){$1 \times \fB(A,B)$};
\node (B) at (1,0) {$\fB(A,B)$};
\node (D) at (0,0) {$\fB(A,B) \times \fB(A,A)$};
\node (C) at (.3,.3){$\Downarrow\rho$};
\draw[->] (A) to node[above]{$\iso$} (B);
\draw[->] (D) to node[below]{$\otimes$} (B);
\draw[->] (A) to node[left]{$\id \times I_B$} (D);
\end{tikzpicture}
\end{aligned}
\end{equation*}
    \end{center}
  \item The usual axioms for bicategories~\cite{maclane} hold.
 \end{enumerate}
\end{defn}

The 2-category $\cDbl$ has finite products, which gives the monoidal structure. 
When we apply the definition above to $\cDbl$, we see that a locally cubical bicategory has objects, hom-double categories {of 1-cells, tight 2-cells, loose 2-cells, and 3-cells which fit in a square of tight and loose 2-cells. We have double functors $\comp$, $I_A$ and tight isomorphisms $\alpha, \rho, \lambda$.
  
\begin{rmk}\label{rmk:lcbop}
  Recall from \autoref{rmk:dblop} that a double category has three opposites $\lD\lop$, $\lD\ttop$, and $\lD\tlop$ depending on whether we reverse the loose 1-cells, tight 1-cells, or both.
  Applied homwise, these operations induce three duals of a locally cubical bicategory, which we denote $\fB\lco$, $\fB\tco$, and $\fB\tlco$, and which respectively reverse the loose 2-cells, tight 2-cells, and both.
  (If we added the possibility to reverse the 1-cells we would obtain seven different ``opposites'' such as $\fB\op$, $\fB\lcoop$, and so on, but we will have no use for that duality in this paper.)
\end{rmk}

Note that a locally cubical bicategory with one object is precisely a monoidal double category in the sense of \cref{sec:symm-mono-double}, and thus the explicit description of the coherences in \cref{sec:symm-mono-double} can also be applied here.
We will primarily be concerned with the following two examples.

\begin{eg}
  Any category can be regarded as a double category in the \emph{loose} direction, with only identity tight morphisms and identity 2-cells.
  (Note that this is a \emph{strict} double category, not a pseudo one.)
  This operation preserves products, and thereby any strict 2-category can be regarded as a locally cubical bicategory.
  In particular, we will regard the 2-category $\cDblf$ as a locally cubical bicategory in this way; thus it has only identity tight 2-cells and identity 3-cells. We write $\fDblf$ for this locally cubical bicategory.
\end{eg}

\begin{eg}
  By~\cite[Corollary 12]{gg:ldstr-tricat}, there is a locally cubical bicategory $\fBicat$ defined by the following data:
  \begin{itemize}
  \item Its objects are bicategories.
  \item Its morphisms are functors.
  \item Its loose 2-cells are pseudo natural transformations.
  \item Its tight 2-cells are icons~\cite{lack:icons}.
  \item Its 3-cells are \emph{cubical modifications}, defined below:
  \end{itemize}
\end{eg}

\begin{defn}[{\cite[Definition 13]{gg:ldstr-tricat}}]
Let $F,G,H,K: \cD \rightarrow \cE$ be pseudo functors; let $\alpha: F \Rightarrow G$, $\beta: H \Rightarrow K$ be pseudo transformations; let $\gamma: F \Rightarrow H$, $\delta: G \Rightarrow K$ be icons. A \textbf{cubical modification}
\[
\begin{tikzpicture}[scale=1.5]
\node (tl) at (0,1) {$F$};
\node (tr) at (1,1) {$G$};
\node (bl) at (0,0) {$H$};
\node (br) at (1,0) {$K$};
\draw[doubletight] (tl) to node[above]{$\alpha$} (tr);
\draw[doubletight] (bl) to node[below]{$\beta$} (br);
\draw[doubletight] (tl) to node[left]{$\gamma$} (bl);
\draw[doubletight] (tr) to node[right]{$\delta$} (br);
\node at (.5,.5) {$\DDownarrow \Gamma$};
\end{tikzpicture}
\]
is given by a family of 2-cells $\Gamma_A: \alpha_A \RRightarrow \beta_A$ such that for every 1-cell $f:A \rightarrow B$ of $\cD$, the following equality holds.

 \begin{equation}
 \begin{aligned}
 \begin{tikzpicture}[scale=1.5]
 \node (tl) at (-1,1) {$FA$};
 \node (tm) at (0,1) {$FB$};
 \node (tr) at (1,1) {$GB$};
 \node (bl) at (-1,0) {$FA$};
 \node (bm) at (0,0) {$GA$};
 \node (br) at (01,0) {$GB$};
 \node (bl1) at (-1,-.7){$HA$};  
 \node (bm1) at (0,-.7) {$KA$};
 \node (br1) at (1,-.7) {$KB$}; 
 \draw[doubletight] (tm)  to node[above]{$\alpha_B$} (tr);
 \draw[doubleeq] (bm) to (bm1);
 \draw[doubletight] (bm) to node[above] {$Gf$}(br);
 \draw[doubleeq] (tr) to (br);
 \draw[doubleeq] (tl)  to  (tm);
 \draw[doubleeq] (tl) to (bl);
 \draw[doubletight] (tl) to node[above]{$Ff$}(tm);
 \draw[doubletight] (bl) to node[above]{$\alpha_A$}(bm);
 \node at (0,.5) {\footnotesize $\Downarrow \alpha_f$}; 
 \node at (0.5,-.3) {\footnotesize $\Downarrow \delta_f$}; 
  \node at (-0.5,-.3) {\footnotesize $\Downarrow \Gamma_A$};
 \draw[doubletight] (bl1)  to node[above]{$\beta_A$} (bm1);
 \draw[doubletight] (bm1) to  node[above]{$Kf$}(br1);
 \draw[doubleeq] (bl)  to (bl1);
 \draw[doubleeq] (br)  to (br1);
 \end{tikzpicture}
 \end{aligned}
 =
\begin{aligned}
 \begin{tikzpicture}[scale=1.5]
 \node (tl) at (-1,1) {$FA$};
 \node (tm) at (0,1) {$FB$};
 \node (tr) at (1,1) {$GB$};
 \node (bl) at (-1,0) {$HA$};
 \node (bm) at (0,0) {$HB$};
 \node (br) at (01,0) {$KB$};
 \node (bl1) at (-1,-.7){$HA$};  
 \node (bm1) at (0,-.7) {$KA$};
 \node (br1) at (1,-.7) {$KB$}; 
 \draw[doubletight] (tm)  to node[above]{$\alpha_B$} (tr);
 \draw[doubleeq] (tm) to (bm);
 \draw[doubletight] (bm) to node[above] {$\beta_B$}(br);
 \draw[doubleeq] (tr) to (br);
 \draw[doubleeq] (tl)  to  (tm);
 \draw[doubleeq] (tl) to (bl);
 \draw[doubletight] (tl) to node[above]{$Ff$}(tm);
 \draw[doubletight] (bl) to node[above]{$Hf$}(bm);
 \node at (-0.5,.5) {\footnotesize $\Downarrow \gamma_f$}; 
 \node at (0.5,.5) {\footnotesize $\Downarrow \Gamma_B$}; 
 \draw[doubletight] (bl1)  to node[above]{$\beta_A$} (bm1);
 \draw[doubletight] (bm1) to  node[above]{$Kf$}(br1);
 \draw[doubleeq] (bl)  to (bl1);
 \draw[doubleeq] (br)  to (br1);
 \node at (0,-0.3) {\footnotesize $\DDownarrow \beta_f$}; 
 \end{tikzpicture}
 \end{aligned}
\end{equation}

\end{defn}

The double functor $\transid_A: * \rightarrow \fBicat(A,A)$ maps the cells in the trivial bicategory $*$ to the identity cells and morphisms of $\fBicat(A,A)$. 
The double functor $\comp$ is defined on functors of bicategories by composition. On pseudo transformations and icons it is given by the Godement product. On cubical modifications it is defined below:

\begin{equation*}
\begin{aligned}
 \begin{tikzpicture}[scale=2]
 \node (tl) at (-1,1) {$FF'A$};
 \node (tm) at (0,1) {$GF'A$};
 \node (tr) at (1,1) {$GG'A$};
 \node (bl) at (-1,0) {$HF'A$};
 \node (bm) at (0,0) {$KF'A$};
 \node (br) at (01,0) {$KG'A$};
 \node (bl1) at (-1,-1){$HH'A$};  
 \node (bm1) at (0,-1) {$KH'A$};
 \node (br1) at (1,-1) {$KK'A$}; 
 \draw[doubletight] (tm)  to node[above]{$G(\alpha'_A)$} (tr);
 \draw[doubleeq] (tm) to (bm);
 \draw[doubletight] (bm) to node[above] {$K(\alpha'_A)$}(br);
 \draw[doubleeq] (tr) to (br);
 \draw[doubleeq] (tl)  to  (tm);
 \draw[doubleeq] (tl) to (bl);
  \draw[doubleeq] (bm) to (bm1);
 \draw[doubletight] (tl) to node[above]{$\alpha_{F'A}$}(tm);
 \draw[doubletight] (bl) to node[above]{$\beta_{F'A}$}(bm);
 \node at (-0.5,.5) {\footnotesize $\Downarrow \Gamma_{F'A}$}; 
 \node at (0.5,.5) {\footnotesize $\Downarrow \delta_{\alpha'_A}$}; 
 \draw[doubletight] (bl1)  to node[above]{$\beta_{H'A}$} (bm1);
 \draw[doubletight] (bm1) to  node[above]{$K(\beta'A)$}(br1);
 \draw[doubleeq] (bl)  to (bl1);
 \draw[doubleeq] (br)  to (br1);
 \node at (-.5,-0.5) {\footnotesize $=$}; 
\node at (.5,-0.5) {\footnotesize $\DDownarrow K\Gamma'_A$}; 
\end{tikzpicture}
\end{aligned}
\end{equation*}

Functoriality follows from naturality of the icons. Note that there are several equivalent ways to define this composition on cubical modifications, by choosing different versions of the Godement product.  

Now, we can similarly obtain the notion of \emph{functor} between locally cubical bicategories from the definition of a $\mathcal{V}$-enriched functor between $\mathcal{V}$-enriched bicategories.

\begin{defn}\label{def:lcbcfunc}
Let $\mathcal{V}$ be a monoidal 2-category. Let ${\fB,\fC}$ be $\mathcal{V}$-enriched bicategories. A \textbf{$\mathcal{V}$-enriched functor} $F: {\fB} \rightarrow {\fC}$ consists of the following data:
\begin{enumerate}
\item An assignment on objects that sends each object $A$ of ${\fB}$ to an object $F A$ of ${\fC}$.
\item For each two objects $A,B$ of ${\fB}$, a 1-cell ${\fB}(A,B) \rightarrow {\fC}(F(A),F(B))$ of $\mathcal{V}$.
\item For every triple of objects $A,B,C$ of ${\fB}$, a 2-cell of $\mathcal{V}$ 
\begin{align} 
\begin{tikzpicture}
\node(1) at (0,0) {${\fB}(A,B) \otimes {\fB}(B,C)$};
\node(2) at (5,0) {${\fC}(F(A),F(C))$};
\draw[->] (1) to[in=155, out=25] node[above]{$\fB \comp $} (2); 
\draw[->] (1) to[in=-155, out=-25] node[below]{$ \comp (F,F)$} (2); 
\node at (2.5,0) {$\Downarrow \phi \iso$};
\end{tikzpicture}
\end{align}
\item For every object $A$ of ${\fB}$ a 2-cell of $\mathcal{V}$
\begin{align}
\begin{tikzpicture}[xscale=.5, yscale=.3]
\node(1) at (0,0) {$*$};
\node(2) at (5,0) {${\fC}(A,A)$};
\node(3) at (5,-5) {${\fB}(F(A),F(A))$};
\draw[->] (1) to node[above]{$\looseid_{A}$} (2); 
\draw[->] (1) to node[below]{$\looseid_{F(A)}$} (3);
\draw[->] (2) to node[right]{$F$} (3); 
\node at (3.5,-1.5) {$\Downarrow \phi_u \iso$};
\end{tikzpicture}
\end{align}
\item The usual coherence diagrams, Definition 10 of~\cite{nick:tricatsbook} commute.
\end{enumerate}
\end{defn}

When we apply this definition to the monoidal 2-category $\cDbl$, we see that a {\bf functor of locally cubical bicategories} $F$ consists of a map of objects $A \mapsto F A$; pseudo double functors ${\fB}(A,B) \rightarrow {\fC}(F(A),F(B))$ for each two objects $A,B$; and tight transformations $\phi$ for each two objects $A,B$ of ${\fB}$, and $\phi_u$ for every object $A$ of ${\fB}$, plus axioms.

\begin{thm}\label{thm:h-functor}
There is a functor of locally cubical bicategories $\cL \maps \fDblf\to \fBicat$.
\end{thm}
\begin{proof}
Note that the 1-functor $\cH:\mathbf{Dbl}\to\mathbf{Bicat}$ has a left adjoint that regards a bicategory as a double category with only identity tight 1-cells.
Thus, the functor of bicategories $\cDblf(\D,\E) \rightarrow \cBicat(\cH(\D), \cH (\E)) = \cH( \fBicat(\cH(\D),\cH(\E)))$ has an adjunct pseudo double functor $\fDblf(\D,\E) \rightarrow \fBicat(\cH\D, \cH \E)$. Consequently, the first two requirements in \autoref{def:lcbcfunc} are satisfied by Theorems \ref{thm:1-func} and \ref{thm:h-locfr}.
Since $\cH$ strictly preserves composition of 1-cells,
the third requirement amounts to the existence of a tight transformation $\phi\maps \behat * \alhat \iso \widehat{\be*\al}$ for every pair of transformations with loosely strong companions 

  \[\vcenter{\xymatrix{\lC \rtwocell^F_G{\al} & \lD \rtwocell^H_K{\be}
      & \lE}}\]
      
      such that 
 \begin{equation}
        \vcenter{\xymatrix@-.5pc{
        1_{{\cH}H \odot {\cH}F} \ar[r]\ar[d]_{=} &
        \hat{1}_{H}* \hat{1}_{F}\ar[d]^{\phi}\\
        1_{{\cH}(H \odot F)}\ar[r] &
        \widehat{1_{H} * 1_{F}}}} \quad\text{and}\quad       
    \vcenter{\xymatrix@-.5pc{
        \widehat{\gm\al}* \widehat{\de\be} \ar[r]\ar[d]_\phi &
        (\gmhat* \dehat)\circ(\alhat* \behat)\ar[d]^{\phi * \phi}\\
        \widehat{\gm\al* \de\be}\ar[r] &
        (\widehat{\gm* \de})\circ(\widehat{\al* \be})}}
  \end{equation}
commute. 
Here, we use the 'Godement product' $*$ of 2-cells in $\cDbl$.  

  Now by Lemmas \ref{thm:comp-compose} and
  \ref{thm:comp-func}, $(\behat *\alhat)_A = \behat_{GA} \circ
  H(\alhat_A)$ is a companion of $(\be*\al)_A = \be_{GA} \circ
  H(\al_A)$.  Therefore, we take the component $(\phi_{\alpha,\beta})_A$ to be
  \[\theta_{\behat_{GA} \circ H(\alhat_A),\, \widehat{\be*\al}_A}.\]
 As the other morphisms in the diagrams above are also $\theta$-isomorphisms, the equations hold by Lemma~\ref{thm:theta-compose-vert}.
For the tight transformation  $\phi_u$ we can simply take the identity, since $\cH$ is strictly unital.
The coherence equations hold by Lemma~\ref{thm:theta-compose-vert}
\end{proof}

Our goal is to enhance this functor to act on ``monoidal objects''.
It is well-known that ``monoidal functors preserve monoid objects'', so our approach will be to categorify this: we will show that the functor $\cH$ is monoidal, in an appropriate sense, and that monoidal functors of this sort preserve monoidal objects of the appropriate sort.

In fact, the monoidality of $\cH$ is easy to describe, because the monoidal structures of $\fDblf$ and $\fBicat$ are cartesian and very strict.

In general, if $\mathcal{V}$ is a monoidal 2-category with strict 2-categorical finite products (such as \cDbl), we say that a $\mathcal{V}$-enriched bicategory ${\fB}$ has \textbf{finite products} when for each two objects $C,D \in {\fB}$ there is an object $C\times D$ with projections $C\times D\to C$ and $C\times D\to D$ (i.e.\ morphisms $I\to \fB(C\times D,C)$ and $I\to \fB(C\times D,D)$ in \cV) inducing an \emph{isomorphism} in $\mathcal{V}$ (not merely an equivalence):
\begin{align}
\fB(A, C \times D) \xrightarrow{\cong} \fB(A,C) \times \fB(A,D)
\end{align}
and similarly there is a strict terminal object $\ast$ such that $\fB(A,\ast)$ is strictly terminal in \cV\ for all $A$.
This holds for \fBicat\ and \fDblf, because cartesian products of bicategories and double categories are simply componentwise, and all the morphisms in \fBicat\ and \fDblf\ (no matter how weak) are defined in terms of data in their targets.

Similarly, we say that a functor $F$ of \cV-enriched bicategories \textbf{preserves products} if it takes the terminal object to a terminal object and pairs of product projections $A \leftarrow A\times B \to B$ to pairs of product projections (in the above strict sense).

\begin{thm}
The functor of locally cubical bicategories $\cH: \fDblf \rightarrow \fBicat$ preserves products.
\end{thm}
\begin{proof}
Since $\cH$ merely forgets a part of the double categories and double functors, we have simple equalities
$\cH(\mathbb{D} \times \mathbb{E}) = \cH(\mathbb{D}) \times \cH(\mathbb{E})$, and the product projections are likewise preserved.
The case of the terminal object is likewise easy.
\end{proof}


\section{Monoidal objects in locally cubical bicategories}
\label{sec:mono-objects}

We now move on to define an appropriate abstract sort of ``monoidal objects'', that will be preserved by the product-preserving functor $\cH$, and that specializes to monoidal double categories and to monoidal bicategories. 
Such structures are obtained by abstracting the definitions of monoidal, braided, and symmetric structure for bicategories given in~\cite{nick:tricatsbook},~\cite{mccrudden:bal-coalgb}, and~\cite{gg:ldstr-tricat}. 
Under this translation pseudonatural transformations become \emph{loose} 2-cells and modifications become globular 3-cells.
Recall that the loose 2-cells in \fDblf\ (which has no nonidentity tight 2-morphisms) are the (tight) transformations, while those in \fBicat\ are exactly the pseudonatural transformations (its tight 2-morphisms are icons).
We will define monoidal, braided, sylleptic and symmetric objects, 1-cell, 2-cells, icons and 3-cells by taking instances of bicategories, functors, pseudonatural transformations, and modifications in the original definitions as objects, 1-cells, 2-cells, and 3-cells, respectively. 

We start by unfolding the structure of a locally cubical bicategory and fixing our notation. In addition to objects, locally cubical bicategories have hom-double-categories whose objects, loose 1-cells, tight 1-cells, and 2-morphisms we will call '1-cells', 'loose 2-cells', 'tight 2-cells' and '3-cells' of the locally cubical bicategory, respectively.  We write ``$\onecell$" to denote 1-cells, ``$\looseRightarrow$" to denote the loose 2-cells, ``$\Rightarrow$" to denote the tight 2-cells and ``$\RRightarrow$" to denote $3$-cells.

There are three types of composition, for which we will adopt the notation introduced for intercategories~\cite{gp:intercategories-i} (of which locally cubical bicategories are a special case). Firstly, we have loose composition ``$\horc$" within the hom-double-categories, giving composition of loose 2-cells along a 1-cell boundary and of 3-cells along a tight 2-cell boundary. We write this in diagrammatic order: $\alpha \horc \beta$ means ``$\beta$ after $\alpha$". We write $\looseid_{f}$ and $\looseid_{\alpha}$ for the loose identity on a 1-cell $f$ and a tight 2-cell $\alpha$, respectively. Loose composition is associative and unital up to isomorphism.  We denote its associator and left and right unitors by $\hora, \horr, \horl$, respectively; they are globular 3-cells.
However, for readability we will implicitly strictify the hom-double-categories, as we did with the double categories in Section~\ref{sec:symm-mono-double}, except where we prove that this structure is preserved by monoidal cells. As a consequence, we omit the associativity and unit constraints for $\horc$ in various places.

Secondly, we have tight composition ``$\verc$" in the hom-double-categories. This gives us composition of tight 2-cells along a 1-cell boundary and tight composition of 3-cells along a loose 2-cell boundary, written in the conventional order: $f \verc g$ denoting ``$f$ after $g$". We write $\tightid_f$ and $\tightid_{\alpha}$ for the tight identity on a 1-cell $f$ and a loose 2-cell $\alpha$, respectively. Tight composition is strictly associative and has strict identities.

Thirdly, there is composition ``$\comp$" of 1-cells, 2-cells, and 3-cells along a 0-cell boundary, given by the enriched structure. We write this composition in the conventional order: $f \comp g$ meaning ``$f$ after $g$". When it is clear from the context, we omit the composition symbol ``$\comp$",  and write the juxtaposition of 1-cells instead. The identities for this composition are denoted by ``$\transid$". This composition is weakly associative and unital in general, but we will restrict our input locally cubical bicategories by assuming that composition of \emph{1-cells} along 0-cells is strict, in the following precise sense:

\begin{defn}\label{defn:1strict}
  A locally cubical bicategory is \textbf{1-strict} if the tight transformations $\alpha$, $\lambda$, and $\rho$ from \cref{def:lcbc} are globular, i.e.\ their tight 1-cell components are identities.
  Similarly, a functor of locally cubical bicategories is \textbf{1-strict} if the tight transformations $\phi$ and $\phi_u$ from \cref{def:lcbcfunc} are globular.
\end{defn}

Our primary examples $\fDblf$ and $\fBicat$ are 1-strict, as is the functor $\cH:\fDblf\to\fBicat$.
However, while $\fDblf$ is completely strict (being a strict 2-category regarded as a locally cubical bicategory in the 2-loose direction), the composition $\comp$ in $\fBicat$ is not fully strict, i.e.\ the transformations $\alpha,\lambda,\rho$ are not identities: their nonidentity 2-cell components arise from the pseudonaturality constraints of the loose 2-cells in $\fBicat$.
In general, composition of loose 2-cells along 0-cells in a 1-strict locally cubical bicategory is only weakly associative and unital, with an associator and left and right unitors that we denote by by $\compa, \compr, \compl$, respectively; they are globular 3-cells.

A 1-strict locally cubical bicategory is almost exactly an abstraction of the structure and level of strictness possessed by bicategories, functors, pseudonatural transformations, and icons.
This enables an almost-literal transposition of definition from the world of bicategories, as long as the monoidal constraints and coherence axioms are reformulated in equational form instead of their original diagrammatic expression of 2-cell composites in the target bicategory.
(Such a reformulation is entirely algorithmic, and could be automated by formulating an ``internal type theory'' of locally cubical bicategories.)

There is one other difference between $\fBicat$ and a general locally cubical bicategory: the latter contains a unique basic operation of ``composing two loose 2-cells along a 0-cell'', whereas when composing two pseudonatural transformations of bicategories
\[ \xymatrix{ \cA \rtwocell^F_G{\alpha} & \cB \rtwocell^H_K{\beta} & \cC } \]
we have to choose whether the component of $\beta\circ\alpha$ at $a\in \cA$ should be the top-right composite or the left-bottom composite in the following square:
\[
  \begin{tikzcd}
    H F a \ar[r,"H \alpha_a"] \ar[d,"\beta_{F a}"'] \ar[dr,phantom,"\cong"] & H G a \ar[d,"\beta_{G a}"]\\
    K F a \ar[r,"K \alpha_a"'] & K G a.
  \end{tikzcd}
\]
The two choices are isomorphic (the square is a pseudonaturality constraint for $\beta$), but neither of them is canonically ``the'' choice.
This ambiguity is not in itself a problem; the real issue it reflects is that the locally cubical bicategory $\fBicat$ does not have a basic way to express the operation of ``whiskering'' a pseudonatural transformation by a functor on either side.
We can compose  along an object with a loose identity 2-cell, but since bicategories are only weakly unital, the result is only \emph{isomorphic} to the ``basic'' operation of whiskering for bicategories.
In particular, bicategories are strict in one way that is invisible to $\fBicat$: whiskering is also strictly associative, e.g.\ $H(F(\gamma)) = (H\circ F)(\gamma)$, but $\fBicat$ can only express the isomorphism $\looseid_{H}\circ (\looseid_{F}\circ \gamma) \cong \looseid_{H \circ F} \circ \gamma$ which is \emph{not} a strict equality.
This means that certain additional isomorphisms have to be inserted in our definitions of monoidal objects that are not present in the explicit definitions of monoidal bicategories.

One solution to this problem would be to define a locally-cubical-type structure analogous to the \emph{bicategory-enriched categories} of~\cite{verity:base-change}, which use a multicategory structure to express whiskering as a basic operation (just as is done in the strict case with Gray-categories~\cite{gps:tricats}).
However, probably a better use of time would be to extend our definitions to locally cubical bicategories that are not even 1-strict, such as the locally cubical bicategory of \emph{monoidal} bicategories and its relatives (i.e.\ the \emph{outputs} of our construction).
This is more subtle than it may appear, since the coherence axioms of a monoidal bicategory involve pasting together loose 2-cells along boundary 1-cells that are strictly equal in a 1-strict locally cubical bicategory, but only related by a \emph{tight} isomorphism in the general case; thus even to write down the definitions the general case would already require some fibrancy in the hom-double-categories.

One further note about our definitions is that certain constraint 2-cells are suppressed in~\cite{nick:tricatsbook},~\cite{mccrudden:bal-coalgb}, and~\cite{gg:ldstr-tricat} to simplify notation, but in the equational form below, this cannot be done.
Thus the coherence equations look more intimidating, but in fact are not more complicated than their bicategorical counterparts.

Let \fB\ be a 1-strict locally cubical bicategory with products.
In the following definitions, we will give the coherence axioms for monoidal objects, 1-cells, 2-cells and 3-cells explicitly (in Appendix \ref{ap:coherence}), but we generally leave the axioms for braided, symmetric and sylleptic monoidal structures as an exercise for the reader.
They can be obtained by translating the standard axioms for the appropriate kinds of monoidal bicategory into the language of locally cubical bicategories, in an entirely analogous way.

\begin{defn}
A {\bf monoidal object} in \fB\ is an object $A$, equipped with 1-cells $\otimes: A \times A \onecell A$ and $I_A: * \onecell A$, and loose 2-cell equivalences
\begin{itemize} 
\item $\alpha: \mathord{\ten} (\mathord{\tens} \times \id) \looseRightarrow{} \mathord{\ten}  (\id \times \mathord{\tens}) $
\item $l: \mathord{\ten} (I \times \transid) i_2 \looseRightarrow{} \transid$ and $r: \mathord{\ten} (\transid \times I) i_1 \looseRightarrow{} \transid$ 
\end{itemize}
where $*$ is the terminal object of $\fB$ and $i_1: A \rightarrow A \times *$, $i_2: A \rightarrow 
* \times A$ are the canonical isomorphisms. Finally, it must be equipped with the invertible globular 3-cells $\pi, \mu, \lambda, \rho$ depicted below, which generalise the Mac Lane pentagon and the three other coherence diagrams given in Definition 4.1 of~\cite{nick:tricatsbook}. The horizontally drawn equality is added for readability. These 3-cells must satisfy the three axioms~\eqref{eq:monobjeq1},~\eqref{eq:monobjeq2}~\eqref{eq:monobjeq3} in Appendix~\ref{ap:coherence}, which generalise the coherence axioms for monoidal bicategories~\cite[Definition 4.1]{nick:tricatsbook}.

\begin{equation}\label{eq:pi}
\begin{aligned}
 \begin{tikzpicture}[xscale=3.4]
\node (t0) at (0,2) {\scriptsize $\tens (\tens \times \transid)(\tens \times \transid \times \transid)$};
\node (t1) at (1,2) {\scriptsize $\tens(\tens \times \transid)(\transid \times
\tens \times \transid)$};
\node (t2) at (2,2) {\scriptsize $\tens(\transid \times \tens)(\transid \times
\tens \times \transid)$};
\node (t3) at (3,2) {\scriptsize $\tens (\transid \times \tens)(\transid \times \transid \times \tens)$};
\node (m0) at (0,1) {\scriptsize $\tens (\tens \times \transid)(\tens \times \transid \times \transid)$};
\node (m1) at (1,1) {\scriptsize $\tens(\transid \times \tens )(\tens \times \transid \times \transid)$};
\node (m2) at (2,1) {\scriptsize $\tens(\tens \times \transid)(\transid \times \transid \times
\tens)$};
\node (m3) at (3,1) {\scriptsize $\tens (\transid \times \tens)(\transid \times \transid \times \tens)$};
\draw[doubletighteq] (t3) to (m3);
\draw[doubletighteq] (t0) to (m0);
\draw[doubleloose] (t0) to node[above]{\scriptsize $\looseid (\alpha \times \looseid)$} (t1);
\draw[doubleloose] (t1) to node[above]{\scriptsize $\alpha \looseid$}
(t2);
\draw[doubleloose] (t2) to node[above]{\scriptsize $\looseid (\looseid \times \alpha)$} (t3);
\draw[doubleloose] (m0) to node[above]{\scriptsize $\alpha \looseid$} (m1);
\draw[doubletighteq] (m1) to (m2);
\draw[doubleloose] (m2) to node[above]{\scriptsize $\alpha \looseid$} (m3);
\node at (1.5,1.5) {\scriptsize $\DDownarrow \pi$};
\end{tikzpicture}
\end{aligned}
\end{equation}

\begin{equation}\label{eq:mu}
\begin{aligned}
{\scriptstyle
\begin{tikzpicture}[xscale=3.5]
\node (t0) at (0,2) {\scriptsize $\tens (\transid \times \transid)$};
\node (t1) at (1,2) {\scriptsize$\tens(\tens \times \transid)(\transid \times
I \times \transid)$};
\node (t2) at (2,2) {\scriptsize$\tens(\transid \times \tens)(\transid \times
I \times \transid)$};
\node (t3) at (3,2) {\scriptsize$\tens (\transid \times \transid)$};
\node (m0) at (0,1) {\scriptsize$\tens$};
\node (m3) at (3,1) {\scriptsize$\tens $};
\draw[doubletighteq] (t3) to (m3);
\draw[doubletighteq] (t0) to (m0);
\draw[doubleloose] (t0) to node[above]{\scriptsize $\looseid (r^{-1} \times \looseid)$} (t1);
\draw[doubleloose] (t1) to node[above]{\scriptsize$\alpha$}
(t2);
\draw[doubleloose] (t2) to node[above]{\scriptsize$\looseid (\looseid \times l)$} (t3);
\draw[doubleloose] (m0) to node[above]{\scriptsize$\looseid$}(m3);
\node at (1.5,1.5) {\scriptsize$\DDownarrow \mu$};
\end{tikzpicture}}
\end{aligned}
\end{equation}

\begin{equation}\label{eq:lcbclambda}
\begin{aligned}
{\scriptstyle \begin{tikzpicture}[xscale=3.5]
\node (t0) at (0,2) {\scriptsize$\tens (\tens \times \transid)(I \times \transid \times \transid)$};
\node (t1) at (1,2) {\scriptsize$\tens(\transid \times \transid)$};
\node (t3) at (3,2) {\scriptsize$\transid \tens$};
\node (m0) at (0,1) {\scriptsize$\tens (\tens \times \transid)(I \times \transid \times \transid)$};
\node (m1) at (1,1) {\scriptsize$\tens(\transid \times \tens )(I \times \transid \times \transid)$};
\node (m2) at (2,1) {\scriptsize$\tens(I \times \transid)\iota_2 \tens$};
\node (m3) at (3,1) {\scriptsize$\transid \tens$};
\draw[doubletighteq] (t3) to (m3);
\draw[doubletighteq] (t0) to (m0);
\draw[doubleloose] (t0) to node[above]{\scriptsize $\looseid (l \times \looseid)$} (t1);
\draw[doubletighteq] (t1) to  (t3);
\draw[doubleloose] (m0) to node[above]{\scriptsize$\alpha \looseid$} (m1);
\draw[doubletighteq] (m1) to (m2);
\draw[doubleloose] (m2) to node[above]{\scriptsize$l \looseid$} (m3);
\node at (1.5,1.5) {\scriptsize$\DDownarrow \lambda$};
\end{tikzpicture}}
\end{aligned}
\end{equation}

\begin{equation}\label{eq:lcbcrho}
\begin{aligned}
{\scriptstyle \begin{tikzpicture}[xscale=3.5]
\node (t0) at (0,2) {\scriptsize$\transid \tens $};
\node (t1) at (1.5,2) {\scriptsize$\tens(\transid \times \transid)$};
\node (t3) at (3,2) {\scriptsize$\tens ( \transid \times \tens)( \transid \times \transid \times I)$};
\node (m0) at (0,1) {\scriptsize$\transid \tens $};
\node (m1) at (1,1) {\scriptsize$\tens(\transid \times I )\iota_1\tens$};
\node (m2) at (2,1) {\scriptsize$\tens(\tens \times \transid)(\transid \times \transid \times I)$};
\node (m3) at (3,1) {\scriptsize$\tens ( \transid \times \tens)( \transid \times \transid \times I)$};
\draw[doubletighteq] (t3) to (m3);
\draw[doubletighteq] (t0) to (m0);
\draw[doubletighteq] (t0) to  (t1);
\draw[doubleloose] (t1) to node[above]{\scriptsize$\looseid (\looseid \times r^{-1})$} (t3);
\draw[doubleloose] (m0) to node[above]{\scriptsize$r^{-1} \looseid$} (m1);
\draw[doubletighteq] (m1) to (m2);
\draw[doubleloose] (m2) to node[above]{\scriptsize$\alpha \looseid$} (m3);
\node at (1.5,1.5) {\scriptsize$\DDownarrow \rho$};
\end{tikzpicture}}
\end{aligned}
\end{equation}

A monoidal object is {\bf braided} if, in addition, it is equipped with a loose 2-cell $\sigma_A: \tens \looseRightarrow{} \mathord{\ten} \tau$, where $\tau: A \times A \rightarrow A \times A$ interchanges the two copies of $A$; and if there are invertible globular 3-cells 
\begin{equation}\label{eq:R}
\begin{aligned}
\begin{tikzpicture}[xscale=2.7, yscale=1.5]
\node (t0) at (0,2) {\scriptsize$\tens (\tens \times \transid)$};
\node (t1) at (1,2) {\scriptsize$\tens(\transid \times \tens)$};
\node (t2) at (2,2) {\scriptsize$\tens \tau (\transid \times \tens)$};
\node (t3) at (3,2) {\scriptsize\begin{tabular}{c}$ \tens (\tens \times \transid)$\\ $\tau((\transid \times \transid) \times \transid))$\end{tabular}};
\node (t4) at (4,2) {\scriptsize\begin{tabular}{c}$\tens (\transid \times \tens)$\\ $\tau((\transid \times \transid) \times \transid)$\end{tabular}};
\node (m0) at (0,1) {\scriptsize$\tens (\tens \times \transid)$};
\node (m1) at (1,1) {\scriptsize$\tens (\tens \tau \times \transid)$};
\node (m2) at (2,1) {\scriptsize\begin{tabular}{c}$\tens (\tens \times \transid)$\\$(\tau(\transid \times \transid) \times \transid)$\end{tabular}};
\node (m3) at (3,1) {\scriptsize\begin{tabular}{c}$\tens (\transid \times \tens)$\\$(\tau(\transid \times \transid) \times \transid)$\end{tabular}};
\node (m4) at (4,1) {\scriptsize\begin{tabular}{c}$\tens (\transid \times \tens \tau)$\\$(\tau(\transid \times \transid) \times \transid)$\end{tabular}};
\draw[doubletighteq] (t4) to (m4);
\draw[doubletighteq] (t0) to (m0);
\draw[doubleloose] (t0) to node[above]{\scriptsize$\alpha$} (t1);
\draw[doubleloose] (t1) to node[above]{\scriptsize$\sigma \looseid$}
(t2);
\draw[doubletighteq] (t2) to (t3);
\draw[doubleloose] (t3) to node[above]{\scriptsize$\alpha \looseid$} (t4);
\draw[doubleloose] (m0) to node[below]{\scriptsize$\looseid(\sigma \times  \looseid)$} (m1);
\draw[doubletighteq] (m1) to (m2);
\draw[doubleloose] (m2) to node[below]{\scriptsize$\alpha \looseid$} (m3);
\draw[doubleloose] (m3) to node[below]{\scriptsize$\looseid (\looseid \times \sigma)$} (m4);
\node at (2,1.5) {\scriptsize$\DDownarrow R$};
\end{tikzpicture}
\end{aligned}
\end{equation}
\begin{equation}\label{eq:S}
\begin{aligned}
\begin{tikzpicture}[xscale=2.7, yscale=1.5]
\node (t0) at (0,2) {\scriptsize$\tens(\transid \times \tens) $};
\node (t1) at (1,2) {\scriptsize$\tens (\tens \times \transid)$};
\node (t2) at (2,2) {\scriptsize$ \tens \tau (\tens \times \transid)$};
\node (t3) at (3,2) {\scriptsize\begin{tabular}{c} $\tens (\transid \times \tens) $\\$ \tau((\transid \times \transid) \times \transid))$\end{tabular}};
\node (t4) at (4,2) {\scriptsize \begin{tabular}{c} $\tens (\tens \times \transid)$\\ $\tau((\transid \times \transid) \times \transid)$ \end{tabular}};
\node (m0) at (0,1) {\scriptsize$ \tens(\transid \times \tens)$};
\node (m1) at (1,1) {\scriptsize$\tens(\transid \times \tens \tau)$};
\node (m2) at (2,1) {\scriptsize\begin{tabular}{c}$ \tens(\transid \times \tens)$\\$(\transid \times \tau(\transid \times \transid) )$ \end{tabular}};
\node (m3) at (3,1) {\scriptsize\begin{tabular}{c}$ \tens (\tens \times \transid)$\\
$(\transid \times \tau(\transid \times \transid) $\end{tabular}};
\node (m4) at (4,1) {\scriptsize\begin{tabular}{c} $ \tens ( \tens \tau \times \transid)$\\$(\transid \times \tau(\transid \times \transid)$\end{tabular}};
\draw[doubletighteq] (t4) to (m4);
\draw[doubletighteq] (t0) to (m0);
\draw[doubleloose] (t0) to node[above]{\scriptsize$ \alpha^{-1}$} (t1);
\draw[doubleloose] (t1) to node[above]{\scriptsize$ \sigma \looseid$}
(t2);
\draw[doubletighteq] (t2) to (t3);
\draw[doubleloose] (t3) to node[above]{\scriptsize$\alpha^{-1} \looseid$} (t4);
\draw[doubleloose] (m0) to node[below]{\scriptsize$\looseid(\looseid \times \sigma)$} (m1);
\draw[doubletighteq] (m1) to (m2);
\draw[doubleloose] (m2) to node[below]{\scriptsize$\alpha^{-1} \looseid$} (m3);
\draw[doubleloose] (m3) to node[below]{\scriptsize$ \looseid (\sigma \times \looseid )$} (m4);
\node at (2,1.5) {\scriptsize$\DDownarrow S$};
\end{tikzpicture}
\end{aligned}
\end{equation}
satisfying axioms analogous to (BA1), (BA2), (BA3), and (BA4) given in~\cite[p136--139]{mccrudden:bal-coalgb} . 
It is {\bf sylleptic} when it is additionally equipped with an invertible globular 3-cell
\begin{equation}
\begin{aligned}
 \begin{tikzpicture}
 \node (tl) at (-2,1) {\scriptsize$\ten$};
 \node (tr) at (2,1) {\scriptsize$\ten$};
 \node (b) at (0,-.25) {\scriptsize$\tens \tau$};
 \draw[double] (tl)  -- (tr);
 \draw[doubleloose] (tl) to node[left, yshift=-5pt]{\scriptsize$\sigma$} (b);
 \draw[doubleloose] (b) to node[right, yshift=-5pt] {\scriptsize$\sigma$}(tr);
 \node at (0,0.5) {\scriptsize $\DDownarrow \upsilon \iso$}; 
 \end{tikzpicture}
 \end{aligned}
 \end{equation}
  satisfying the axioms analogous to (SA1), (SA2) on~\cite[p144--145]{mccrudden:bal-coalgb}. It is {\bf symmetric} if in addition, it satisfies the axiom analogous to the equation given on~\cite[p91]{mccrudden:bal-coalgb}.
\end{defn}

We now define monoidal 1-cells, 2-cells, and 3-cells.
As in the 1-categorical case, a morphism (1-cell) between monoidal objects can be lax, colax, or strong, corresponding to the direction of the comparison morphisms relating $fa \ten fb$ to $f(a\ten b)$.
But there is no need to define the two independently, as they are related by the duality $\fB\lco$ from \autoref{rmk:lcbop} that reverses loose 2-cells.

\begin{defn}\label{defn:lax-1cell}
Let $A,B$ be monoidal objects in \fB. A 1-cell $f:A \onecell B$ is {\bf lax monoidal} when it is equipped with the following loose 2-cells:
\begin{itemize}
\item $\chi: \mathord{\ten} (f \times f) \looseRightarrow{} f  \mathord{\otimes}$
\item $\iota: I_B \looseRightarrow{} fI_A $
\end{itemize}
as well as globular invertible 3-cells 
\begin{align*}
& \omega:  \looseid_{\tens}(\chi \times \looseid_f)  \horc  \chi\looseid_{\tens \times \transid} \horc  \looseid_f \alpha \RRightarrow \alpha\looseid_{f \times f \times f}  \horc \looseid_{\tens}(\looseid \times \chi)  \horc \chi \looseid_{\transid \times \tens}  \\
 &\gamma: \looseid_{\tens}(\iota_f \times \looseid_f) \looseid_{i_2} \horc \chi \looseid_{I \times \transid} \looseid_{i_2} \horc \looseid_f l\RRightarrow l \looseid_f \\
 &\delta:  \looseid_f r^{-1} \RRightarrow r^{-1} \looseid_f \horc \looseid_{\tens} (\looseid \times \iota) \looseid_{i_1 f} \horc \chi \looseid_{(\transid \times I)i_1}
\end{align*}
which satisfy the two axioms~\eqref{eq:laxfunc1},~\eqref{eq:laxfunc2} in Appendix~\ref{ap:coherence}, expressing the usual associativity and unitality conditions~\cite[Definition 4.10]{nick:tricatsbook}.

If $A$ and $B$ are braided, then $f$ is called {\bf braided} when there is a globular 3-cell 
\begin{equation*}
u: \sigma_B \looseid_{f \times f} \horc \chi  \looseid_{\tau} \looseRightarrow \chi \horc (\looseid_f \sigma_A)
\end{equation*}
satisfying braiding axioms analogous to (BHA1) and (BHA2) given in~\cite[p141-142]{mccrudden:bal-coalgb}. 
If $A$ and $B$ are sylleptic or symmetric, then $f$ is {\bf sylleptic} or {\bf symmetric}, respectively, if it is braided and the 3-cells defining the braided monoidal structure satisfy an additional axiom analogous to (SHA1) of~\cite[p145]{mccrudden:bal-coalgb}.

A \textbf{(braided, sylleptic, or symmetric) colax monoidal 1-cell} in \fB  is simply a lax one in $\fB\lco$.
Thus it is equipped with loose 2-cells $\bar{\chi}$ and $\bar{\iota}$ which go in the opposite direction of $\chi$ and $\iota$, respectively, and appropriate invertible 3-cells.


If $f$ is both lax and colax monoidal, the associated loose 2-cells $\chi$ and $\iota$ form adjoint equivalences with their colax counterparts, and the 3-cells correspond to their colax counterparts as mates under the adjoint equivalence structure, it is {\bf strong monoidal}. Details on the mate correspondences for $\omega$, $\gamma$ and $\delta$ aregiven in the second author's PhD thesis~\cite[B0.7,B0.8,B0.9]{westerthesis}.
\end{defn}

By construction, these definitions give the expected results in \fBicat. 
In \fDblf, where there are no nonidentity 3-cells, they reduce to the definitions from section~\ref{sec:symm-mono-double}; and in particular, every syllepsis is a symmetry.

Note that since \autoref{defn:lax-1cell} contains no tight 2-cells and all the 3-cells are isomorphisms, lax monoidal 1-cells in the tight 2-cell dual $\fB\tco$ coincide with those in $\fB$ itself.

For ordinary monoidal categories, once a type of 1-cell (lax, colax, strong) is fixed, there is only one kind of 2-cell between them.
But for monoidal bicategories, a 2-cell can also be lax, colax, or strong at the 2-dimensional level, independently of whether its domain and codomain 1-cells are lax, colax, or strong.
(This has nothing to do with laxity of its underlying natural transformation; for us these are always strong, i.e.\ pseudonatural.)
But again, there is no need to define them independently, as they are related by the other duality $\fB\tco$ that reverses the tight 2-cells.

\begin{defn}\label{Def:mon1cell}
Let $f, g:A \onecell B$ be lax monoidal 1-cells in \fB. A {\bf lax monoidal 2-cell} $\beta: f \looseRightarrow g$ is a loose 2-cell in \fB\ that is equipped with globular 3-cells
\begin{itemize}
\item $\Pi:  \looseid_{\ten}(\beta \times \beta) \horc \chi_g \RRightarrow{} \chi_f \horc \beta  \looseid_{\ten} $
\item $M:  \looseid_{I} \horc \iota_g \RRightarrow{} \iota_f \horc \beta  \looseid_{I_A}$
\end{itemize}
such that coherence equations \eqref{eq:mon2cell1}, \eqref{eq:mon2cell2}, and \eqref{eq:mon2cell3} in Appendix~\ref{ap:coherence} hold. Applied to the special case of bicategories this gives us equations (TA2), (TA3) and (TA4) of~\cite{gg:ldstr-tricat}.

Such a 2-cell is {\bf braided}, {\bf sylleptic} or {\bf symmetric} when $f,g$ are so and a coherence axiom analogous to~\cite[(BTA1) p143]{mccrudden:bal-coalgb} holds.

A {\bf colax monoidal 2-cell} between lax monoidal 1-cells in \fB\ is a lax monoidal 2-cell in $\fB\tco$; thus it is equipped with morphisms $\bar{\Pi}, \bar{M}$ in the opposite direction from $\Pi$ and $M$.
A {\bf strong monoidal 2-cell} is a lax (or equivalently colax) one whose $\Pi$ and $M$ are isomorphisms.

If instead $f,g$ are colax monoidal 1-cells in \fB, then a \textbf{lax monoidal 2-cell} between them is a lax monoidal 2-cell between $f,g$ regarded as lax monoidal 1-cells in $\fB\lco$, while a \textbf{colax monoidal 2-cell} is a lax monoidal 2-cell between $f,g$ regarded as lax monoidal 2-cells in $\fB\tlco$.
\end{defn}



When interpreted in $\fBicat$, monoidal 2-cells correspond to monoidal pseudonatural transformations. Such transformations are less general than tritansformations between one-object tricategories; they are obtained by restricting the constraints $\Pi$ and $M$ in a tritransformation to be given by icons instead of adjoint equivalences, leading to the notion of a pseudo-icon~\cite[Definition 5]{gg:ldstr-tricat} from which the definition above is derived. This notion equals, up to canonical isomorphisms, that of a monoidal pseudo transformation given in~\cite{sp:thesis}. In $\fDblf$, monoidal 2-cells correspond to monoidal double transformations.

As remarked above, we will actually construct a locally cubical bicategory of monoidal objects.
The monoidal 2-cells will be the loose 2-cells therein; we now define the tight ones and the 3-cells.

\begin{defn}\label{Def:monicon}
  Let $f, g:A \rightarrow B$ be lax monoidal 1-cells in \fB.
  A \textbf{monoidal icon} $\beta: f \Rightarrow g$ is a (tight) 2-morphism in \fB\ that is equipped with (non-globular!)\ 3-cells
\begin{equation}
\begin{aligned}
 \begin{tikzpicture}[scale=2]
 \node (tl) at (0,1) {\scriptsize$I_B$};
 \node (tr) at (1,1) {\scriptsize$f I_A$};
 \node (bl) at (0,0) {\scriptsize$I_B$};
 \node (br) at (01,0) {\scriptsize$g I_A$}; 
 \draw[doubleloose] (tl)  to node[above]{\scriptsize$\iota_f$} (tr);
 \draw[doubleeq] (tl) to (bl);
 \draw[doubleloose] (bl) to node[below] {\scriptsize$\iota_g$}(br);
  \draw[doubletight] (tr) to node[right] {\scriptsize$\beta \tightid_I$}(br);
 \node at (0.5,0.5) {\scriptsize$\DDownarrow N^{\beta}$}; 
 \end{tikzpicture}
 \end{aligned}
 \hspace{.5cm}
 \begin{aligned}
  \begin{tikzpicture}[scale=2]
 \node (tl) at (0,1) {\scriptsize$\ten (f \times f)$};
 \node (tr) at (1,1) {\scriptsize$f \ten$};
 \node (bl) at (0,0) {\scriptsize$\ten(g \times g)$};
 \node (br) at (01,0) {\scriptsize$g  \ten$}; 
 \draw[doubleloose] (tl)  to node[above]{\scriptsize$\chi_f$} (tr);
 \draw[doubletight] (tl) to node[left]{\scriptsize$\tightid_{\ten} (\beta \times \beta)$} (bl);
 \draw[doubleloose] (bl) to node[below] {\scriptsize$\chi_g$}(br);
  \draw[doubletight] (tr) to node[right] {\scriptsize$\beta \tightid_{\ten}$}(br);
 \node at (0.5,0.5) {\scriptsize $\DDownarrow \Sigma^{\beta}$}; 
 \end{tikzpicture}
\end{aligned}
\end{equation}
such that the coherence axioms~\eqref{eq:monicon1},~\eqref{eq:monicon2},~\eqref{eq:monicon3} in Appendix~\ref{ap:coherence} hold. In the case of bicategories, these axioms specialize to (TA2), (TA3) and (TA4) of~\cite{gg:ldstr-tricat}.

A monoidal icon is {\bf braided}, {\bf sylleptic} or {\bf symmetric} when $f,g$ are so and the coherence axiom~\eqref{eq:bricon} in Appendix~\ref{ap:coherence} holds.
(We give this axiom explicitly, even though it pertains to the braided case which we are generally omitting, because braided monoidal icons do not seem to appear anywhere in the existing literature.)

If $f, g:A \rightarrow B$ are instead colax monoidal 1-cells, a \textbf{monoidal icon} between them is a monoidal icon between $f,g$ regarded as lax monoidal 1-cells in $\fB\lco$.
\end{defn}

Monoidal icons interpreted in $\fBicat$ correspond to ico-icons~\cite[Definition 2]{gg:ldstr-tricat} between one-object tricategories. These are doubly degenerate colax tritransformations between lax trihomomorphisms which agree on both 0- and 1-cells.
In $\fDblf$, since there are no nonidentity tight 2-cells, there are no nonidentity monoidal icons.
Note that a monoidal icon $f\to g$ in $\fB\tco$ is the same as a monoidal icon $g\to f$ in $\fB$.

\begin{defn}
  Let $f,g,f',g': A \rightarrow B$ be lax monoidal 1-cells, let $\alpha: f \looseRightarrow{} g$, $\beta: f' \looseRightarrow{} g'$ be lax monoidal 2-cells, and let $\gamma: f \Rightarrow f'$, $\delta: g \Rightarrow g'$ be monoidal icons. A \textbf{monoidal 3-cell} is a 3-cell 
  
\begin{equation*}
\begin{aligned}
 \begin{tikzpicture}[scale=2]
 \node (tl) at (0,1) {\scriptsize$f$};
 \node (tr) at (1,1) {\scriptsize$g$};
 \node (bl) at (0,0) {\scriptsize$f'$};
 \node (br) at (01,0) {\scriptsize$g'$}; 
 \draw[doubleloose] (tl)  to node[above]{\scriptsize$\alpha$} (tr);
 \draw[doubletight] (tl) to node[left]{\scriptsize$\gamma$} (bl);
 \draw[doubleloose] (bl) to node[below] {\scriptsize$\beta$}(br);
  \draw[doubletight] (tr) to node[right] {\scriptsize$\delta$}(br);
 \node at (0.5,0.5) {\scriptsize$\DDownarrow \Gamma$}; 
 \end{tikzpicture}
 \end{aligned}
 \end{equation*}
such that the coherence axioms~\eqref{eq:mon3cell1} and~\eqref{eq:mon3cell2} in Appendix \ref{ap:coherence} hold.


If $f,g,f',g', \alpha, \beta, \gamma$ and $\delta: g \Rightarrow g'$ are braided, sylleptic, or symmetric monoidal 1-cells, 2-cells, and icons, respectively, then a \textbf{braided, sylleptic, or symmetric monoidal 3-cell} $\Gamma$ as depicted above, is simply a monoidal 3-cell. 
\end{defn}

Monoidal 3-cells in $\fBicat$ correspond to cubical pseudo-icon modifications~\cite[Definition 16]{gg:ldstr-tricat} between one-object tricategories. When the icon components are identities, this reduces to a monoidal modification~\cite[Definition 2.8]{sp:thesis}. 

We will show that monoidal objects, lax monoidal 1-cells, lax monoidal 2-cells, monoidal icons, and monoidal 3-cells in a locally cubical bicategory \fB\ form a locally cubical bicategory $\cM on_{ll}\fB$.
The first subscript $l$ indicates the laxity of the 1-cells, the second that of the loose 2-cells; in the colax case we write $c$ and in the strong case $p$ (for ``pseudo'' --- the subscript $s$ is usually used in 2-category theory instead for ``strict'').
We write $v,w \in \{l,c,p\}$ for unspecified ``laxness annotations''.
Thus by definition we will have
\begin{align*}
  \cMon_{cl}(\fB) &= \cMon_{ll}(\fB\lco)\\
  \cMon_{lc}(\fB) &= \cMon_{ll}(\fB\tco)\tco\\
  \cMon_{cc}(\fB) &= \cMon_{ll}(\fB\tlco)\tco.
\end{align*}
Similarly, braided, sylleptic and symmetric cells will form locally cubical bicategories $\cB r \cM on_{vw}\fB$, $\cS yl \cM on_{vw}\fB$, and $\cS ym \cM on_{vw}\fB$.

\begin{prop}\label{prop:dc}
Let $A,B$ be monoidal objects in a 1-strict locally cubical bicategory with products. The hom-spaces $\cM on_{vw}\fB (A,B)$, $\cB r \cM on_{vw}\fB(A,B)$, $\cS yl \cM on_{vw}\fB(A,B)$, and $\cS ym \cM on_{vw}\fB(A,B)$ are double categories for $v,w \in \{l,c,p\}$.
\end{prop}

\begin{proof}
By duality, the colax cases are implied by the lax ones.
The strong 2-cell cases ($w=)p$ are generally trivial (the property of invertibility is always preserved), but the strong 1-cell cases ($v=p$) require carrying through the adjoint equivalence data.

First we show that 1-cells and icons in the respective hom-spaces form a category. For each lax monoidal 1-cell $f:A \rightarrow B$, the identity icon $\tightid_f$ is a lax monoidal icon with the 3-cells $N^{\tightid_f} := \tightid_{\iota_f}$ and $\Sigma^{\tightid_f} := \tightid_{\chi_f}$. This is well-defined, because the functor ``$\comp$" preserves tight identities. The coherence equations are trivially satisfied.  For each two lax monoidal 1-cells $f,g$ and lax monoidal icons $\alpha, \beta: f \Rightarrow g$, the composite icon $\beta \verc \alpha$ can be equipped with the lax monoidal structure given by the composites $N^{\beta \verc \alpha} := N^{\beta} \verc N^{\alpha}$ and $\Sigma^{\beta \verc \alpha} : = \Sigma^{\beta} \verc \Sigma^{\alpha}$.  We have a strict interchange law between $\verc$ and $\comp$, induced by functoriality of $\comp$, so these 3-cells are well-defined. The coherence conditions~\eqref{eq:monicon1}--\eqref{eq:monicon3} hold by componentwise application of the coherence equalities for $N^{\beta \verc \alpha}$ and $\Sigma^{\beta \verc \alpha}$.

For strong 1-cells and icons we need to verify that the lax structure cells $N^{\tightid_f}$, $N^{\beta \verc \alpha}, \Sigma^{\tightid_f}$, and $\Sigma^{\beta \verc \alpha} $ are inverse in the loose direction to their colax counterparts. For $N^{\tightid_f}$ and $\Sigma^{\tightid_f}$, this follows from functoriality of ``$\horc$". For $N^{\beta \verc \alpha}$ and $\Sigma^{\beta \verc \alpha}$, this follows from the the fact that the statement is true for their components combined with the exchange law between ``$\horc$" and ``$\verc$" and strictness of ``$\verc$".
When $f$ and $g$ are braided, sylleptic or symmetric, the same data satisfies the coherence equation for braided monoidal icons.

We also need to show that lax monoidal 2-cells and monoidal 3-cells form a category. For every lax monoidal 2-cell $\alpha: f \looseRightarrow{} g$, the identity 3-cell $\tightid_{\alpha}$ in $\fB$  is lax monoidal. The required two equations~\eqref{eq:mon3cell1},~\eqref{eq:mon3cell2} are trivially satisfied.
For any two monoidal 3-cells $L: \alpha \Rightthreecell \beta$, $K:\beta \Rightthreecell \gamma$, the composition $K \verc L$ in $\fB$ is a monoidal 3-cell. The equations for monoidal 3-cells~\eqref{eq:mon3cell1},~\eqref{eq:mon3cell2} hold by sequential application of the respective equations for $L$ and $K$.
The colax and strong cases follow automatically, and similarly the braided, sylleptic, and symmetric cases since in this case the 3-cells have no additional data.

Now we describe the loose structure; we need to show that $\horc$ and $\looseid$ are well-defined as the functors which give the loose structure in the new double category. To see this, recall that $\horc$ and $\looseid$ correspond to the functors $\odot$ and $U$, respectively, given  in Definition~\ref{def:dblcat}.
Let $f$ be a lax monoidal 1-cell. The loose identity 2-cell $\looseid_f$ is a strong monoidal 2-cell with monoidal structure given by the composites of coherence cells $\horl$, $\horr$ shown in Figure~\ref{fig:idf}

\begin{figure}
\begin{equation}
M^{\looseid_f}:=
\begin{aligned}
 \begin{tikzpicture}[yscale=1.5, xscale=3]
 \node (tl) at (0,0) {\scriptsize$I_B$};
\node (tr) at (1,0) {\scriptsize$f   I_A$};
 \node (tm) at (.5,0) {\scriptsize$f  I_A$};
 \node (bl) at (0,1) {\scriptsize$I_B$};
 \node (bm) at (.5,1) {\scriptsize$I_B$};
 \node (br) at (1,1) {\scriptsize$f I_A$}; 
 \draw[doubleloose] (tl)  to node[below]{\scriptsize$\iota_f$} (tm);
  \draw[doubleloose] (tm)  to node[below]{\scriptsize$\looseid_f \looseid_I$} (tr);
 \draw[doubleeq] (tl) to (bl);
  \draw[doubleloose] (bl) to node[above] {\scriptsize$\looseid_I$}(bm);
 \draw[doubleloose] (bm) to node[above] {\scriptsize$\iota_f$}(br);
  \draw[=] (tr) to (br);
 \node at (0.5,0.5) {\scriptsize $\DDownarrow \iso $}; 
 \end{tikzpicture}
 \end{aligned}
\end{equation}
\begin{equation}
 \Pi^{\looseid_f}:=
 \begin{aligned}
  \begin{tikzpicture}[yscale=1.5, xscale=5]
 \node (tl) at (0,0) {\scriptsize$\ten  (f \times f)$};
 \node (tr) at (1,0) {\scriptsize$f  \ten$};
 \node (bl) at (0,1) {\scriptsize$\ten  (f \times f)$};
 \node (br) at (01,1) {\scriptsize$f \ten$}; 
 \node(tm) at (.5,0) {\scriptsize$f \ten$};
 \node (bm) at (.5,1) {\scriptsize$\ten (f\times f)$};
 \draw[doubleloose] (tl)  to node[below]{\scriptsize$\chi_f$} (tm);
  \draw[doubleloose] (tm)  to node[below]{\scriptsize$\looseid_f \looseid{\ten}$} (tr);
 \draw[=] (tl) to (bl);
  \draw[doubleloose] (bl) to node[above] {\scriptsize$\looseid_{\ten}(\looseid_f \times \looseid_f)$}(bm);
 \draw[doubleloose] (bm) to node[above] {\scriptsize$\chi_f$}(br);
  \draw[=] (tr) to (br);
 \node at (0.5,0.5) {\scriptsize$ \DDownarrow$ $\iso$}; 
 \end{tikzpicture}
\end{aligned}
\end{equation}
\caption{Strong monoidal constraints for $\looseid_f$}\label{fig:idf}
\end{figure}

The conditions for monoidal 3-cells follow from the naturality conditions of the coherence cells. 
Let $\gamma$ be a monoidal icon, one can verify that the loose identity 3-cell $\looseid_{\gamma}$ is a monoidal 3-cell: The loose source and target 2-cells of $\looseid_{\gamma}$, are loose identities; hence, the coherence condition holds by naturality of $\horl$ and $\horr$. 

Let $\alpha:f \looseRightarrow{} g$ and $\beta: g \looseRightarrow{} h$ be two lax monoidal 2-cells. Their composition $\alpha \horc \beta$ is lax monoidal with the structure 3-cells given in Figure~\ref{fig:horc}.
The coherence equations are satisfied by sequential application of the respective equation for $\alpha$ and $\beta$, applications of the exchange law between loose and tight composition, together with simple manipulations of coherence cells.

\begin{figure}
\begin{equation}
\begin{aligned}
M^{\alpha \horc \beta} := 
\end{aligned}
\begin{aligned}
 \begin{tikzpicture}[yscale=1.5, xscale=3]
 \node (tl) at (0,-1) {\scriptsize$I_B$};
\node (tr) at (1,-1) {\scriptsize$g   I_A$};
 \node (tm) at (.5,-1) {\scriptsize$f  I_A$};
 \node (bl) at (0,0) {\scriptsize$I_B$};
 \node (bm) at (0.5,0) {\scriptsize$I_B$};
 \node (br) at (1,0) {\scriptsize$g I_A$}; 
 \node (trr) at (1.5,-1) {\scriptsize$h I_A$};
 \node (brr) at (1.5,0) {\scriptsize$h I_A$};
 \node (bbr) at (1.5,1) {\scriptsize$hI_A$};
  \node (bbm1) at (.5,1) {\scriptsize$I_B$};
 \node (bbm) at (1,1) {\scriptsize$I_B$};
 \node(bbl) at (0,1) {\scriptsize$I_B$};
    \draw[doubleloose] (tm) to[in=220, out=-60] node[below]{\scriptsize$(\alpha \horc \beta)\looseid_{I}$} (trr);
 \draw[doubletight] (brr) to node[right] {} (bbr);
 \draw[doubleeq] (bl) to (bbl);
  \draw[doubleloose] (bbl) to node [below]{\scriptsize$\looseid_{I}$} (bbm1);
    \draw[doubleloose] (bbm1) to node [below]{\scriptsize$\looseid_{I}$} (bbm);
 \draw[doubleloose] (bbm) to node [below]{\scriptsize$\iota_{h}$} (bbr);
 \draw[doubleloose] (tr) to node[below]{\scriptsize$\beta \looseid_I$} (trr);
  \draw[doubleloose] (br) to node[below]{\scriptsize$\beta \looseid_I$}(brr);
  \draw[doubleeq] (trr) to (brr);
 \draw[doubleloose] (tl)  to node[below]{\scriptsize$\iota_f$} (tm);
  \draw[doubleloose] (tm)  to node[below]{\scriptsize$\alpha \looseid_I$} (tr);
 \draw[doubleeq] (tl) to (bl);
  \draw[doubleloose] (bl) to node[above] {\scriptsize$\looseid_I$}(bm);
 \draw[doubleloose] (bm) to node[above] {\scriptsize$\iota_g$}(br);
 \draw[doubleloose] (bbl) to[in=120, out=60] node[above]{$\looseid_I$} (bbm);
  \draw[doubleeq] (tr) to (br);
   \draw[doubleeq] (bm) to (bbm1);
 \node at (0.5,-.5) {\scriptsize $M^{\alpha} \DDownarrow  $}; 
  \node at (1,.5) {\scriptsize $M^{\beta} \DDownarrow $}; 
 \node at (1.25,-.5) {\scriptsize $=$}; 
 \node at (1,-1.25) {\scriptsize$\iso$};
 \node at (0.5,1.25) {\scriptsize$\iso$};
 \end{tikzpicture}
 \end{aligned}
\end{equation}

\begin{equation}
 \Pi^{\alpha \horc \beta}:=
 \begin{aligned}
  \begin{tikzpicture}[yscale=1.5, xscale=5]
 \node (tl) at (0,-1) {\scriptsize$\ten  (f \times f)$};
 \node (tr) at (1,-1) {\scriptsize$g \ten$};
 \node (bl) at (0,0) {\scriptsize$\ten  (f \times f)$};
 \node (br) at (01,0) {\scriptsize$g \ten$}; 
 \node(tm) at (.5,-1) {\scriptsize$f \ten$};
 \node (bm) at (.5,0) {\scriptsize$\ten (g\times g)$};
 \node (trr) at (1.5,-1) {\scriptsize$h \ten$};
  \node (brr) at (1.5,0) {\scriptsize$h \ten$};
  \node (bbl) at (0,1) {\scriptsize$\ten (f \times f)$};
  \node (bbm) at (.5,1) {\scriptsize$\ten (g \times g)$}; 
  \node (bbr) at (1,1) {\scriptsize$\ten (h \times h)$};
  \node (bbrr) at (1.5,1) {\scriptsize$h \ten $};
 \draw[doubleloose] (tl)  to node[below]{\scriptsize$\chi_f $} (tm);
  \draw[doubleloose] (tm)  to node[below]{\scriptsize$\alpha \looseid_{\ten}$} (tr);
 \draw[doubleeq] (tl) to (bl);
  \draw[doubleloose] (bl) to node[above] {\scriptsize$\looseid_{\ten} (\alpha \times \alpha)$}(bm);
 \draw[doubleloose] (bm) to node[above] {\scriptsize$\chi_g$}(br);
  \draw[doubleeq] (tr) to (br); 
 \draw[doubleeq] (trr) to (brr);
 \draw[doubleloose] (tr) to node[below]{\scriptsize$\beta \looseid_{\ten}$} (trr);
 \draw[doubleloose] (br) to node[below]{\scriptsize$\beta \looseid_{\ten}$} (brr);
 \draw[doubleloose] (bbr) to node[below]{\scriptsize$\chi_h$} (bbrr);
 \draw[doubleeq] (bl) to (bbl);
 \draw[doubleeq] (bm) to (bbm);
 \draw[doubleeq] (brr) to (bbrr);
 \draw[doubleloose] (bbl) to node[below]{\scriptsize$\looseid_{\ten} (\alpha \times \alpha)$} (bbm);
  \draw[doubleloose] (bbm) to node[below]{\scriptsize$\looseid_{\ten} (\beta \times \beta)$} (bbr);
   \draw[doubleloose] (tm) to[in=220, out=-60] node[below]{\scriptsize$(\alpha \horc \beta)\looseid_{\ten}$} (trr);
   \draw[doubleloose] (bbl) to[in=120, out=60] node[above]{\scriptsize$\looseid_{\ten} \comp (\alpha \horc \beta)\times (\alpha \horc \beta)$} (bbr);
    \node at (0.5,-.5) {\scriptsize $\DDownarrow  \Pi^{\alpha}$};
  \node at (1.25,-.5) {\scriptsize $=$};
  \node at (0.25,.5) {\scriptsize $=$};
  \node at (1,.5) {\scriptsize $\DDownarrow  \Pi^{\beta}$};
  \node at (1,-1.2) {\scriptsize$\iso$};
 \node at (.5,1.2) {\scriptsize$\iso$};
 \end{tikzpicture}
\end{aligned}
\end{equation}
\caption{Monoidal constraints for $\alpha\horc\beta$}\label{fig:horc}
\end{figure}

Let $\Gamma$ and $\Delta$ be monoidal 3-cells. Their composite $\Gamma \horc \Delta$ is again monoidal. Again, the conditions for monoidal 3-cells follow directly from the conditions on the monoidal 3-cells $\Gamma$ and $\Delta$, applications of the exchange law between loose and tight composition, and simple manipulations of coherence cells.
The colax and strong cases follow automatically.

Let $f$ be a braided, sylleptic or symmetric monoidal 1-cell. The loose identity $\looseid_f$ is a braided, sylleptic or symmetric monoidal 2-cell, respectively, as the coherence equation~\cite[(BTA1) p143]{mccrudden:bal-coalgb} merely states that the 3-cell $u$ pasted with coherence 3-cells equals itself. Let $\alpha, \beta$ be braided, sylleptic, or symmetric monoidal 2-cells, the loose composition $\alpha \horc \beta$ is braided, sylleptic, or symmetric monoidal, respectively. One can verify that~\cite[(BTA1) p143]{mccrudden:bal-coalgb} holds by applying the exchange law between loose and tight composition, manipulation of coherence cells, and sequential application of the respective equations for $\alpha$ and $\beta$.  Braided, sylleptic and symmetric monoidal 3-cells are simply monoidal 3-cells; therefore, it follows that the images of $\horc$ and $\looseid$ of braided, sylleptic, or symmetric monoidal cells are well-defined in $\cB r \cM on_{vw}\fB(A,B)$, $\cS yl \cM on_{vw}\fB(A,B)$, and $\cS ym \cM on_{vw}\fB(A,B)$, respectively.

Functoriality of $\horc$ and $\looseid$ in $\cM on \fB(A,B)$ follows from their functoriality in $\fB(A,B)$. 
The unitality and associativity cells $\hora$, $\horl$, and $\horr$ are monoidal 3-cells. Unfolding the constraints for monoidal 3-cells shows that the 3-cells pasted together with coherence cells need to equal themselves. It follows from coherence of the functor $\horc$, that this is the case. Therefore,  $\cM on_{vw}\fB(A,B)$ is a double category for $v,w \in \{l,c,p\}$. 
Since braided, symmetric, or sylleptic monoidal 3-cells require no extra data, it follows that $\cB r \cM on_{vw}\fB(A,B)$, $\cS yl \cM on_{vw}\fB(A,B)$, and $\cS ym \cM on_{vw}\fB(A,B)$ are double categories for $v,w \in \{l,c,p\}$.
\end{proof}

\begin{thm}\label{thm:lcbc}
Let \fB\ be a 1-strict locally cubical bicategory with products. Monoidal objects, lax monoidal 1-cells, lax monoidal 2-cells, monoidal icons, and monoidal 3-cells in  \fB\ form a locally cubical bicategory $\cMon_{ll}\fB$, and similarly we have $\cMon_{vw}\fB$ for $v,w\in \{l,c,p\}$.
When the objects and cells are braided, sylleptic or symmetric,  we obtain the locally cubical bicategories $\cB r \cM on_{vw}\fB$, $\cS yl \cM on_{vw}\fB$, and $\cS ym \cM on_{vw}\fB$.
\end{thm}

\begin{proof}
We have established in Proposition~\ref{prop:dc}, that the respective hom-double-categories $\cM on_{vw} \fB(A,B)$, $\cB r \cM on_{vw}\fB$, $\cS yl \cM on_{vw}\fB$, and $\cS ym \cM on_{vw}\fB$ exist.

We need to check that the unit $\transid_A$ is a well-defined pseudo double functor from the trivial double category to the respective hom-double-categories, as well as braided, sylleptic and symmetric cells. 
The unit 1-cells $\transid_A$ are monoidal for all objects $A \in$ \fB,  with the monoidal structure $\xi, \iota$ given by the unitor 2-cells, and $\gamma, \delta$, and $\omega$ by coherence cells for the structure of the double category. The constraints boil down to simple manipulations of coherence cells, which hold by coherence of double categories. By functoriality of $\transid$, its image on the loose 2-cell  is isomorphic to the loose identity $\looseid_{\transid_A}$. This isomorphism gives rise to the lax monoidal structure on $\transid_{\transid_A}$ shown in Figure~\ref{fig:transid}.

\begin{figure}
\begin{equation}
M^{\transid_{\transid_A}}:=
\begin{aligned}
 \begin{tikzpicture}[yscale=1.8, xscale=3]
 \node (tl) at (0,0) {\scriptsize$I_A$};
\node (tr) at (1,0) {\scriptsize$\transid_A   I_A$};
 \node (tm) at (.5,0) {\scriptsize$\transid_A  I_A$};
 \node (bl) at (0,1) {\scriptsize$I_A$};
 \node (bm) at (.5,1) {\scriptsize$I_A$};
 \node (br) at (1,1) {\scriptsize$\transid_A I_A$}; 
 \draw[doubleloose] (tl)  to node[below]{\scriptsize$\iota_{\transid_A}$} (tm);
  \draw[doubleloose] (tm) to[in=220, out=-60] node[below] {\scriptsize$\transid_{\transid_A} \looseid_I$} (tr);
 \draw[doubleloose] (tm)  to node[above]{\scriptsize$\looseid_{\transid_A} \looseid_I$} (tr);
 \draw[doubleeq] (tl) to (bl);
  \draw[doubleloose] (bl) to node[above] {\scriptsize$\looseid_I$}(bm);
 \draw[doubleloose] (bm) to node[above] {\scriptsize$\iota_{\transid_A}$}(br);
  \draw[=] (tr) to (br);
 \node at (0.5,0.5) {\scriptsize $\DDownarrow$ $\iso $}; 
   \node at (0.75,-.1) {\scriptsize $ \DDownarrow$ $\iso$}; 
 \end{tikzpicture}
 \end{aligned}
\end{equation}
\begin{equation}
 \Pi^{\transid_{\transid_A}}:=
 \begin{aligned}
  \begin{tikzpicture}[yscale=1.8, xscale=5]
 \node (tl) at (0,0) {\scriptsize $\ten  (\transid_A \times \transid_A)$};
 \node (tr) at (1,0) {\scriptsize $\transid_A  \ten$};
 \node (bl) at (0,1) {\scriptsize $\ten  (\transid_A \times \transid_A)$};
 \node (br) at (01,1) {\scriptsize $\transid_A \ten$}; 
 \node(tm) at (.5,0) {\scriptsize $\transid_A \ten$};
 \node (bm) at (.5,1) {\scriptsize $\ten (\transid_A \times \transid_A)$};
 \draw[doubleloose] (tl)  to node[below]{\scriptsize $\chi_{\transid_A} $} (tm);
 \draw[doubleloose] (tm) to[in=220, out=-60] node[below]{\scriptsize $\transid_{\transid_A} \looseid_{\ten}$} (tr);
  \draw[doubleloose] (tm)  to node[above]{\scriptsize $\looseid_{\transid_A} \looseid{\ten}$} (tr);
 \draw[=] (tl) to (bl);
  \draw[doubleloose] (bl) to node[below] {\scriptsize $\looseid_{\ten}(\looseid_{\transid_A} \times \looseid_{\transid_A})$}(bm);
          \draw[doubleloose] (bl) to[in=120, out=60] node[above]{\scriptsize $\looseid_{\ten}(\transid_{\transid_A} \times \transid_{\transid_A})$}(bm); 
 \draw[doubleloose] (bm) to node[below] {\scriptsize $\chi_{\transid_A}$}(br);    
  \draw[=] (tr) to (br);
 \node at (0.5,0.5) {\scriptsize $ \DDownarrow$ $\iso$}; 
  \node at (0.75,-.1) {\scriptsize  $ \DDownarrow$ $\iso$}; 
    \node at (0.25,1.1) {\scriptsize  $ \DDownarrow$ $\iso$}; 
 \end{tikzpicture}
\end{aligned}
\end{equation}
\caption{Monoidal structure on $\transid_{\transid_A}$}\label{fig:transid}
\end{figure}

Coherence equations~\ref{eq:mon2cell1}, \ref{eq:mon2cell2}, \ref{eq:mon2cell3}, and~\cite[(BTA1) p143]{mccrudden:bal-coalgb} hold by simple manipulations of the isomorphisms. This makes $\transid_{\transid_A}$ a braided, sylleptic, or symmetric monoidal 2-cell if $\transid_A$ is braided, sylleptic, or symmetric, respectively, and the colax and pseudo cases are automatic.

By functoriality, the image of $\transid$ on the  tight 2-cell and 3-cell equal $\tightid_{\transid_A}$ and $\looseid_{\tightid_{\transid_A}} = \tightid_{\looseid_{\transid_A}}$, respectively. These cells are braided, symmetric or sylleptic; lax, colax, or strong monoidal, depending on $A$. It follows that $\transid_A$ is a well-defined functor from the trivial double category to the respective hom-double-categories.

Next, we need to show that monoidal structure is preserved by the composition along a 0-cell boundary.
For any two lax monoidal 1-cells $f:A \rightarrow B$, $g:B \rightarrow C$, the composite $g \comp f$ is monoidal with $\chi^{g \comp f}$ and $\iota^{g \comp f}$ defined below. 
\begin{align}
\chi_{g \comp f} &: \hspace{.5cm} &\otimes (gf \times gf) \xlooseRightarrow{\chi_g \looseid_{f \times f}} g \otimes (f \times f) \xlooseRightarrow{\looseid_g \chi_f} gf \tens \\
\iota_{g \comp f} & : \hspace{.5cm} &I_C \xlooseRightarrow{\iota_g} g I_B \xlooseRightarrow{\looseid_g \iota_f} gfI_A
\end{align}

The structure 3-cell $\gamma^{g\circ f}$ for is given in Figure~\ref{fig:gamma-circ}. For readability, we write "$\hspace{5pt}\widehat{}\hspace{5pt}$" to express that we have left out some coherence 3-cells, which are needed to make the source and target of the 3-cell well defined. 
The 3-cells $\delta^{g \comp f}$ and $\omega^{g \comp f}$ are defined similarly, and so is $u^{g \comp f}$ when $g, f$ are braided monoidal 1-cells. 

\begin{figure}
\begin{equation}
\begin{aligned}
 \begin{tikzpicture}[yscale=1.5, xscale=4.8]
 \node (t0) at (0,2) {\scriptsize $\tens(I_C \times gf)i_2$};
 \node (t1) at (.5,2) {\scriptsize $\tens(gI_B \times gf)i_2$};
\node (t2) at (1,2) {\scriptsize $g \tens (I_B \times f)i_2$};
 \node (t3) at (1.5,2) {\scriptsize $g \tens (fI_A \times f)i_2$};
  \node (t4) at (2,2) {\scriptsize $gf \tens (I_A \times \transid)i_2$};
 \node (t5) at (2.3,2) {\scriptsize $gf$};
  \node (m0) at (0,1) {\scriptsize $\tens(I_C \times g)i_2f$};
 \node (m1) at (.5,1) {\scriptsize $\tens(gI_B \times g)i_2f$};
\node (m2) at (1,1) {\scriptsize $g \tens (I_B \times \transid)i_2f$};
 \node (m5) at (2.3,1) {\scriptsize $gf$};
  \node (b0) at (0,0) {\scriptsize $\tens(I_C \times \transid)i_2 gf$};
 \node (b5) at (2.3,0) {\scriptsize $gf$};
  \draw[doubleloose] (0.2,2.2) to[in=120, out=60] node[above]{\scriptsize $\looseid_{\tens} (\iota_{gf} \times \looseid_{gf})\looseid_{(I_A \times \transid)i_2} \horc \chi_{gf}\looseid_{i_2}$} (t4);
 \draw[doubleloose] (t0)  to node[above]{\scriptsize $\looseid_{\tens}(\iota_g \times \looseid_{gf})\looseid_{i_2}$} (t1);
  \draw[doubleloose] (t1)  to node[above]{\scriptsize $\chi_g\looseid_{I_A \times f}\looseid_{i_2}$} (t2);
\draw[doubleloose] (t2) to node[above]{\scriptsize $\looseid_{g\tens }(\iota_f \times \looseid_{f})\looseid_{i_2}$} (t3);
  \draw[doubleloose] (t3) to node[above]{\scriptsize $\looseid_g \chi_f \looseid_{(I_A \times \transid)i_2}$}(t4);
  \draw[doubleloose] (t4) to node[above]{\scriptsize $\looseid_{gf}l_I$}(t5);
  \draw[doubleloose] (m0)  to node[above]{\scriptsize $\looseid_{\tens}(\iota_g \times \looseid_{g})\looseid_f$} (m1);
  \draw[doubleloose] (m1)  to node[above]{\scriptsize $\chi_g\looseid_{(I_B \times \transid)i_2 f}$} (m2);
   \draw[doubleloose] (m2) to node[below]{\scriptsize $ \looseid_g l \looseid_f$}(m5); 
    \draw[doubleloose] (b0) to node[above]{\scriptsize $ l \looseid_g \looseid_f$}(b5); 
       \draw[doubleloose] (b0) to[in=220, out=-60] node[above]{\scriptsize $l \looseid_{gf}$}(b5); 
  \draw[doubleeq] (t0) to (m0);
    \draw[doubleeq] (t2) to (m2);
  \draw[doubleeq] (t5) to (m5);
  \draw[doubleeq] (m0) to (b0);
    \draw[doubleeq] (m5) to (b5);
    \node at (.5,1.5) {\footnotesize $=$}; 
   \node at (1.75,1.5) {\footnotesize $\widehat{ \tightid_g \gamma^f}$}; 
   \node at (1.25,.5) {\footnotesize $\widehat{  \gamma^g \tightid_{\looseid}}$}; 
      \node at (1,2.35) {\footnotesize $\iso$}; 
  \node at (1.25,-.35) {\footnotesize $\iso$}; 
 \end{tikzpicture}
 \end{aligned}
\end{equation}
\caption{The structure 3-cell $\gamma^{g\circ f}$}\label{fig:gamma-circ}
\end{figure}

Let $f,h: A \rightarrow B $ and $g,i: B \rightarrow C$ be lax monoidal 1-cells and let $\alpha: f \rightarrow h$, $\beta: g \rightarrow i$ be lax monoidal 2-cells, the composite $\beta \comp \alpha$ is lax monoidal with the structure 3-cells shown in Figure~\ref{fig:comp}.
Colax and strong monoidal structure 2-cells are obtained dually.

\begin{figure}
\begin{equation}
M^{\beta \comp \alpha} := 
\begin{aligned}
 \begin{tikzpicture}[yscale=1.5, xscale=4]
  \node (t0) at (0,-2) {\scriptsize$I_C$};
\node (t2) at (1,-2) {\scriptsize$g f  I_A$};
 \node (t4) at (2,-2) {\scriptsize$i h I_A$};
 \node (m0) at (0,-1) {\scriptsize$I_C$};
 \node (m1) at (.5,-1) {\scriptsize$g I_B$}; 
\node (m2) at (1,-1) {\scriptsize$h I_B$};
\node (m3) at (1.5,-1) {\scriptsize$h f I_A$};
\node (m4) at (2,-1) {\scriptsize$h k I_A$};
 \node (b0) at (0,0) {\scriptsize$I_C$};
 \node (b1) at (.5,0) {\scriptsize$I_C$}; 
\node (b2) at (1,0) {\scriptsize$h I_B$};
\node (b3) at (1.5,0) {\scriptsize$h I_B$};
\node (b4) at (2,0) {\scriptsize$h k I_A$};
\node (bb0) at (0,1) {\scriptsize$I_C$};
 \node(bb2) at (1,1) {\scriptsize$I_C$};
   \node(bb4) at (2,1) {\scriptsize$hk I_A$};
   \draw[doubleloose] (t0)  to node[above]{\scriptsize$\iota_{g f}$} (t2);
  \draw[doubleloose] (t2)  to node[above]{\scriptsize$\beta \alpha$} (t4);
\draw[doubleloose] (m0) to node[above]{\scriptsize$\iota_g $} (m1);
  \draw[doubleloose] (m1) to node[above]{\scriptsize$\beta \looseid_{I}$}(m2);
  \draw[doubleloose] (m2) to node[above]{\scriptsize$\looseid_h \iota_f $}(m3);
  \draw[doubleloose] (m3) to node[above]{\scriptsize$\looseid_h \alpha \looseid_{I}$}(m4);
  \draw[doubleloose] (b0) to node[above]{\scriptsize$\looseid$} (b1);
  \draw[doubleloose] (b1) to node[above]{\scriptsize$\iota_h$} (b2);
  \draw[doubleloose] (b2) to node[above]{\scriptsize$\looseid_h \looseid_{I}$}(b3);
  \draw[doubleloose] (b3) to node[above]{\scriptsize$\looseid_h \iota_k$}(b4);
  \draw[doubleloose] (bb0)  to node[above]{\scriptsize$\looseid_{I}$} (bb2);
  \draw[doubleloose] (bb2)  to node[above]{\scriptsize$\iota_{hk}$} (bb4);
    \draw[doubleeq] (t0) to (m0);  
   \draw[doubleeq] (m0) to (b0);
      \draw[doubleeq] (b0) to (bb0);
    \draw[doubleeq] (t4) to (m4);  
   \draw[doubleeq] (m4) to (b4);
      \draw[doubleeq] (b4) to (bb4);
   \draw[doubleeq] (m2) to (b2);
 \node at (1,.5) {\scriptsize $\iso$}; 
  \node at (.5,-.5) {\scriptsize $\DDownarrow M^{\beta} $}; 
    \node at (1.5,-.5) {\scriptsize $\DDownarrow \widehat{\tightid_{I} M^{\alpha}} $}; 
   \node at (1,-1.5) {\scriptsize $\iso$};  
 \end{tikzpicture}
 \end{aligned}
\end{equation}

\begin{equation}
\Pi^{\beta \comp \alpha} := 
\begin{aligned}
  \begin{tikzpicture}[yscale=1.5, xscale=5]
 \node (t0) at (0,-2) {\scriptsize $\tens (gf \times gf)$};
\node (t2) at (1,-2) {\scriptsize$gf \tens $};
 \node (t4) at (2,-2) {\scriptsize$hk \tens $};
 \node (m0) at (0,-1) {\scriptsize$\tens (gf \times gf)$};
 \node (m1) at (.5,-1) {\scriptsize$g \tens (f\times f)$}; 
\node (m2) at (1,-1) {\scriptsize$h \tens (f \times f)$};
\node (m3) at (1.5,-1) {\scriptsize$hf \tens $};
\node (m4) at (2,-1) {\scriptsize$hk \tens $};
 \node (b0) at (0,0) {\scriptsize$\tens (gf  \times gf)$};
 \node (b1) at (.5,0) {\scriptsize$\tens (hf \times hf)$}; 
\node (b2) at (1,0) {\scriptsize$h \tens (f \times f)$};
\node (b3) at (1.5,0) {\scriptsize$h \tens (k \times k)$};
\node (b4) at (2,0) {\scriptsize$hk \tens $};
\node (bb0) at (0,1) {\scriptsize$\tens (gf \times gf)$};
 \node(bb2) at (1,1) {\scriptsize$\tens (hk \times hk)$};
   \node(bb4) at (2,1) {\scriptsize$hk \tens $};
 \draw[doubleloose] (t0)  to node[above]{\scriptsize$\chi_{gf} $} (t2);
  \draw[doubleloose] (t2)  to node[above]{\scriptsize$\beta \alpha \looseid_{\tens}$} (t4);
\draw[doubleloose] (m0) to node[above]{\scriptsize$\chi_g \looseid_{f \times f}$} (m1);
  \draw[doubleloose] (m1) to node[above]{\scriptsize$\beta \looseid_{\tens (f \times f)}$}(m2);
  \draw[doubleloose] (m2) to node[above]{\scriptsize$\looseid_{h} \chi_f$}(m3);
  \draw[doubleloose] (m3) to node[above]{\scriptsize$\looseid_h \alpha \looseid_{\tens}$}(m4);
  \draw[doubleloose] (b0) to node[above]{\scriptsize$\looseid_{\tens} (\beta \times \beta) \looseid_{f \times f}$} (b1);
  \draw[doubleloose] (b1) to node[above]{\scriptsize$\chi_h \looseid_{f \times f}$} (b2);
  \draw[doubleloose] (b2) to node[above]{\scriptsize$\looseid_{h \tens} (\alpha \times \alpha) $}(b3);
  \draw[doubleloose] (b3) to node[above]{\scriptsize$\looseid_h \chi_k $}(b4);
  \draw[doubleloose] (bb0)  to node[above]{\scriptsize$\looseid_{\tens} (\beta \alpha \times \beta \alpha) $} (bb2);
  \draw[doubleloose] (bb2)  to node[above]{\scriptsize$\chi_{hk} $} (bb4); 
  \draw[doubleeq] (t0) to (m0);  
   \draw[doubleeq] (m0) to (b0);
      \draw[doubleeq] (b0) to (bb0);
    \draw[doubleeq] (t4) to (m4);  
   \draw[doubleeq] (m4) to (b4);
      \draw[doubleeq] (b4) to (bb4);
   \draw[doubleeq] (m2) to (b2);
 \node at (1,.5) {\scriptsize$\iso$}; 
  \node at (.5,-.5) {\scriptsize$\DDownarrow \widehat{\Pi_{lax}^{\beta} \tightid_{\looseid_{f \times f}}}$}; 
    \node at (1.5,-.5) {\scriptsize $\DDownarrow \widehat{\tightid_{\looseid_{h}} \Pi_{lax}^{\alpha} }$}; 
   \node at (1,-1.5) {\scriptsize $\iso$}; 
 \end{tikzpicture}
 \end{aligned}
\end{equation}
\caption{Monoidal structure on $\beta \comp \alpha$}\label{fig:comp}
\end{figure}

Let $f,h: A \rightarrow B $ and $g,i: B \rightarrow C$ be lax monoidal 1-cells and let $\alpha: f \rightarrow h$, $\beta: g \rightarrow i$ be monoidal icons, the composite $\beta \comp \alpha$ is monoidal with $N^{\beta \comp \alpha}:= N^{\beta} \horc \looseid_{\beta} N^{\alpha}$ and $\Sigma^{\beta \comp \alpha}:= \Sigma^{\beta}\looseid_{\alpha \times \alpha} \horc \looseid_{\beta} \Sigma^{\alpha}$.

When $g,f$ are strong monoidal, the maps $\chi_{g \comp f}$ and $\overline{\chi_{g \comp f}}$ are an adjoint equivalence, constructed by the enriched composition from the adjoint equivalence of the pairs $\chi_g, \overline{\chi_g}$ and $\chi_f, \overline{\chi_f}$. Similarly, $\iota_{g \comp f}$ and $\overline{\iota_{g \comp f}}$ form an adjoint equivalence. One can check that the required pairs of 3-cells correspond to eachother as mates by componentwise application of the adjoint equivalences for the composites of $\iota_{g \comp f}$ and $\chi_{g \comp f}$.

In all coherence equations between 3-cells for the monoidal and braided, sylleptic and symmetric structure of composition above, each 3-cell consists of a component for the first composite  composed with the identity on the second composite, and a component for the second composite composed with the identity on the product of the first composite with itself. This means that the coherence equations for $g \comp f$  can be established by componentwise application of the equations for $g$ and $f$. Some 3-cells also contain coherence cells, but these equally break up in a part concerning the first, and a part concerning the second component. Manipulation of these coherence cells results in the required equalities. Note that rewriting the 1-cells and composites of loose 2-cells is necessary in several of the steps. A similar argument holds for coherence equations for braided, sylleptic and symmetric cells.

Let $\Gamma$ and $\Delta$ be two composable monoidal 3-cells. It is easy to see that the composition $\Gamma \comp \Delta$ satisfies the two equations for monoidal 3-cells. This is a matter of applying the equations for $\Gamma$ and $\Delta$ sequentially.
\end{proof}

\begin{rmk}
  Although $\fB$ is 1-strict, $\cMon_{ll}\fB$ is not.
  In particular, our construction cannot be ``iterated'' without generalizing our definitions.
  However, if such a generalization were made, it would be natural to conjecture that $\cMon_{pp}\cMon_{pp}\fB \simeq \cBr\cMon_{pp}\fB$, analogously to~\cite[Remark 5.1]{js:braided}, and similarly $\cMon_{pp}\cMon_{pp}\cMon_{pp}\fB \simeq \cSyl\cMon_{pp}\fB$ and $\cMon_{pp}\cMon_{pp}\cMon_{pp}\cMon_{pp}\fB \simeq \cSym\cMon_{pp}\fB$.
  Similarly, a monoidal object in $\cMon_{ll}\fBicat$ ought to yield a notion of ``duoidal bicategory'' categorifying~\cite{am:mfsha,bm:duoidal}.
\end{rmk}

 \begin{lem}\label{lem:funcmonob}
   Let $F: T \rightarrow S$ be a 1-strict functor of 1-strict locally cubical bicategories with products.
   If $F$ preserves products, it preserves  monoidal objects, 1-cells, 2-cells, icons and 3-cells as well as any braided, sylleptic or symmetric structure on the objects, 1-cells, 2-cells, icons and 3-cells.
 \end{lem}
 
 \begin{proof}
Let $A$ be a monoidal object. As the functor $F$ preserves products, we have a product $F(A) \times_{F} F(A) = F(A \times A)$. As a consequence $\ten\maps
  A\times A\to A$ induces 1-cells $\ten_{F} \maps
 FA\times_{F} F A\to FA$ and $I_{F}:= F(I_A)$. 
 
Since $\phi$ and $\phi_u$ are globular, we have an equality $F(f \comp g) = F(f) \comp F(g)$ for all $f$ and $g$ and for the identity 1-cell we have an equality $F(\transid_A) = \transid_{F(A)}$. The loose associativity 2-cell of $A$ gives rise to a loose 2-cell
  \[\vcenter{\xymatrix@C=6pc{F(A)\times F(A)\times F(A) \rtwocell^{\ten_{F}
        (\Id\times\ten_{F})}_{\ten_{F}(\ten_{F}\times\Id)}{\hspace{.2cm}\fa_{F}\eqv} &F(A) }}\]
  which simply equals $F(\alpha)$ together with the invertible 2-cells.
  
  Likewise, the unit constraints $l, r$ as well as the constraints for (braided) monoidal 1-cells $\sigma$ induce 1-cells $l_{F}, r_{F}$, and $\sigma_{F}$, respectively. Note that the swap functor $\tau$ is mapped by $F$ to the swap functor for the product $\times_{F}$, so $\sigma_{F}$ is well-defined.
  
  Furthermore, the invertible 3-cell filling the Mac Lane pentagon lifts to the invertible 3-cell of the Mac Lane pentagon for $F(A)$ shown in Figure~\ref{fig:F-pent}, which is simply its image under $F$ composed the natural transformations $F_{\odot}$, $\phi$, and $\phi_u$ to ensure that it has the right type.
   Note that there may be several way to paste these 3-cells, but by coherence of enriched pseudo functors, the result is the same. Likewise, the invertible 3-cells $\mu, \lambda,\rho$, as well as $R,S$, and $v$ witnessing the braiding and syllepsis, lift to the appropriate 3-cells for $F(A)$.
   Finally,  the three equations between pasting composites of $\pi_{F}, \mu_{F}, \lambda_{F}, \rho_{F}$ hold by coherence of enriched pseudo functors.

  
\begin{figure}
\begin{tikzpicture}[yscale=1.5, xscale=3]
\node(tl) at (0,1) {$\ten (\ten \times \transid)(\ten \times \transid \times \transid)$};
\node(t) at (1.5,2) {$\ten (\ten \times \transid)(\transid \times \ten \times \transid)$};
\node(tr) at (3,1) {$\ten (\transid \times \ten )(\transid \times \ten \times \transid)$};
\node(br) at (3,0) {$\ten (\transid \times \ten )(\transid \times \transid \times \ten )$};
\node(b) at (1.5,-1) {$\ten (\ten \times \transid)(\transid \times \transid \times \ten )$};
\node(bl) at (0,0) {$\ten (\transid \times \ten )(\ten \times \transid \times \transid)$};
\draw[->] (tl) to node[left, yshift=1pt] {$\looseid (\alpha \times \looseid)$} (t);
\draw[->] (t) to node[right, yshift=1pt] {$\alpha \looseid$} (tr);
\draw[->] (tr) to node[right] {$\looseid (\looseid \times \alpha)$} (br);
\draw[->] (tl) to node[left] {$\alpha \looseid$} (bl);
\draw[->] (bl) to node[left,yshift=-1pt] {$\looseid$} (b);
\draw[->] (b) to node[right,yshift=-1pt] {$\alpha \looseid$} (br);
\draw[->] (tl) to [in=155, out=5] (br);
\draw[->] (tl) to [in=180, out=-10] (br);
\node at (1.5,.6) {$\Downarrow F(\pi) \iso$};
\node at (2,1) {$\iso$};
\node at (1,0) {$\iso $};
\end{tikzpicture}  
\caption{The pentagonator for $F(A)$}\label{fig:F-pent}
\end{figure}


        

Similarly, one can show that for a monoidal transformation $f$, $F(f)$ is monoidal with structure loose 2-cells $F(\iota_f)$ and $F(\chi_f)$ and 3-cells obtained from $F(\omega),  F(\gamma)$, and $F(\delta)$ analogously to $\alpha$. If $f$ is braided, the braiding of $F(f)$ is witnessed by the 3-cell obtained from $F(u)$. Likewise, $F$ preserves any monoidal, braided, sylleptic and symmetric structure of tight 2-morphisms, loose 2-cells and 3-cells. 
 \end{proof}

\begin{thm}\label{thm:lcbcfunctor}
  Let $F: \fB  \rightarrow \fC$ be a 1-strict functor between 1-strict locally cubical bicategories with products.
  If $F$ preserves products, then it lifts to the functors below between locally cubical bicategories, for $v,w \in \{l,c,p\}$.
\begin{align*}
\cM on_{vw}F &: \cM on_{vw} \fB  \rightarrow \cM on_{vw}\fC\\ 
\cBr \cM on_{vw} F &: \cBr \cM on_{vw} \fB  \rightarrow \cBr \cM on_{vw}\fC\\
 \cSyl \cM on_{vw}F &: \cSyl \cM on_{vw} \fB  \rightarrow \cSyl \cM on_{vw}\fC\\
\cSym \cM on_{vw}F &: \cSym \cM on_{vw} \fB  \rightarrow \cSym \cM on_{vw}\fC
\end{align*}
\end{thm}

\begin{proof}
By Lemma~\ref{lem:funcmonob}, the assignment of $F$ on objects and higher cells is well-defined in the respective categories.
We wil show that functor $F$ gives rise to a pseudo double functor $\cM on_lF: \cM on_l \fB(A, B)  \rightarrow \cM on_l\fC(\cM on_lF(A), \cM on_lF(B))$. The other functors are derived in a similar way. We need to verify that $F(N^{\alpha \verc \beta}) = F(N^{\alpha}) \verc F(N^{\beta})$ and $F(N^{\tightid_f}) = N^{\tightid_{Ff}}$ in $\cM on_l\fC$ and likewise for $\Sigma$. This follows from functoriality of $F$, the fact that $N^{\alpha \verc \beta} = N^{\alpha \verc \beta}$ and $N^{\tightid_f} = \tightid_{\iota_f}$, $F(\iota_f) = \iota_{Ff}$.
The natural transformations $F_{\odot}$ and $F_U$ are well-defined 3-cells in $\cM onD$; the respective equations hold by coherence of the pseudo double functor $F$.
Finally, we need to prove that $\phi$ and $\phi_u$ are well-defined tight transformations in $\cM on \fC(F A, F B)$. Since their components are globular 3-cells, we only need to check that these are monoidal. The equations hold by coherence of enriched functors after expanding the definitions of $M^{F\alpha \comp F \beta}, M^{F(\alpha \comp \beta)}, M^{F(\looseid_f)}$, and $M^{\looseid_{F f}}$ and similarly for $\Pi$.
\end{proof}

\begin{thm}\label{thm:H}
  The functor $\cH: \fDblf \rightarrow \fBicat$ lifts to the functors
  \begin{align*}
\cM on_{vw}\cH &: \cM on_{vw} \fDblf  \rightarrow \cM on_{vw}\fBicat\\ 
\cBr \cM on_{vw} \cH &: \cBr \cM on_{vw} \fDblf  \rightarrow \cBr \cM on_{vw}\fBicat\\
 \cSyl \cM on_{vw}\cH &: \cSyl \cM on_{vw} \fDblf  \rightarrow \cSyl \cM on_{vw}\fBicat\\
\cSym \cM on_{vw}\cH &: \cSym \cM on_{vw} \fDblf  \rightarrow \cSym \cM on_{vw}\fBicat.
\end{align*}
\end{thm}

\begin{proof}
  Since $\cH$ is 1-strict, the result then follows from Theorem~\ref{thm:lcbcfunctor}.
\end{proof}

\begin{rmk}\label{rmk:loosely-strong}
  Recall that the notation $\fDblf$ means that we restrict the 2-cells to have \emph{loosely strong} companions.
  As noted in \autoref{thm:iso-strong}, the ``loosely strong'' restriction is automatic for \emph{invertible} 2-cells, so it does not appear at the level of monoidal objects: any monoidal double category with companions has a monoidal loose bicategory.
  However, in order for a lax monoidal double functor $F:\lD\to\lE$ to induce a monoidal functor of loose bicategories, it does not suffice for $\lD$ and $\lE$ to have companions: the constraints $FA \tens FB \to F(A\tens B)$ and $I_\lE \to F(I_\lD)$ must induce pseudonatural, rather than colax natural, transformations.
  But, again, for a \emph{strong} monoidal double functor, this extra condition is automatic.
\end{rmk}




\section{Examples}\label{sec:Alg} 

We now consider a number of examples, to illustrate the utility of our results.

\subsection{Monoids and bimodules in double categories}
\label{sec:mod}

To start with, we apply Theorem~\ref{thm:H} to prove that the family of bicategories $\cAlg[{\mathbb{D}}]$ of monoids, bimodules and bimodule homomorphisms in a monoidal double category $\mathbb{D}$ is monoidal. The result of this section builds on work of the first author in~\cite[Theorem 11.5]{shulman:frbi} and some extentions of this theorem in the second author's PhD thesis~\cite[Chapter 5]{westerPhDthesis}.

In addition to the well-known bicategory $\cMod$ of rings and bimodules,
the family $\cAlg[{\mathbb{D}}]$ contains the well-known bicategories $2\cVect$, $2\cHilb$ introduced in~\cite{kv:2cat-zam,baez2004higher} and also the bicategory $\cProf$ of categories and profunctors~\cite{benabou}. Additionally, it contains various subcategories that are relevant in the field of categorical quantum mechanics, such as the equivariant completion of a braided monoidal bicategory~\cite{carquevillerunkel}, which is a tool for finding topological quantum field theories, as well as the bicategory $2[CP^*[{\bf C}]]$ defined in~\cite{heunenvicarywester} as the mathematical foundation of a diagrammatic language for quantum protocols. 

A {\bf monoid} $(a, A, \mult, \unit)$ in a monoidal double category $\mathbb{D}$ consists of an object $a$, a loose 1-cell $A: a \mapsto a$ and globular 2-cells $\mult: A \odot A \rightarrow A$, $\unit: I_a \rightarrow A$, such that $\mult \circ (\mult \odot \id) = \mult \circ (\id \odot \mult)$ and $\mult \circ (\unit \odot \id) = \id = \mult \circ (\id \odot \mult)$.
In other words, it is a monoid in the usual sense in the monoidal category $\mathbb{D}(a,a)$. A {\bf monoid homomorphism} $(A,\tinymult[gray dot], \tinyunit[gray dot]) \mapsto (B,\tinymult[black dot], \tinyunit[black dot])$ is a pair $(f, {\bf f})$ of a 1-morphism $f:a\rightarrow b$ in $\mathbb{D}$ and a 2-cell ${\bf f}:A \rightarrow B$ in $\mathbb{D}$ that respects the multiplication, $f \circ \tinymult[gray dot] = \tinymult[black dot] \circ (f \odot f)$, as well as the unit, $\tinyunit[black dot] = f \circ \tinyunit[gray dot]$. $f: A \rightarrow B$. A {\bf bimodule} $(a,A,\tinymult[gray dot], \tinyunit[gray dot]) \mapsto (b, B,\tinymult[black dot], \tinyunit[black dot])$ is a pair (${\bf M}, M)$ of a $1-cell$ $M$ in $\doub{D}$ and a globular 2-cell ${\bf M}:A \times M \times B \rightarrow M$ in $\doub{D}$, with the structure of an $A$-$B$- bimodule. We will simply write $M$ for the bimodule $({\bf M}, m)$. This is also called an {\bf A-B-bimodule}. Note that $S(M) = T(A)$ and $T(M) = S(B)$. Let $\phi: (a,A,\tinymult[gray dot], \tinyunit[gray dot]) \rightarrow (c,C,\tinymult[black dot], \tinyunit[black dot])$ and $\psi: (b,B,\tinymult[gray dot], \tinyunit[gray dot]) \rightarrow (d,D,\tinymult[black dot], \tinyunit[black dot])$ be monoid homomorphisms and let $M$ and $N$ be an $A$-$B$-bimodule and a $C$-$D$-bimodule, respectively. A {\bf $(\phi, \psi)$-equivariant map} is a 2-morphism $f:M \rightarrow N$ in ${\doub{D}}$ such that ${\bf M} \circ (\phi \tens f \tens \psi) = f \circ {\bf N}$.

In other sources~\cite{westerthesis}~\cite{jamiesbook}~\cite{heunenvicarywester},~\cite{carquevillerunkel},  monoids are called algebras and equivariant maps correspond to extended bimodule homomorphisms.


\begin{defn}
Let ${\mathbb{D}}$ be a double category. We can define a new double category $\lAlg[{\mathbb{D}}]$ consisting of the elements listed below.

\begin{itemize}
\item 0-cells are {\bf monoids} $(a, A,\tinymult[gray dot], \tinyunit[gray dot])$, in the monoidal category ${\mathbb{D}_1}$. 
\item tight 1-cells $(A,\tinymult[gray dot], \tinyunit[gray dot]) \mapsto (B,\tinymult[black dot], \tinyunit[black dot])$ are {\bf monoid homomorphism}
\item loose 1-cells $(A,\tinymult[gray dot], \tinyunit[gray dot]) \mapsto (B,\tinymult[black dot], \tinyunit[black dot])$ are {\bf bimodules}
\item 2-cells $({\bf M},M) \rightarrow ({\bf N},N)$ from a $A$-$B$-bimodule to a $C$-$D$-bimodule are $(\phi, \psi)$-equivariant maps, where $\phi:  A\rightarrow C$, $\psi: B \rightarrow D$ are monoid homomorphisms.  
\end{itemize}
Structural data regarding this construction, such as loose composition, is described in~\cite{shulman:frbi}, where the double category is called $\mathbb{M}od(\mathbb{D})$. A more detailed description in the case that $\mathbb{D}$ is a monoidal category is given in~\cite{westerPhDthesis}.
\end{defn}

For the double category to be well-defined, we need certain coequalizers in $\mathbb{D}$ to exist. To this end we recall the definition below from~\cite[Definition 11.4]{shulman:frbi}.

\begin{defn}
A double category $\mathbb{D}$ has {\bf local coequalizers} if for each objects $a,b$, the hom-category $\mathbb{D}(a,b)$ has all coequalizers and the coequalizers are preserved by $\odot$. We write $\cDbl^l$ for the 2-category of double categories with local coequalizers and $\cDblf^l$ for the 2-category of double categories with local coequalizers with 2-cells that have loosely strong companions.
\end{defn}

\begin{prop}\label{thm:eqcomp}
Let $\mathbb{D}$ be a monoidal double category with local coequalizers and companions and conjoints, and such that the tensor product $\tens$ preserves the local coequalizers. The loose bicategory $\cAlg[\mathbb{D}]$ of $\lAlg[\mathbb{D}]$ is monoidal; it is braided or symmetric whenever $\mathbb{D}$ is braided or symmetric.
\end{prop}

\begin{proof}
By~\cite[Examples 9.2]{shulman:frbi}, the double category $\cAlg[\mathbb{D}]$ is symmetric monoidal and has companions and conjoints. Therefore, we may apply Theorem~\ref{thm:lcbcfunctor}, which states that its loose bicategory $\cAlg[\mathbb{D}]$ is symmetric monoidal.
\end{proof}

An explicit description of the monoidal structure for the special case when ${\mathbb{D}}$ is a monoidal category was given in~\cite{westerPhDthesis}. 

\begin{cor}
The bicategory $2\cVect$ is symmetric monoidal.
\end{cor}

\begin{proof}
We obtain $2\cVect$ as $\cAlg[\cat{FVect}]$ from the braided monoidal category $\cat{FVect}$ of finite dimensional vector spaces and linear maps. The category $\cat{FVect}$ is symmetric monoidal, contains local coequalisers and has a tensor product that preserves coequalisers. The result followis from Proposition~\ref{prop:funcAlg}.
\end{proof}

\begin{cor}
The bicategory $\cProf(\cat{C})$ of internal categories and profunctors in a category $\cat{C}$ with finite limits and coequalizers preserved by pullback is symmetric monoidal.
\end{cor}

\begin{proof}
The bicategory $\cProf(\cat{C})$ can be constructed as the category $\cAlg[\lSpan[\cat{C}]]$ from the monoidal double category of spans~\cite[Examples 4.2]{shulman:frbi}.
This category has local coequalisers and companions and is symmetric monoidal~\cite[Examples 4.15, 9.2]{shulman:frbi}. The result follows from Proposition~\ref{prop:funcAlg}.
\end{proof}

\begin{cor}
The bicategory $\cat{V}\cProf$ of categories and profunctors enriched over a cocomplete closed symmetric monoidal category $\cat{V}$ is symmetric monoidal.
\end{cor}

\begin{proof}
The bicategory $\cat{V}\cProf$ can be constructed as the category $\cAlg[\lMat[\cat{V}]]$ from the monoidal double category of matrices~\cite[Examples 11.8]{shulman:frbi}.
This category has local coequalisers and companions and is symmetric monoidal, so the result again follows from Proposition~\ref{prop:funcAlg}.
\end{proof}

Similarly, monoidal functors $F,G:\mathbb{D} \rightarrow \mathbb{E}$ between double categories with companions lift to monoidal functors of the form $\lAlg F,\lAlg G: \cAlg[\mathbb{D}] \rightarrow \cAlg[\mathbb{E}]$ and monoidal transformations $\alpha: F \Rightarrow G$ lift to monoidal transformations of the form $\lAlg\alpha: \lAlg F \Rightarrow \lAlg G$, as long as all the constraints are loosely strong (which in \autoref{rmk:loosely-strong} we noted is automatic when they are invertible). In fact, the $\lAlg$ construction gives rise to a functor, as shown below. 

\begin{prop}\label{prop:funcAlg}
The $\bAlg$ construction gives rise to the following functors of bicategories:
\begin{align*}
\mathcal{M}on_{pp}\cDblf^l &\rightarrow \mathcal{M}on\mathcal{B}icat\\
\mathcal{B}r\mathcal{M}on_{pp}\cDblf^l &\rightarrow \mathcal{B}r\mathcal{M}on\mathcal{B}icat\\
\mathcal{S}ym\mathcal{M}on_{pp}\cDblf^l &\rightarrow \mathcal{S}ym\mathcal{M}on\mathcal{B}icat.
\end{align*}
\end{prop}
\begin{proof}
By~\cite[Proposition 11.22]{shulman:frbi}, $\lAlg$ is a functor $\mathcal{M}on_{pp}\cDblf^l \rightarrow \mathcal{M}on_{pp} \cDblf^l$, and similarly for the braided and symmetric versions. We compose this with the functor from Theorem~\ref{thm:H} to obtain the result.
\end{proof}

\begin{eg}\label{thm:span-functor}
  If $F:\cat{C}\to\cat{D}$ is a functor between categories with finite limits and coequalizers preserved by pullback, and $F$ preserves finite limits and coequalizers, then it induces a strong symmetric monoidal double functor $\lSpan(\cat{C}) \to \lSpan(\cat{D})$, hence $\lAlg[\lSpan(\cat{C})] \to \lAlg[\lSpan(\cat{D})]$, and thus a symmetric monoidal functor of bicategories $\cProf(\cat{C}) \to \cProf(\cat{D})$.
\end{eg}

\begin{eg}\label{thm:cat-functor}
  If $F:\cat{V}\to\cat{W}$ is a cocontinuous strong monoidal functor between closed symmetic monoidal categories, it induces a strong symmetric monoidal double functor $\lMat(\cat{V})\to\lMat(\cat{W})$, hence $\lAlg[\lMat(\cat{V})]\to\lAlg[\lMat(\cat{W})]$, and thus a symmetric monoidal functor of bicategories $\cat{V}\cProf \to \cat{W}\cProf$.
\end{eg}

Similar results apply to the ``enriched indexed categories'' of~\cite{shulman:eicats}.

\subsection{Applications in Quantum Theory}
Frobenius algebras and modules play an important role in quantum theory. 

\begin{defn}
A {\bf Frobenius algebra} in a monoidal double category ${\cat{C}}$ is a monoid $(A, \tinymult, \tinyunit)$ together with a comonoid 
$(a, A, \tinycomult, \tinycounit)$ that satisfies the equation $
(\id \tens \tinymult[gray dot]) \circ \alpha \circ (\tinycomult[gray dot] \tens \id) = (\tinymult[gray dot] \tens \id) \circ \alpha^{-1} \circ (\id \tens \tinycomult[gray dot] )$.
A monoid in a braided monoidal category {\cat C} is called {\bf commutative} when $
\tinymult[gray dot] \circ \sigma = \tinymult[gray dot]$;
it is {\bf symmetric} if the weaker condition 
$\tinycounit[gray dot] \circ \tinymult[gray dot] \circ \sigma = \tinycounit[gray dot] \circ \tinymult[gray dot]
$ holds.
A pair of a monoid  $\tinymult[gray dot]$ and a comonoid $\tinycomult[black dot]$  is called {\bf special} when the equation $ \tinymult[gray dot] \circ \tinycomult[black dot] = \id$ holds.
\end{defn}

Recently, the bicategory of Frobenius algebras, bimodules and bimodule homomorphisms in a monoidal bicategory $\cat{B}$ was introduced in~\cite{carquevillerunkel} as the {\it equivariant completion }of $\cat{B}$. This is a tool for generating topological quantum field theories, which is a non-trivial process in general.

 
 
 



Let $\doub{D}$ be a double category. The {\bf equivariant completion} of the loose bicategory $\cH(\doub{D})$ is the loose bicategory of the double category $\mathbb{E}q[{\mathbb{D}}]$ defined below.  

\begin{defn}
Let ${\mathbb{D}}$ be a double category with local coequalisers and loosely strong companions, where the tensor product preserves coequalisers. We can define a new double category $\mathbb{E}q[{\mathbb{D}}]$ consisting of the elements listed below.

\begin{itemize}
\item 0-cells are {\bf Frobenius algebras} $(a, A,\tinymult[gray dot], \tinyunit[gray dot])$, in the monoidal category ${\mathbb{D}_1}$. 
\item tight 1-cells $(A,\tinymult[gray dot], \tinyunit[gray dot]) \mapsto (B,\tinymult[black dot], \tinyunit[black dot])$ are {\bf monoid homomorphism}
\item loose 1-cells $(A,\tinymult[gray dot], \tinyunit[gray dot]) \mapsto (B,\tinymult[black dot], \tinyunit[black dot])$ are {\bf bimodules}
\item 2-cells $({\bf M},M) \rightarrow ({\bf N},N)$ from a $A$-$B$-bimodule to a $C$-$D$-bimodule are $(\phi, \psi)$-equivariant maps, where $\phi:  A\rightarrow C$, $\psi: B \rightarrow D$ are monoid homomorphisms.  
\end{itemize}

Structural data regarding this construction, such as loose composition, is as for $\lAlg[{\mathbb{D}}]$.
\end{defn}

\begin{cor}
The equivariant completion of a monoidal bicategory $\cat{B}$ is monoidal. It is braided whenever $\cat{B}$ is braided and it is symmetric whenever $\cat{B}$ is symmetric.
\end{cor}

\begin{proof}
This follows from Proposition~\ref{thm:eqcomp} applied to the double category $\mathbb{E}q[\cat{B}]$, where \cat{B} is regarded as a double category with trivial tight 1-cells.
\end{proof}

Another related example in quantum theory is the bicategory of $\cAlg[\mathbb{D}]$ special dagger Frobenius algebras, bimodules and bimodule homomorphisms in a monoidal bicategory~\cite{heunenvicarywester}, which was introduced as a mathematical foundation for a diagrammatic language of quantum protocols. Examples are the well-known bicategory $\cat{2Hilb}$ and the bicategory $2[CP[{\bf C}]]$ of mixed quantum states. The monoidal structure is essential for such applications, as it enables the description of compound quantum systems, as well as parallel quantum protocols.

\begin{defn}
Let ${\bf C}$ be a category. A dagger $\dagger: {\cat C} \rightarrow {\cat C}$, is a contravariant functor which is the identity on objects such that $\dagger(\dagger(f)) = f$. 
A {\bf \index{dagger monoidal category}dagger monoidal category} is a monoidal category that is equipped with a dagger $\dagger$, such that the equalities below hold.
 \begin{mathpar}
 (f \otimes g)^{\dagger} = g^{\dagger} \otimes f^{\dagger}\and
 \alpha^{\dagger} = \alpha^{-1} \and
  \rho^{\dagger} = \rho^{-1} \and
   \lambda^{\dagger} = \lambda^{-1} 
 \end{mathpar}
A {\bf \index{dagger braided monoidal category}dagger braided monoidal category} is a dagger monoidal category with a braiding that satisfies the equality below.
 \begin{equation*}
    \sigma^{\dagger} = \sigma^{-1} \\
 \end{equation*}
 A Frobenius algebra in a dagger braided monoidal category $\cat{C}$ is a {\bf dagger Frobenius algebra} when the comonoid is the dagger image of the monoid. 
A bimodule is called a {\bf dagger bimodule} when the equation below holds, where the comonoid is the dagger of the monoid and we denote $\dagger({\bf M})$ by ${\bf M}^{\dagger}$.
\begin{equation}
{\bf M}^{\dagger} = (\id \tens {\bf M }\tens \id) \circ (\tinycomult[white dot] \tens \id \tens \tinycomult[black dot]) \circ (\tinyunit[white dot] \tens \id \tens \tinyunit[black dot])
\end{equation}
\end{defn}

\begin{defn}
Let $f: A\rightarrow B$ be an algebra homomorphism between dagger Frobenius algebras in a dagger category. The {\bf conjugate} $f_*$ of $f$ is defined as $f_* := (\id \otimes \tinycounit[gray dot]) \circ (\id \otimes \tinymult[gray dot]) \circ (\id \otimes f^{\dagger} \otimes \id) \circ (\tinycomult[gray dot] \otimes \id) \circ (\tinyunit[gray dot] \otimes \id)$
An algebra homomorphism $f$ is {\bf self-conjugate} if $f=f_*$.
\end{defn}

To prove that the $2[-]$ construction preserves braided monoidal structure, one could define a suitable notion of dagger double category and prove that Proposition~\ref{thm:eqcomp} can be specialised to the case where monoids are special dagger Frobenius algebras and bimodules are dagger bimodules. Instead, we make use of the direct proof of the fibrant and monoidal structure of the double category $\ltwo[\cat{C}]$ defined in~\cite{westerPhDthesis}, for the special case that $\cat{C}$ is a monoidal dagger category,~\cite[Prop 5.4.25]{westerPhDthesis}. This double category consists of special dagger Frobenius algebras, self-conjugate algebra homomorphisms, dagger bimodules and equivariant maps.

\begin{thm}
The assignment $2[-]$ gives rise to the functors of locally cubical bicategories below,
\begin{align*}
 \mathcal{B}r\mathcal{M}on\mathcal{C}at \rightarrow \mathcal{B}r\mathcal{M}on\mathcal{B}icat\\
 \mathcal{S}ym\mathcal{M}on\mathcal{C}at \rightarrow \mathcal{S}ym\mathcal{M}on\mathcal{B}icat
\end{align*}
Where $\mathcal{B}r\mathcal{M}on\mathcal{C}at$ and $\mathcal{S}ym\mathcal{M}on\mathcal{C}at$ are the locally cubical bicategories of braided monoidal bicategories and symmetric monoidal bicategories, respectively, where all coequalizers exist and the tensor product preserves coequalizers. These locally cubical bicategories have only trivial loose 2-cells and 3-cells.
\end{thm}

\begin{proof}
This follows directly from~\cite[Proposition 5.4.25]{westerPhDthesis} and Theorem~\ref{thm:eqcomp}.
\end{proof}

Similarly, one can show that the result holds for the bicategory of commutative or symmetric dagger Frobenius algebras, dagger bimodules and bimodule homomorphisms.

\begin{cor}
The bicategory $2\cHilb$ is symmetric monoidal.
\end{cor}

\begin{proof}
As shown in~\cite[Section 3.6.3]{westerthesis}, $2\cHilb$ is equivalent to $2[\cat{FHilb}]$, where $\cat{FHilb}$ is the symmetric monoidal category of finite Hilbert spaces and linear maps, which contains all coequalizers.
\end{proof}

\subsection{Black-boxing of open Markov processes}
\label{sec:markov}

The authors of~\cite{bc:markov} construct a symmetric monoidal double category $\mathbb{M}\mathsf{ark}$ whose objects are finite sets and whose loose 1-cells are ``open Markov processes''.
The loose composition and the tensor product then give two related ways to put together smaller open Markov processes into larger ones (and eventually into closed ones).
In addition, they constructed a symmetric monoidal \emph{functor}, called \textbf{black-boxing}, from this double category to a double category $\mathbb{L}\mathsf{inRel}$ of linear relations.

The authors of~\cite{bc:markov} also showed that both of these double categories have companions for tight isomorphisms, and then used our results to conclude that their loose bicategories $\mathcal{M}\mathit{ark}$ and $\mathcal{L}\mathit{inRel}$ are symmetric monoidal bicategories.
At the time of writing~\cite{bc:markov}, only the earlier version~\cite{shulman:smbicat} of this paper was available, which constructed monoidal bicategories but not monoidal functors between them; thus, the authors of~\cite{bc:markov} were only able to conjecture that their black-boxing double functor induced a symmetric monoidal functor of bicategories.
However, with Theorem~\ref{thm:lcbcfunctor} now in hand, we can prove their conjecture:

\begin{thm}[{\cite[Conjecture 6.7]{bc:markov}}]
  There exists a strong symmetric monoidal functor of bicategories $\blacksquare : \mathcal{M}\mathit{ark} \to \mathcal{L}\mathit{inRel}$ that maps
  \begin{enumerate}
  \item any finite set $S$ to the vector space $\blacksquare(S) = \mathbb{R}^S \oplus \mathbb{R}^S$,
  \item any open Markov process $S \xto{i} (X,H) \xot{o} T$ to the linear relation
    \[ \blacksquare(S \xto{i} (X,H) \xot{o} T) \subseteq \mathbb{R}^S\oplus \mathbb{R}^S\oplus \mathbb{R}^T\oplus \mathbb{R}^T \]
    consisting of all 4-tuples $(i^*(v),I,o^*(v),O)$ where $v\in \mathbb{R}^X$ is some steady state with inflows $I$ and outflows $O$ (see~\cite[Definition 2.7]{bc:markov}), and
  \item any globular morphism of open Markov processes
    \[
      \begin{tikzcd}
        & {(X,H)} \arrow[dd, "p"] &  \\
        S \arrow[ru, "i_1"] \arrow[rd, "i_1'"'] &  & T \arrow[lu, "o_1"'] \arrow[ld, "o_1'"] \\
        & {(X',H')} & 
      \end{tikzcd}
    \]
    to the inclusion $\blacksquare(X,H) \subseteq \blacksquare(X',H')$.
  \end{enumerate}
\end{thm}
\begin{proof}
  Apply Theorem~\ref{thm:lcbcfunctor} to the strong symmetric monoidal double functor of~\cite[Theorem 5.5]{bc:markov}, noting again that the constraints of any \emph{strong} monoidal functor are automatically loosely strong.
\end{proof}

\subsection{Cartesian double categories and cartesian bicategories}
\label{sec:cartesian}

A \emph{cartesian bicategory}, as defined in~\cite{cw:cart-bicats-i,ckww:cartbicats-ii}, is designed as an axiomatization of properties of bicategories like $\cSpan$ in which the objects have ``cartesian products'', but it is complicated because these cartesian products are not products in the ordinary bicategorical sense (only when restricted to the subcategory of ``maps'', i.e.\ left adjoints).
By contrast, a \emph{cartesian double category}~\cite{aleiferi2018cartesian} is easily defined as a cartesian object in the 2-category $\cDbl$, i.e.\ one whose diagonal $\lD \to \lD\times \lD$ and projection $\lD\to *$ have right adjoints, which also has companions and conjoints.
Our results make it clear that:

\begin{thm}
  If $\lD$ is a cartesian double category, then $\cH(\lD)$ is a symmetric monoidal bicategory.
\end{thm}
\begin{proof}
  A cartesian object in a 2-category is automatically a symmetric monoid object therein, so we can apply Theorem~\ref{thm:lcbcfunctor}.
\end{proof}

One might hope that the functor $\cH$ should preserve the cartesian structure as well, implying that $\cH(\lD)$ is a cartesian object in $\fBicat$.
However, as noted above, a cartesian bicategory does not have cartesian products in the ordinary bicategorical sense, so this must fail.
The reason is instructive: the transformations exhibiting the adjunctions in $\cDbl$ making $\lD$ cartesian do not have \emph{loosely strong} companions; thus their loose lifts are only colax natural.
If our functor $\cH$ could be extended to this case somehow, then it might also imply automatically that the loose bicategory of a cartesian double category is a cartesian bicategory.

\subsection{Monoidal fibrations and parametrized spectra}
\label{sec:mfps}

In~\cite{shulman:frbi} the first author also introduced the notion of \emph{monoidal fibration}, and showed that any sufficiently well-behaved symmetric monoidal fibration gives rise to a symmetric monoidal double category with all companions and conjoints.
We will not recall the definition of monoidal fibration here, but only note that the double categories $\lSpan(\cat{C})$ and $\lMat(\cat{V})$ can be constructed in this way (see~\cite[Examples 15.3 and 15.4]{shulman:frbi}).
Moreover, the construction is functorial, so the functors from Examples~\ref{thm:span-functor} and~\ref{thm:cat-functor} can be obtained from morphisms of monoidal fibrations, as well as their enhancement to the enriched indexed categories of~\cite{shulman:eicats}.

One final family of examples obtained in this way are the bicategories of \emph{parametrized spectra} from~\cite{maysig:pht}.

\begin{thm}
  The bicategory $\mathcal{E}\mathit{x}$ of parametrized spectra from~\cite[Chapter 17]{maysig:pht} is symmetric monoidal.
\end{thm}
\begin{proof}
  By~\cite[Example 14.6]{shulman:frbi}, this bicategory can be constructed as the loose bicategory of the double category obtained from a symmetric monoidal fibration of parametrized spectra.
  (The latter monoidal fibration is mostly constructed in~\cite{maysig:pht}; see~\cite{malkiewich:ex-lowtech} for a minor correction.)
  Thus, by~\cite[Theorem 14.2]{shulman:frbi} and Theorem~\ref{thm:lcbcfunctor}, it is symmetric monoidal.
\end{proof}

This symmetric monoidal structure is observed without details in~\cite[Remark 17.17]{maysig:pht}.
Our stronger functoriality results are also relevant to~\cite{maysig:pht}, for instance:

\begin{thm}
  For a map $f:A\to B$ of $G$-spaces, the functor of bicategories $f^*: G\mathcal{E}\mathit{x}^{\mathit{fib}}_B \to G\mathcal{E}\mathit{x}^{\mathit{fib}}_A$ from~\cite[Proposition 19.3.4]{maysig:pht} is a symmetric monoidal functor.
\end{thm}
\begin{proof}
  It is straightforward to verify (indeed, the proof of~\cite[Proposition 19.3.4]{maysig:pht} essentially does this) that it is induced by a strong morphism of symmetric monoidal fibrations.
  Thus, the latter yields a strong symmetric monoidal double functor, and therefore by Theorem~\ref{thm:lcbcfunctor} a symmetric monoidal functor of bicategories.
\end{proof}

The pushforward $i_!$ of~\cite[Proposition 19.3.1]{maysig:pht}, however, is only a colax functor of bicategories.
Thus, although it does come from a colax functor of double categories, our machinery does not apply to it.
This provides an additional motivation to extend our results to lax and colax functors and transformations.


\appendix
\section{Coherence Equations for Locally-Double Bicategories}
\label{ap:coherence}

We give the coherence diagrams for monoidal objects, 1-cells, 2-cells and icons in a locally cubical bicategory, defined in Section~\ref{sec:mono-objects}. For readability, we sometimes write $S(\alpha)$ and $T(\alpha)$ for the composite 2-cells that form the source and target, respectively, of a 3-cell $\alpha$.
Furthermore, we simplify notation for coherence constraints for double categories, by writing $\iso$ to indicate composites of such coherence cells.

For monoidal objects, 1-cells and 2-cells, the coherence axioms live entirely in the loose bicategory of the hom-double-category, since all the tight 2-cells are identities. Therefore, we write them as
ordinary bicategorical pasting diagrams, omitting the identity
2-cells.  

\subsubsection*{Monoidal Object}

\input{coherence/MonObjEq1} \newpage
{\small
\begin{equation*}\hspace{-2cm}
\begin{tikzpicture}[xscale=2.25, yscale=1.5]
\node (A1) at (-.5,7){$\substack{\tens (\tens \times \transid)\\ (\transid \times \transid \times \transid)}$};
\node (A2) at (0,9){$\substack{\tens (\tens \times \transid)\\ ([\tens(\transid \times I] \times \transid \times \transid )}$};
\node (A3) at (1,9){$\substack{\tens (\tens \times \transid) \\  (\tens \times \transid \times \transid) \\ (\transid \times I \times \transid \times \transid)}$};
\node (A5) at (2.5,10){$\substack{\tens (\tens \times \transid) \\ ( \transid \times \tens \times \transid) \\ (\transid \times I \times \transid \times \transid)}$};
\node (A6) at (3.5,10){$\substack{\tens (\tens \times \transid) \\  (\transid \times [\tens (I \times \transid)] \times \transid)}$};
\node (A7) at (5,9.5){$\substack{\tens (\tens \times \transid) \\ (\transid \times \transid \times \transid)}$};
\node (A75) at (6,8.5){$\substack{\tens (\transid \times \tens) \\ (\transid \times \transid \times \transid)}$};
\node (A8) at (6,5){$\substack{\tens (\transid \times \tens)}$};
\draw[doubleloose] (A1) to node[above, xshift=-20]{$\substack{\looseid \looseid \\ (r^{-1} \times \looseid \times \looseid) }$} (A2);
\draw[doubletighteq] (A2) to  (A3);
\draw[doubleloose] (A3) to node[above, yshift=5pt]{$\substack{\looseid (\alpha \times \looseid) \looseid }$} (A5);
\draw[doubletighteq] (A5) to (A6);
\draw[doubleloose] (A6) to node[above, xshift=8pt]{$\substack{\looseid \looseid \\ (\looseid \times l \times \looseid) }$} (A7);
\draw[doubleloose] (A7) to node[above]{$\substack{\alpha}$} (A75);
\draw[doubletighteq] (A75) to (A8);
\node (B6) at (3,8){$\substack{\tens (\transid \times \tens) \\ ( \transid \times \tens \times \transid) \\ (\transid \times I \times \transid \times \transid)}$};
\node (B7) at (4,8){$\substack{\tens (\transid \times \tens) \\ ( \transid \times [\tens (I \times \transid)] \times \transid)}$};
\draw[doubleloose] (A5) to node[right]{$\substack{\alpha \looseid \looseid  }$} (B6);
\draw[doubletighteq] (B6) to  (B7); 
\draw[doubleloose] (B7) to node[below, xshift=10pt, yshift=-5pt]{$\substack{\looseid \looseid (\looseid \times l \times \looseid) }$} (A75); 
\node (C7) at (3,6){$\substack{\tens (\transid \times \tens) \\ ( \transid \times \transid \times \tens) \\ (\transid \times I \times \transid  \times \transid)}$};
\node (C75) at (4,5){$\substack{\tens (\transid \times [\tens (I \times \transid)]) \\ (\transid \times \transid \times \tens) }$};
\node (C8) at (5,5){$\substack{\tens (\transid \times \transid) \\ (\transid \times \transid \times \tens) }$};
\draw[doubleloose] (B6) to node[right]{$\substack{\looseid (\looseid \times \alpha) \\ \looseid}$} (C7); 
\draw[doubletighteq] (C7) to (C75);
\draw[doubleloose] (C75) to node[above]{$\substack{ \looseid (\looseid \times l ) \\ \looseid}$} (C8); 
\draw[doubletighteq] (C8) to  (A8); 
\node (D4) at (0,4){$\substack{\tens ([\tens  (\transid \times I)] \times \tens) \\ ( \transid \times \transid \times \transid) }$};
\node (D5) at (1,5){$\substack{\tens (\transid \times \tens) \\ ( \tens \times \transid \times \transid) \\ (\transid \times I \times \transid \times \transid) }$};
\node (D6) at (2,5){$\substack{\tens (\tens \times \transid) \\ ( \transid \times \transid \times \tens) \\ (\transid \times I \times \transid \times \transid)}$};
\draw[doubleloose] (A3) to node[left]{$\substack{\alpha \looseid \looseid }$} (D5);
\draw[doubletighteq] (D4) to (D5);
\draw[doubletighteq] (D5) to (D6);
\draw[doubleloose] (D6) to node[above,xshift=-10pt]{$\substack{\alpha \looseid    \looseid}$} (C7); 
\node (E1) at (-.5,3){$\substack{\tens (\transid \times \tens)\\ (\transid \times \transid \times \transid)}$};
\node (E3) at (0,2){$\substack{\tens (\transid \times \tens)}$};
\draw[doubleloose] (A1) to node[left]{$\substack{\alpha}$} (E1);
\draw[doubleloose] (E1) to node[right]{$\substack{ \looseid ( r^{-1} \times \looseid ) \looseid}$}  (D4);
\draw[doubletighteq] (E1) to (E3);
\draw[doubletighteq] (E3) to[out= 0, in=270] (A8); 
\node at (0.5,6.5) {$\substack{\DDownarrow \iso } $};
\node at (5.5,7.5) {$\substack{\DDownarrow \iso } $};
\node at (4.5,6.7) {$\substack{\DDownarrow \tightid \lambda}$};
\node at (4,5.8) {$\substack{\DDownarrow \iso } $};
\node at (2.5,8.5){$\substack{\DDownarrow \iso } $};
\node at (2,7){$\substack{\DDownarrow \pi \tightid}$};
\node at (1.5,6){$\substack{\DDownarrow \iso } $};
\node at (4,8.5){$\substack{\DDownarrow \iso } $};
\node at (2.5,4.5){$\substack{\DDownarrow  \iso}$};
\node at (3.5,3){$\substack{\DDownarrow  \looseid (\looseid \times \mu)}$};
\draw[doubleloose] (A3) to[out= -35, in=115] node[right]{$\substack{S(\pi) \looseid}$} (C7); 
\draw[doubleloose] (A3) to[out= -75, in=180] node[left]{$\substack{T(\pi) \\ \looseid}$} (C7);
\draw[doubleloose] (B6) to[out= -25, in=115] node[above, xshift=10pt]{$\substack{\looseid S(\lambda)}$} (A8); 
\draw[doubleloose] (B6) to[out= -55, in=160] node[below, xshift=-10pt]{$\substack{\looseid T(\lambda) }$} (A8);
\draw[doubleloose] (E1) to[out= 25, in=225] node[above]{$\substack{\looseid (\looseid \times S(\mu))}$} (A8); 
\end{tikzpicture}\hspace{-2cm}
\end{equation*}
\begin{equation}\label{eq:monobjeq2}
=
\end{equation}
\begin{equation*}\hspace{-2cm}
\begin{tikzpicture}[xscale=3, yscale=1]
\node (A1) at (-.5,7){$\substack{\tens (\tens \times \transid)\\ (\transid \times \transid \times \transid)}$};
\node (A2) at (0,9){$\substack{\tens (\tens \times \transid)\\ ([\tens(\transid \times I] \times \transid \times \transid )}$};
\node (A3) at (1,10){$\substack{\tens (\tens \times \transid) \\  (\tens \times \transid \times \transid) \\ (\transid \times I \times \transid \times \transid)}$};
\node (A5) at (2.5,10){$\substack{\tens (\tens \times \transid) \\ ( \transid \times \tens \times \transid) \\ (\transid \times I \times \transid \times \transid)}$};
\node (A6) at (3.5,9){$\substack{\tens (\tens \times \transid) \\  (\transid \times [\tens (I \times \transid)] \times \transid)}$};
\node (A7) at (4,7){$\substack{\tens (\tens \times \transid) \\ (\transid \times \transid \times \transid)}$};
\node (A75) at (4,6){$\substack{\tens (\transid \times \tens) \\ (\transid \times \transid \times \transid)}$};
\node (A8) at (3.5,5){$\substack{\tens (\transid \times \tens)}$};
\draw[doubleloose] (A1) to node[above, xshift=-20]{$\substack{\looseid \looseid \\ (r^{-1} \times \looseid \times \looseid) }$} (A2);
\draw[doubletighteq] (A2) to  (A3);
\draw[doubleloose] (A3) to node[above, yshift=5pt]{$\substack{\looseid (\alpha \times \looseid) \looseid }$} (A5);
\draw[doubletighteq] (A5) to (A6);
\draw[doubleloose] (A6) to node[right, xshift=8pt]{$\substack{\looseid \looseid \\ (\looseid \times l \times \looseid) }$} (A7);
\draw[doubleloose] (A7) to node[right]{$\substack{\alpha}$} (A75);
\draw[doubletighteq] (A75) to (A8);
\node (E1) at (-.5,6){$\substack{\tens (\transid \times \tens)\\ (\transid \times \transid \times \transid)}$};
\node (E3) at (0,5){$\substack{\tens (\transid \times \tens)}$};
\draw[doubleloose] (A1) to node[left]{$\substack{\alpha}$} (E1);
\draw[doubletighteq] (E1) to (E3);
\draw[doubletighteq] (E3) to (A8); 
\node at (1.75,9){$\substack{\DDownarrow  \iso}$};
\node at (1.75,7.5){$\substack{\DDownarrow  \tightid (\mu \times \tightid)}$};
\node at (1.75,6){$\substack{\DDownarrow  \iso}$};
\draw[doubleloose] (A1) to node[above]{$\substack{\looseid (\looseid \times \looseid)}$} (A7); 
\draw[doubleloose] (A1) to[out= 55, in=125] node[above]{$\substack{\looseid (S(\mu) \times \looseid)}$} (A7); 
\end{tikzpicture}\hspace{-2cm}
\end{equation*}} \newpage
{\small
\begin{equation*}\hspace{-2cm}
\begin{aligned}
\begin{tikzpicture}[xscale=2.75, yscale=2]
\node (A1) at (0,5){$\substack{\tens (\tens \times \transid)}$};
\node (A3) at (0,8){$\substack{\tens (\transid \times \tens) }$};
\node (A35) at (0.5,9){$\substack{\tens (\transid \times \tens) \ (\transid \times \transid \times \transid)}$};
\node (A4) at (1.5,10){$\substack{\tens (\transid \times \tens) \ (\transid \times [\tens(\transid \times I)] \times \transid)}$};
\node (A5) at (2.5,10){$\substack{\tens (\transid \times \tens) \\ ( \transid \times \tens \times \transid) \\ (\transid \times \transid  \times I \times \transid)}$};
\node (A6) at (5,10){$\substack{\tens (\transid \times \tens)\\(\transid \times \transid \times \tens)\\ (\transid \times \transid  \times I \times \transid)}$};
\node (A7) at (5.5,9){$\substack{\tens (\transid \times \tens)\\(\transid \times \transid \times [\tens( I \times \transid)])}$};
\node (A8) at (5.5,7){$\substack{\tens (\transid \times \tens)}$};
\draw[doubleloose] (A1) to node[left]{$\substack{\alpha }$} (A3);
\draw[doubletighteq] (A3) to (A35);
\draw[doubleloose] (A35) to node[left, yshift=10pt]{$\substack{\looseid \looseid \\(\looseid \times r^{-1} \times  \looseid) }$} (A4);
\draw[doubletighteq] (A4) to (A5);
\draw[doubleloose] (A5) to node[above]{$\substack{\looseid (\looseid \times \alpha) \looseid }$} (A6);
\draw[doubletighteq] (A6) to (A7);
\draw[doubleloose] (A7) to node[right]{$\substack{\looseid \looseid \\ (\looseid \times \looseid \times l)}$} (A8);
\node (B1) at (.5,7.5){$\substack{\tens (\tens \times \transid) \\ (\transid \times \transid \times \transid) }$};
\node (B2) at (1.5,8.5){$\substack{\tens (\tens \times \transid) \\ (\transid \times [\tens(\transid \times I] \times \transid) }$};
\node (B3) at (2.5,9){$\substack{\tens (\tens \times \transid) \\ (\transid \times \tens \times \transid) \\ (\transid \times \transid \times I \times \transid)}$};
\draw[doubletighteq] (A1) to (B1);
\draw[doubleloose] (B1) to node[above, xshift=-15pt]{$\substack{\looseid \looseid \\(\looseid \times r^{-1} \times  \looseid)}$} (B2);
\draw[doubletighteq] (B2) to (B3);
\draw[doubleloose] (B3) to node[left]{$\substack{\alpha \looseid \looseid  }$} (A5);
\node (C0) at (1.5,5.5){$\substack{\tens (\transid \times \transid) \\(\tens \times \transid \times \transid) }$};
\node (C1) at (2.5,6){$\substack{\tens ([\tens (\transid \times I)] \times \transid) \\(\tens \times \transid \times \transid) }$};
\node (C15) at (3,7){$\substack{\tens (\tens \times \transid) \\(\transid \times I \times \transid) \\(\tens \times \transid \times \transid)}$};
\node (C2) at (3,8){$\substack{\tens (\tens \times \transid) \\(\tens \times \transid \times \transid) \\(\transid \times \transid \times I \times \transid) \\}$};
\draw[doubletighteq] (A1) to (C0);
\draw[doubleloose] (C0) to node[above]{$\substack{\looseid (r^{-1} \times \looseid)\\ \looseid}$} (C1);
\draw[doubletighteq] (C1) to (C15);
\draw[doubletighteq] (C15) to  (C2);
\draw[doubleloose] (C2) to node[above, xshift=15pt, yshift=-10pt]{$\substack{\looseid (\alpha \times \looseid) \\ \looseid }$} (B3);
\node (D3) at (4,8){$\substack{\tens (\transid \times \tens) \\ (\tens \times \transid \times  \transid) \\ (\transid \times \transid \times I \times \transid)}$};
\node (D5) at (5,8){$\substack{\tens (\tens \times \transid) \\ (\transid \times \transid \times \tens) \\ (\transid \times \transid \times I \times \transid)}$};
\node (D6) at (4,7){$\substack{\tens (\transid \times [\tens (I \times \transid)]) \\ (\tens \times \transid \times \transid)}$};
\node (D7) at (4.5,6){$\substack{\tens (\transid \times \transid) \\ (\tens \times \transid \times \transid)}$};
\draw[doubleloose] (C2) to node[above]{$\substack{\alpha \looseid \looseid }$} (D3);
\draw[doubletighteq] (D3) to (D5);
\draw[doubletighteq] (D3) to (D6);
\draw[doubleloose] (D5) to node[left]{$\substack{\alpha  \looseid \looseid  }$} (A7);
\node (E7) at (5,5){$\substack{\tens (\tens \times \transid)}$};
\draw[doubleloose] (D6) to node[right]{$\substack{\looseid (\looseid \times l) \looseid  }$} (D7);
\draw[doubletighteq] (D7) to (E7);
\draw[doubleloose] (E7) to node[above]{$\substack{\alpha }$} (A8);
\draw[doubleloose] (C0) to[in=135, out=0] node[above]{$ S(\mu) \looseid$} (D7);
\draw[doubleloose] (C0) to[in=225, out=-35] node[above]{$ T(\mu) \looseid$} (D7);
\draw[doubleloose] (A1) to[in=200, out=-20] node[above]{$\substack{\looseid}$} (E7);
\node at (1,9) {$\substack{\DDownarrow \iso  }$};
\node at (4.5,9.75) {$\substack{\DDownarrow \iso}$};
\node at (4,9) {$\substack{\DDownarrow \pi \tightid }$};
\node at (3.5,8.5) {$\substack{\DDownarrow \iso}$};
\node at (3.5,6.5) {$\substack{\DDownarrow \iso }$};
\node at (5,7) {$\substack{\DDownarrow \iso }$};
\node at (3,5.75) {$\substack{\DDownarrow \mu \tightid}$};
\node at (.75,7.25) {$\substack{\DDownarrow \iso}$};
\node at (1.5,7) {$\substack{\DDownarrow \tightid \rho}$};
\node at (2.25,6.75) {$\substack{\DDownarrow \iso}$};
\node at (2.5,4.75) {$\substack{\DDownarrow \iso  }$};
\draw[doubleloose] (A1) to[in=225, out=60] node[left]{$\substack{\looseid S(\rho)}$} (B3);
\draw[doubleloose] (A1) to[in=260, out=30] node[left]{$\substack{\looseid T(\rho)}$} (B3);
\draw[doubleloose] (A5) to[in=135, out=-30] node[above]{$\substack{ S(\pi)\looseid}$} (A7);
\draw[doubleloose] (A5) to[in=200, out=-60] node[below]{$\substack{T(\pi)\looseid }$} (A7);
\end{tikzpicture} \hspace{-2cm}
\end{aligned}
\end{equation*}
\begin{equation}\label{eq:monobjeq3}
=
\end{equation}
\begin{equation*}\hspace{-2cm}
\begin{aligned}
\begin{tikzpicture}[xscale=2.75, yscale=2]
\node (A1) at (0,7){$\substack{\tens (\tens \times \transid)}$};
\node (A3) at (0,8){$\substack{\tens (\transid \times \tens) }$};
\node (A35) at (0,9){$\substack{\tens (\transid \times \tens) \ (\transid \times \transid \times \transid)}$};
\node (A4) at (0,10){$\substack{\tens (\transid \times \tens) \ (\transid \times [\tens(\transid \times I)] \times \transid)}$};
\node (A5) at (2.5,10){$\substack{\tens (\transid \times \tens) \\ ( \transid \times \tens \times \transid) \\ (\transid \times \transid  \times I \times \transid)}$};
\node (A6) at (5,10){$\substack{\tens (\transid \times \tens)\\(\transid \times \transid \times \tens)\\ (\transid \times \transid  \times I \times \transid)}$};
\node (A7) at (5,9){$\substack{\tens (\transid \times \tens)\\(\transid \times \transid \times [\tens( I \times \transid)])}$};
\node (A8) at (5,8){$\substack{\tens (\transid \times \tens)}$};
\draw[doubleloose] (A1) to node[left]{$\substack{\alpha }$} (A3);
\draw[doubletighteq] (A3) to (A35);
\draw[doubleloose] (A35) to node[left, yshift=10pt]{$\substack{\looseid \looseid \\(\looseid \times r^{-1} \times  \looseid) }$} (A4);
\draw[doubletighteq] (A4) to (A5);
\draw[doubleloose] (A5) to node[above]{$\substack{\looseid (\looseid \times \alpha) \looseid }$} (A6);
\draw[doubletighteq] (A6) to (A7);
\draw[doubleloose] (A7) to node[right]{$\substack{\looseid \looseid \\ (\looseid \times \looseid \times l)}$} (A8);
\node (E7) at (5,7){$\substack{\tens (\tens \times \transid)}$};
\draw[doubleloose] (E7) to node[above]{$\substack{\alpha }$} (A8);
\draw[doubleloose] (A1) to node[above]{$\substack{\looseid}$} (E7);
\draw[doubleloose] (A3) to node[below]{$\substack{\looseid}$} (A8);
\draw[doubleloose] (A35) to[in=195, out=-15] node[above]{$\substack{\looseid T(\mu)}$} (A7);
\draw[doubleloose] (A35) to[in=165, out=15] node[below]{$\substack{\looseid S(\mu)}$} (A7);
\node at (2.5,7.5) {$\substack{\DDownarrow \iso  }$};
\node at (2.5,9.75) {$\substack{\DDownarrow \iso  }$};
\node at (2.5,9) {$\substack{\DDownarrow \mu \tightid}$};
\node at (2.5,8.25) {$\substack{\DDownarrow \iso  }$};
\end{tikzpicture} 
\end{aligned}\hspace{-2cm}
\end{equation*}} \newpage

\subsubsection*{Lax Monoidal 1-cell}

\input{coherence/LaxFunc1} \newpage
{\small
\begin{equation}\label{eq:laxfunc2}\hspace{-2cm}
\begin{tikzpicture}[xscale=2.25, yscale=1.5]
\node (A0) at (0,2) {$\substack{\tens (f\times f)}$};
\node (A3) at (0,6) {$\substack{f \tens }$};
\node (A31) at (.5,7.5) {$\substack{f \tens \\ (\transid \times \transid)}$};
\node (A32) at (2,8.5) {$\substack{f \tens \\ ([\tens(\transid \times I)] \times \transid)}$};
\node (A4) at (3,8.5) {$\substack{f \tens (\tens \times \transid) \\ (\transid \times I\times \transid)}$};
\node (A5) at (4.5,8) {$\substack{f \tens (\transid \times \tens) \\ (\transid \times I \times \transid)}$};
\node (A51) at (5.5,7) {$\substack{f \tens (\transid \times [\tens (I \times \transid)]) }$};
\node (A52) at (6,6) {$\substack{f \tens (\transid \times \transid) }$};
\node (A6) at (6,5) {$\substack{ f \tens }$};
\draw[doubleloose] (A0) to node[left] {$\substack{\chi}$} (A3);
\draw[doubleeq] (A3) to  (A31);
\draw[doubleloose] (A31) to node[above, yshift=10pt] {$\substack{\looseid \looseid (r^{-1} \times \transid)}$} (A32);
\draw[doubleeq] (A32) to (A4);
\draw[doubleloose] (A4) to node[above] {$\substack{\looseid \alpha \looseid}$} (A5);
\draw[doubleeq] (A5) to  (A51);
\draw[doubleloose] (A51) to node[right] {$\substack{\looseid  \looseid (\looseid \times l)}$} (A52);
\draw[doubleeq] (A52) to  (A6);
\node (B2) at (1,5.5) {$\substack{\tens (f \times f)\\(\transid \times \transid) }$};
\node (B3) at (2,7) {$\substack{\tens (f \times f)\\([\tens(\transid \times I)] \times \transid)}$};
\node (B35) at (2.25,7.75) {$\substack{\tens (f \times f)\\(\tens \times \transid) \\ (\transid \times I \times \transid)}$};
\draw[doubleeq] (A0) to (B2);
\draw[doubleloose] (B2) to node[left] {$\substack{\looseid \looseid (r^{-1} \times \looseid)}$} (B3);
\draw[doubleeq] (B3) to (B35);
\draw[doubleloose] (B35) to node[right, xshift=3pt] {$\substack{\chi \looseid \looseid}$} (A4);
\node (C1) at (.75,2) {$\substack{\tens (\transid \times \transid) \\ (f \times f)}$};
\node (C2) at (1.75,2.5) {$\substack{\tens ([\tens ( \transid \times I)]\times \transid) \\ (f \times f)}$};
\node (C3) at (2.5,3) {$\substack{\tens (\tens \times \transid) \\ (f \times I \times f)}$};
\node (C4) at (2.5,5) {$\substack{\tens \\(\tens \times \transid) \\ (f \times fI \times f)}$};
\draw[doubleeq] (A0) to (C1);
\draw[doubleloose] (C1) to node[above] {$\substack{\looseid \\ (r^{-1} \times \looseid)\looseid }$} (C2);
\draw[doubleeq] (C2) to (C3);
\draw[doubleloose] (C3) to node[left] {$\substack{\looseid \looseid \\ (\looseid \times \iota \times \looseid)}$} (C4);
\draw[doubleloose] (B3) to node[right] {$\substack{ \looseid (\chi \times \looseid) \looseid}$} (C4);
\node (D2) at (3.5,5.5) {$\substack{\tens (\transid \times \tens )\\ (f \times fI \times f) }$};
\node (D3) at (4,6.5) {$\substack{\tens \\(f \times [\tens (f\times f)] )\\ (\transid \times I \times \transid) }$};
\node (D4) at (5,5.5) {$\substack{ \tens (f \times f\tens) \\ (\transid \times I \times \transid)}$};
\node (D5) at (5,6.5) {$\substack{ \tens (f \times f) \\ (\transid \times \tens) \\ (\transid \times I \times \transid)}$};
\node (D6) at (5,4.5) {$\substack{ \tens (f \times f) \\ (\transid \times [\tens(I \times \transid)]) }$};
\node (D7) at (5.5,3.5) {$\substack{ \tens (f \times f) \\ (\transid \times  \transid) }$};
\draw[doubleloose] (C4) to node[above] {$\substack{ \alpha \looseid}$} (D2);
\draw[doubleloose] (D3) to node[below, xshift=-3pt] {$\substack{ \looseid \\ (\looseid \times \chi) \looseid}$} (D4);
\draw[doubleeq] (D2) to  (D3);
\draw[doubleeq] (D4) to  (D5);
\draw[doubleeq] (D4) to  (D6);
\draw[doubleloose] (D5) to node[right] {$\substack{\chi \looseid \looseid}$} (A5);
\node (E2) at (3.5,2.5) {$\substack{\tens ( \transid \times \tens) \\ (f \times I \times f)}$};
\node (E3) at (4.25,2) {$\substack{\tens ( \transid \times [\tens(I \times \transid)]) \\ (f \times f)}$};
\node (E4) at (5.25,1.5) {$\substack{\tens ( \transid \times  \transid) \\ (f  \times f)}$};
\node (E5) at (6,2) {$\substack{ \tens (f \times f)}$};
\draw[doubleloose] (C3) to node[below] {$\substack{\alpha \looseid }$} (E2);
\draw[doubleloose] (E2) to node[right] {$\substack{ \looseid \looseid \\ (\looseid  \times \iota \times \looseid)}$} (D2);
\draw[doubleloose] (D6) to node[right] {$\substack{\looseid  \looseid (\looseid \times l)}$} (D7);
\draw[doubleeq] (D7) to  (E5);
\draw[doubleloose] (E5) to node[right] {$\chi$} (A6);
\draw[doubleeq] (E2) to  (E3);
\draw[doubleloose] (E3) to node[above] {$\substack{ \looseid  (\looseid  \times l)}$} (E4);
\draw[doubleeq] (E4) to  (E5);
\draw[doubleloose] (A0) to[out=270, in=285] node[above] {$\substack{\looseid}$} (E5);
\node at (.5,6) {$\DDownarrow \iso$};
\node at (2,6) {$\DDownarrow \iso$};
\node at (1.5,4.5) {$\DDownarrow \tightid \delta$};
\node at (2,3.5) {$\DDownarrow \iso$};
\node at (3.4,8.1) {$\DDownarrow \iso$};
\node at (3.5,7.5
) {$\DDownarrow \omega \tightid$};
\node at (5.5,5) {$\DDownarrow \iso$};
\node at (4.25,4.5) {$\DDownarrow \iso$};
\node at (4.25,2.6) {$\DDownarrow \tightid \gamma$};
\node at (5.2,1.85) {$\DDownarrow \iso$};
\node at (3,3.5) {$\DDownarrow \iso$};
\node at (3,2) {$\DDownarrow \iso$};
\node at (3,1) {$\DDownarrow \mu$};
\node at (3,.25) {$\DDownarrow \iso$};
\draw[doubleloose] (A0) to[in=135, out=65] node[above,xshift=44pt, yshift=36pt]{$\substack{\looseid S(\delta)}$} (C4);
\draw[doubleloose] (A0) to[in=225, out=20] node[left]{$\substack{\looseid T(\delta)}$} (C4);
\draw[doubleloose] (B35) to node[below]{$\substack{S(\omega) \looseid }$} (A5);
\draw[doubleloose] (B35) to[in=245, out=-45] node[below]{$\substack{ T(\omega) \looseid}$} (A5);
\draw[doubleloose] (E2) to node[above,xshift=5pt]{$\substack{S(\gamma) \looseid }$} (E5);
\draw[doubleloose] (E2) to[in=155, out=60] node[above]{$\substack{ T(\gamma) \looseid}$} (E5);
\draw[doubleloose] (C1) to node[below]{$\substack{S(\mu) \looseid }$} (E4);
\draw[doubleloose] (C1) to[in=215, out=-45] node[below]{$\substack{T(\mu) \looseid }$} (E4);
\end{tikzpicture}\hspace{-2cm}
\end{equation}
\begin{equation*}
=
\end{equation*}
\begin{equation*}\hspace{-2cm}
\begin{tikzpicture}[xscale=2.25, yscale=1.5]
\node (A0) at (0,5) {$\substack{\tens (f\times f)}$};
\node (A3) at (0,6) {$\substack{f \tens }$};
\node (A31) at (.5,7.5) {$\substack{f \tens \\ (\transid \times \transid)}$};
\node (A32) at (2,8.5) {$\substack{f \tens \\ ([\tens(\transid \times I)] \times \transid)}$};
\node (A4) at (3,8.5) {$\substack{f \tens (\tens \times \transid) \\ (\transid \times I\times \transid)}$};
\node (A5) at (4.5,8.5) {$\substack{f \tens (\transid \times \tens) \\ (\transid \times I \times \transid)}$};
\node (A51) at (5.5,7.5) {$\substack{f \tens (\transid \times [\tens (I \times \transid)]) }$};
\node (A52) at (6,6.5) {$\substack{f \tens (\transid \times \transid) }$};
\node (A6) at (6,6) {$\substack{ f \tens }$};
\draw[doubleloose] (A0) to node[left] {$\substack{\chi}$} (A3);
\draw[doubleeq] (A3) to  (A31);
\draw[doubleloose] (A31) to node[above, yshift=10pt] {$\substack{\looseid \looseid (r^{-1} \times \transid)}$} (A32);
\draw[doubleeq] (A32) to (A4);
\draw[doubleloose] (A4) to node[above] {$\substack{\looseid \alpha \looseid}$} (A5);
\draw[doubleeq] (A5) to  (A51);
\draw[doubleloose] (A51) to node[right] {$\substack{\looseid  \looseid (\looseid \times l)}$} (A52);
\draw[doubleeq] (A52) to  (A6);
\node (E5) at (6,5) {$\substack{ \tens (f \times f)}$};
\draw[doubleloose] (E5) to node[right] {$\chi$} (A6);
\draw[doubleloose] (A0) to node[above] {$\substack{\looseid}$} (E5);
\node at (4,8) {$\DDownarrow \iso$};
\node at (3,7) {$\DDownarrow \mu$};
\node at (1.5,6.5) {$\DDownarrow \iso$};
\node at (3,5) {$\DDownarrow \iso$};
\draw[doubleloose] (A3) to node[below]{$\substack{ \looseid }$} (A6);
\draw[doubleloose] (A31) to[in=190, out=-35] node[below]{$\substack{T(\mu) \looseid }$} (A52);
\draw[doubleloose] (A31) to[in=145, out=20] node[below]{$\substack{S(\mu) \looseid }$} (A52);
\end{tikzpicture}\hspace{-2cm}
\end{equation*}} \newpage








\subsubsection*{Lax Monoidal 2-cell}

 \begin{equation*}\hspace{-2cm}
\begin{aligned}
\begin{tikzpicture}[xscale=3, yscale=1.5]
\node (t0) at (0,2) {\small $\tens(I_B \times f)$};
\node (t1) at (1,2) {\small $\tens(f I_A \times f)$};
\node (t15) at (2,2) {\small $\tens(f  \times f) (I \times \transid)$};
\node (t2) at (3,2) {\small $f \tens(I_A \times \transid)$};
\node (t3) at (4,2) {\small $f $};
\node (t4) at (5,2) {\small $g $};
\node (m0) at (0,1) {\small $\tens(I_B \times \transid)f$};
\node (b3) at (5,1) {\small $\tens (I_B \times \transid)g$};
\draw[doubleloose] (t0) to node[above]{$\substack {\looseid (\iota \times \looseid)}$} (t1);
\draw[doubleeq] (t1) to (t15);
\draw[doubleloose] (t15) to node[above]{$\substack {\chi \looseid}$} (t2);
\draw[doubleloose] (t2) to node[above]{$\substack {\looseid l}$} (t3);
\draw[doubleloose] (t3) to node[above]{$\substack {\beta}$} (t4);
\draw[doubleloose] (m0) to node[above]{$\substack {l \looseid}$} (t3);
\draw[doubleloose] (m0) to node[above]{$\substack {\looseid  \looseid \beta}$} (b3);
\draw[doubleloose] (b3) to node[right]{$\substack {l \looseid}$} (t4);
\draw[doubletighteq] (t0) to (m0);
\node at (1,1.5) {\small $\DDownarrow \gamma$};
\node at (2,1) {\small $\DDownarrow \tightid$};
\node at (3.5,1.5) {\small $\iso$};
\end{tikzpicture}
\end{aligned}\hspace{-2cm}
\end{equation*}
\begin{equation}\label{eq:mon2cell1}
=
\end{equation}
\begin{equation*}\hspace{-2cm}
\begin{aligned}
\begin{tikzpicture}[xscale=3, yscale=2]
\node (04) at (0,4) {\small $\tens(I_B \times f)$};
\node (14) at (1,5.5) {\small $\tens(f I_A\times f)$};
\node (154) at (1.5,6) {\small $\tens(f \times f) (I_A \times \transid)$};
\node (24) at (3,6.5) {\small $f \tens(I_A \times \transid_A)$};
\node (34) at (4.5,5.5) {\small $f $};
\node (44) at (5,4) {\small $g $};
\node (11) at (.5,2) {\small $\tens(I_B \times g)$};
\node (22) at (2,3) {\small $\tens(g I_A \times g)$};
\node (32) at (3,3) {\small $\tens(g \times g) (I_A \times \transid)$};
\node (33) at (4,4) {\small $g \tens(I_A \times \transid_A)$};
\node (01) at (0,1) {\small $\tens(I_B \times \transid)f$};
\node (10) at (5,1) {\small $\tens(  I_B \times \transid)g$};
\draw[doubleloose] (04) to node[above, xshift=-10]{$\substack {\looseid(\iota \times \looseid) }$} (14);
\draw[doubleeq] (14) to  (154);
\draw[doubleloose] (154) to node[above]{$\substack {\chi \looseid}$} (24);
\draw[doubleloose] (24) to node[above]{$\substack { \looseid l}$} (34);
\draw[doubleloose] (34) to node[above, xshift=3pt]{$\substack {\beta}$} (44);
\draw[doubleloose] (24) to node[right]{$\substack {\beta \looseid\looseid}$} (33);
\draw[doubleloose] (33) to node[above]{$\substack {\looseid l}$} (44);
\draw[doubleloose] (154) to node[below, xshift=-20pt]{$\substack {\looseid (\beta \times \beta) \looseid}$} (32);
\draw[doubleeq] (22) to  (32);
\draw[doubleloose] (32) to node[left]{$\substack {\chi \looseid}$} (33);
\draw[doubleloose] (04) to node[right]{$\substack {\looseid (\looseid \times \beta)}$} (11);
\draw[doubleloose] (11) to node[above,yshift=2pt]{$\substack {\looseid (\iota \times \looseid) }$} (22);
\draw[doubleloose] (10) to node[left]{$\substack {l \looseid} $} (44);
\draw[doubleeq] (04) to  (01);
\draw[doubleloose] (01) to node[above]{$\substack {\looseid \looseid \beta}$} (10);
\draw[doubleloose] (14) to node[left, yshift=5pt]{$\substack { \looseid (\beta \looseid\times \beta )}$} (22);
\node at (3,6) {$\substack {\DDownarrow \iso }$};
\node at (3,5) {$\substack {\DDownarrow \Pi^{\beta}\tightid}$};
\node at (3,3.75) {$\substack {\DDownarrow \iso }$};
\node at (4,5) {$\substack {\DDownarrow \iso }$};
\node at (1.75,4.75) {$\substack {\DDownarrow \iso }$};
\node at (1,4.75) {$\substack {\DDownarrow \iso }$};
\node at (.9,3.75) {$\substack {\DDownarrow \tightid_{\looseid} (M^{\beta} \times \tightid) }$};
\node at (.75,3) {$\substack {\DDownarrow \iso }$};
\node at (3,2) {$\substack {\DDownarrow \gamma^g}$};
\draw[doubleloose] (154) to[in=115, out=0]  node[below, xshift=-5pt]{$\substack{S(\Pi)\tightid}$} (33);
\draw[doubleloose] (154) to[in=180
, out=-45] node[above, xshift=5pt]{$\substack{T(\Pi)\tightid}$} (33);
\draw[doubleloose] (04) to[in=135, out=20]  node[above, xshift=3pt, yshift=5pt]{$\substack{\looseid (\beta \times S(M)) }$} (22);
\draw[doubleloose] (04) to[in=180
, out=-45] node[above, xshift=3pt,yshift=5pt]{$\substack{\looseid (\beta \times T(M))}$} (22);
\end{tikzpicture}
\end{aligned}\hspace{-2cm}
\end{equation*} \newpage

\begin{equation*}\hspace{-2cm}
\begin{aligned}
\begin{tikzpicture}[xscale=3, yscale=1.5]
\node (t0) at (0,2) {\small $\tens(f\times I_B)$};
\node (t1) at (1,2) {\small $\tens(f\times fI_A)$};
\node (t15) at (2,2) {\small $\tens (f \times f) (\transid \times I)$};
\node (t2) at (3,2) {\small $f\tens(\id \times I_A)$};
\node (t3) at (4,2) {\small $f $};
\node (t4) at (5,2) {\small $g $};
\node (m0) at (0,1) {\small $\tens(\transid \times I_B)f$};
\node (b3) at (5,1) {\small $\tens (\transid \times I_B)g$};
\draw[doubleloose] (t0) to node[above]{$\substack{\looseid(\looseid \times \iota)}$} (t1);
\draw[doubleeq] (t1) to  (t15);
\draw[doubleloose] (t15) to node[above]{$\substack{\chi (\looseid \times \looseid)}$} (t2);
\draw[doubleloose] (t2) to node[above]{$\substack{\looseid r}$} (t3);
\draw[doubleloose] (t3) to node[above]{$\substack{\beta}$} (t4);
\draw[doubleloose] (m0) to node[above]{$\substack{r \looseid}$} (t3);
\draw[doubleloose] (m0) to node[above]{$\substack{\looseid_{\tens}(\beta \times \looseid_I)}$} (b3);
\draw[doubleloose] (b3) to node[right]{$\substack{r \looseid}$} (t4);
\draw[doubletighteq] (t0) to (m0);
\node at (1,1.5) {$\substack{\DDownarrow \delta^f}$};
\node at (4,1.5) {$\substack{\iso}$};
\end{tikzpicture}
\end{aligned}\hspace{-2cm}
\end{equation*}
\begin{equation}\label{eq:mon2cell2}
=
\end{equation}
\begin{equation*}\hspace{-2cm}
\begin{aligned}
\begin{tikzpicture}[xscale=3, yscale=1.5]
\node (04) at (0,4) {\small $\tens(f\times I_B)$};
\node (14) at (1,5.5) {\small $\tens(f\times f I_A)$};
\node (154) at (1.5,6) {\small $\tens(f \times f) (\transid \times I_A)$};
\node (24) at (3.5,6.5) {\small $f \tens(\transid \times I_A)$};
\node (34) at (4.5,5.5) {\small $f $};
\node (44) at (5,4) {\small $g $};
\node (11) at (.5,2) {\small $\tens(g\times  I_B)$};
\node (22) at (2,2.5) {\small $\tens(g\times g I_A)$};
\node (32) at (3,3) {\small $\tens(g\times g) (\transid \times  I_A)$};
\node (33) at (4,4) {\small $g \tens(\transid \times I_A)$};
\node (00) at (0,1) {\small $\tens(\transid \times I_B)f$};
\node (10) at (5,1) {\small $\tens(\transid \times  I_B)g$};
\draw[doubleloose] (04) to node[above, xshift=-10pt]{$\substack{\looseid (\looseid \times \iota_f)}$} (14);
\draw[doubleeq] (14) to (154);
\draw[doubleloose] (154) to node[above]{$\substack{\chi \looseid}$} (24);
\draw[doubleloose] (24) to node[above]{$\substack{\looseid r}$} (34);
\draw[doubleloose] (34) to node[above, xshift=3pt]{$\substack{\beta}$} (44);
\draw[doubleloose] (24) to node[right]{$\substack{\beta \looseid \looseid}$} (33);
\draw[doubleloose] (33) to node[above]{$\substack{\looseid r}$} (44);
\draw[doubleloose] (154) to node[below, xshift=-20pt]{$\substack{\looseid (\beta \times \beta)\looseid }$} (32);
\draw[doubleeq] (22) to (32);
\draw[doubleloose] (32) to node[above, xshift=-10pt]{$\substack{\chi \looseid}$} (33);
\draw[doubleloose] (04) to node[right]{$\substack{\looseid (\beta \times \looseid)}$} (11);
\draw[doubleloose] (11) to node[above, xshift=-10pt]{$\substack{\looseid (\looseid \times \iota_g) }$} (22);
\draw[doubleloose] (00) to node[above]{$\substack{\looseid  \looseid \beta}$} (10);
\draw[doubleloose] (10) to node[left]{$\substack{r \looseid }$} (44);
\draw[doubleloose] (14) to node[right]{$\substack{ \looseid (\beta \times \beta \looseid)}$} (22);
\draw[doubletighteq] (04) to (00);
\draw[doubletighteq] (11) to (10);
\node at (4.25,5) {$\substack{\DDownarrow \iso }$};
\node at (3.25,5.75) {$\substack{\DDownarrow\iso }$};
\node at (3,3.75) {$\substack{\DDownarrow\iso }$};
\node at (2,3.75) {$\substack{\DDownarrow\iso }$};
\node at (3,5) {$\substack{\DDownarrow \Pi^{\beta}\tightid}$};
\node at (1,4.75) {$\substack {\DDownarrow \iso }$};
\node at (0.8,3.75) {$\substack {\DDownarrow \tightid_{\looseid} (\tightid \times M^{\beta} ) }$};
\node at (.75,2.75) {$\substack {\DDownarrow \iso }$};
\node at (1,1.5) {$\substack{\DDownarrow \iso}$};
\node at (4,2.5) {$\substack{\DDownarrow \delta^g}$};
\draw[doubleloose] (154) to[in=115, out=0]  node[below, xshift=3pt]{$\substack{S(\Pi)\tightid}$} (33);
\draw[doubleloose] (154) to[in=180
, out=-45] node[below, xshift=-3pt]{$\substack{T(\Pi)\tightid}$} (33);
\draw[doubleloose] (04) to[in=135, out=20]  node[above, xshift=3pt, yshift=5pt]{$\substack{\looseid (\beta \times S(M)) }$} (22);
\draw[doubleloose] (04) to[in=180
, out=-45] node[above, xshift=3pt,yshift=5pt]{$\substack{\looseid (\beta \times T(M))}$} (22);
\end{tikzpicture}
\end{aligned}\hspace{-2cm}
\end{equation*} \newpage

\begin{equation*}\hspace{-2cm}
\begin{aligned}
\begin{tikzpicture}[xscale=3.5, yscale=1.5]
\node (04) at (0,4) {\small$ \tens( \tens \times \transid)(f \times f \times f)$};
\node (054) at (0,5.5) {\small $ \tens([\tens(f  \times f)] \times f)$};
\node (14) at (1,6) {\small $ \tens(f \tens \times f)$};
\node (154) at (1.75,6.5) {\small $ \tens(f \times f)(\tens \times \transid)$};
\node (24) at (3,6) {\small $f \tens(\tens \times \transid)$};
\node (34) at (4,5.5) {\small $f\tens (\transid \times \tens)$};
\node (44) at (4,3.5) {\small $g \tens (\transid \times \tens)$};
\node (03) at (1,4) {\small $\tens( \transid \times \tens)(f \times f \times f)$ };
\node (13) at (1.5,4.5) {\small $\tens( f \times [\tens(f \times f)])$ };
\node (23) at (2.5,4.5) {\small $\tens (f \times f \tens)$ };
\node (33) at (3,5) {\small $\tens (f \times f )(\transid \times \tens)$};
\node (10) at (0,2) {\small $\tens (\tens \times \transid) (g \times g \times g)$};
\node (21) at (1,1) {\small $\tens ( \transid \times \tens) (g \times g \times g)$};
\node (31) at (2,.5) {\small $\tens (g \times [\tens (g \times g) ])$};
\node (32) at (3,1) {\small $\tens (g \times g \tens)$};
\node (41) at (3.5,2.5) {\small $\tens (g \times g) (\transid \times \tens)$};
\draw[doubleeq] (04) to  (054);
\draw[doubleloose] (054) to node[above]{$\substack{\looseid (\chi \times \looseid)}$} (14);
\draw[doubleeq] (14) to  (154);
\draw[doubleloose] (154) to node[above]{$\substack{\chi \looseid}$} (24);
\draw[doubleloose] (24) to node[above]{$\substack{\looseid \alpha}$} (34);
\draw[doubleloose] (34) to node[right]{$\substack{\beta \looseid \looseid}$} (44);
\draw[doubleloose] (04) to node[above]{$\substack{\alpha \looseid_{f \times f \times f}}$} (03);
\draw[doubleeq] (03) to  (13);
\draw[doubleloose] (13) to node[above]{$\substack{\looseid (\looseid \times \chi)}$} (23);
\draw[doubleeq] (23) to  (33);
\draw[doubleloose] (33) to node[above]{$\substack{\chi \looseid}$} (34);
\draw[doubleloose] (13) to node[left, yshift=-10pt, xshift=5pt]{$\substack{\looseid \\ (\beta \times \looseid (\beta \times \beta))}$} (31);
\draw[doubleloose] (23) to node[right]{$\substack{\looseid (\beta \times \beta \looseid)}$} (32);
\draw[doubleloose] (33) to node[left]{$\substack{\looseid (\beta \times \beta) \looseid}$} (41);
\draw[doubleloose] (03) to node[right]{$\substack{\looseid \looseid \\ (\beta \times \beta \times \beta)}$} (21);
\draw[doubleloose] (04) to node[left]{$\substack{\looseid (\looseid \times \looseid)\\(\beta \times \beta \times \beta)}$} (10);
\draw[doubleloose] (10) to node[above]{$\substack{\alpha \looseid}$} (21);
\draw[doubleeq] (21) to  (31);
\draw[doubleloose] (31) to node[above]{$\substack{\looseid (\looseid \times \chi)}$} (32);
\draw[doubleeq] (32) to (41);
\draw[doubleloose] (41) to node[above]{$\substack{\chi \looseid}$} (44);
\node at (2,5.5) {$\substack{\DDownarrow \omega^f}$};
\node at (3.55,3.45) {$\substack{\DDownarrow \iso}$};
\node at (3.5,4.5) {$\substack{\DDownarrow \Pi^{\beta} \tightid}$};
\node at (1.25,3) {$\substack{\DDownarrow \iso}$};
\node at (2.35,4) {$\substack{\DDownarrow \iso}$};
\node at (2.25,2.75) {$\substack{ \DDownarrow \tightid ({\tightid \times \Pi^{\beta})}}$};
\node at (2.1,1.35) {$\substack{\DDownarrow \iso}$};
\node at (3.25,3) {$\substack{\DDownarrow \iso}$};
\node at (.5,2.5) {$\substack{\DDownarrow \iso}$};
\draw[doubleloose] (13) to[in=105, out=-45] node[right, xshift=-25pt, yshift=30pt]{$\substack{\looseid \\(\beta \times S(\Pi))}$} (32);
\draw[doubleloose] (13) to[in=145, out=-70] node[left, yshift=-30pt, xshift=32pt]{$\substack{\looseid \\(\beta \times T(\Pi))}$} (32);
\draw[doubleloose] (33) to[in=105, out=-5] node[right, xshift=-20pt, yshift=20pt]{$\substack{\looseid \\(\beta \times S(\Pi))}$} (44);
\draw[doubleloose] (33) to[in=145, out=-70] node[left, yshift=-15pt, xshift=27pt]{$\substack{\looseid \\(\beta \times T(\Pi))}$} (44);
\end{tikzpicture}
\end{aligned}\hspace{-2cm}
\end{equation*}
\begin{equation}\label{eq:mon2cell3}
=
\end{equation}
\begin{equation*}\hspace{-2cm}
\begin{aligned}
\begin{tikzpicture}[xscale=3.25, yscale=1.5]
\node (04) at (-.5,4) {\small$ \tens( \tens \times \transid)(f \times f \times f)$};
\node (054) at (0,5.5) {\small $ \tens([\tens(f  \times f)] \times f)$};
\node (14) at (1,6) {\small $ \tens(f \tens \times f)$};
\node (154) at (1.75,6.5) {\small $ \tens(f \times f)(\tens \times \transid)$};
\node (24) at (3,6) {\small $f \tens(\tens \times \transid)$};
\node (34) at (4,5.5) {\small $f\tens (\transid \times \tens)$};
\node (44) at (4,4) {\small $g \tens (\transid \times \tens)$};
\node (33) at (3.5,3) {\small $g \tens (\tens \times  \transid)$};
\node (22) at (2,2) {\small $\tens (g  \tens \times g)$};
\node (42) at (3,2) {\small $\tens (g \times g)(\tens \times \transid)$};
\node (10) at (0,2) {\small $\tens (\tens \times \transid) (g \times g \times g)$};
\node (20) at (1,2) {\small $\tens ([\tens (g \times g)] \times g)$};
\node (11) at (1,1) {\small $\tens ( \transid \times \tens) (g \times g \times g)$};
\node (31) at (2,.5) {\small $\tens (g \times [\tens (g \times g) ])$};
\node (32) at (3,1) {\small $\tens (g \times g \tens)$};
\node (41) at (4,2) {\small $\tens (g \times g) (\transid \times \tens)$};
\draw[doubleeq] (04) to  (054);
\draw[doubleloose] (054) to node[above]{$\substack{\looseid(\chi \times \looseid)}$} (14);
\draw[doubleeq] (14) to  (154);
\draw[doubleloose] (154) to node[above]{$\substack{\chi \looseid}$} (24);
\draw[doubleloose] (24) to node[above]{$\substack{\looseid \alpha}$} (34);
\draw[doubleloose] (34) to node[right]{$\substack{ \beta \looseid \looseid}$} (44);
\draw[doubleloose] (154) to node[left]{$\substack{ \looseid (\beta \times \beta) \looseid}$} (42);
\draw[doubleloose] (24) to node[right]{$\substack{\beta \looseid \looseid}$} (33);
\draw[doubleloose] (33) to node[above,yshift=2pt]{$\substack{ \looseid \alpha}$} (44);
\draw[doubleloose] (054) to node[left]{$\substack{\looseid \\ ([\looseid( \beta\times \beta)] \times \beta)}$} (20);
\draw[doubleloose] (14) to node[right]{$\substack{\looseid \\ (\beta \looseid \times \beta)}$} (22);
\draw[doubleeq] (22) to  (42);
\draw[doubleloose] (42) to node[above]{$\substack{\chi \looseid}$} (33);
\draw[doubleeq] (10) to (20);
\draw[doubleloose] (20) to node[above]{$\substack{\looseid (\chi \times \looseid)}$} (22);
\draw[doubleloose] (04) to node[left]{$\substack{\looseid (\looseid \times \looseid)\\(\beta \times \beta \times \beta)}$} (10);
\draw[doubleloose] (10) to node[above]{$\substack{\alpha \looseid}$} (21);
\draw[doubleeq] (21) to  (31);
\draw[doubleloose] (31) to node[above]{$\substack{\looseid (\looseid \times \chi)}$} (32);
\draw[doubleeq] (32) to (41);
\draw[doubleloose] (41) to node[right]{$\substack{\chi \looseid}$} (44);
\node at (3.65,4.5) {$\substack{ \DDownarrow \iso}$};
\node at (2.75, 5.65) {$\substack{ \DDownarrow \iso}$};
\node at (2.5,5) {$\substack{\DDownarrow \Pi^{\beta} \tightid}$};
\node at (3,3.5) {$\substack{ \DDownarrow \iso}$};
\node at (1.75,4.5) {$\substack{\DDownarrow \iso}$};
\node at (2,1.5) {$\substack{\DDownarrow \omega^g}$};
\node at (.75,5) {$\substack{ \DDownarrow \iso}$};
\node at (.85,4) {$\substack{\DDownarrow \tightid (\Pi^{\beta} \times  \tightid)}$};
\node at (1.25,2.5) {$\substack{ \DDownarrow \iso}$};
\node at (0,3.25) {$\substack{\DDownarrow \iso}$};
\draw[doubleloose] (054) to[in=105, out=-45] node[right, xshift=-30pt, yshift=25pt]{$\substack{\looseid \\( S(\Pi) \times \beta)}$} (22);
\draw[doubleloose] (054) to[in=145, out=-70] node[right, yshift=-10pt, xshift=20pt]{$\substack{\looseid \\( T(\Pi) \times \beta)}$} (22);
\draw[doubleloose] (154) to[in=105, out=-45] node[right, xshift=-20pt, yshift=20pt]{$\substack{\beta S(\Pi) }$} (33);
\draw[doubleloose] (154) to[in=145, out=-70] node[right, yshift=0pt, xshift=0pt]{$\substack{\beta  T(\Pi) }$} (33);
\end{tikzpicture}
\end{aligned}\hspace{-2cm}
\end{equation*} \newpage




\newpage
\subsubsection*{Monoidal Icon}

\begin{equation*}\hspace{-2cm}
\begin{aligned}
\begin{tikzpicture}[xscale=4,yscale=2]
\node (02) at (0,2){$\tens(I \times f ) $};
\node (12) at (1,2){$\tens(fI_A \times f) $};
\node (22) at (2,2){$f\tens(I_A \times \transid) $};
\node (32) at (3,2){$f $};
\node (01) at (0,1){$\tens(I_B \times \transid ) f$};
\node (31) at (3,1){$f $};
\node (00) at (0,0){$\tens(I_B \times \transid ) g$};
\node (30) at (3,0){$g$};
\draw[doubleloose] (02) to node[above]{$\looseid_{\tens} (\iota_f \times \looseid_f) $} (12);
\draw[doubleloose] (12) to node[above]{$\chi \looseid_{(I_A \times \transid)}$} (22);
\draw[doubleloose] (22) to node[above]{$\looseid_f l$} (32);
\draw[doubleloose] (01) to node[above]{$l \looseid_f $} (31);
\draw[doubleloose] (00) to node[above]{$l \looseid_g $} (30);
\draw[=] (02) to node[left]{} (01);
\draw[=] (32) to node[left]{} (31);
\draw[doubletight] (01) to node[left]{$\tightid_{\tens} (\tightid_I\times \beta)$} (00);
\draw[doubletight] (31) to node[left]{$\beta$} (30);
\node at (1.5,1.5){$\DDownarrow \gamma^f$};
\node at (1.5,0.5){$\DDownarrow \looseid_{l}\tightid_{\beta}$};
\end{tikzpicture}
\end{aligned}\hspace{-2cm}
\end{equation*}
\begin{equation}\label{eq:monicon1}
  =
\end{equation}
\begin{equation*}\hspace{-2cm}
\begin{aligned}
\begin{tikzpicture}[xscale=4,yscale=2]
\node (02) at (0,2){$\tens(I_B \times f ) $};
\node (12) at (1,2){$\tens(fI_A \times f) $};
\node (22) at (2,2){$f\tens(I_A \times \transid) $};
\node (32) at (3,2){$f $};
\node (01) at (0,1){$\tens(I_B \times g ) $};
\node (11) at (1,1){$\tens(gI_A \times g) $};
\node (21) at (2,1){$g\tens(I_A \times \transid) $};
\node (31) at (3,1){$g $};
\node (00) at (0,0){$\tens(I_B \times \transid ) g$};
\node (30) at (3,0){$g $};
\draw[doubleloose] (02) to node[above]{$\looseid_{\tens} (\iota_f \times \looseid_f) $} (12);
\draw[doubleloose] (12) to node[above]{$\chi \looseid_{(I_A \times \transid)}$} (22);
\draw[doubleloose] (22) to node[above]{$\looseid_f l$} (32);
\draw[doubleloose] (01) to node[above]{$\looseid_{\tens} (\iota_g \times \looseid_g $} (11);
\draw[doubleloose] (11) to node[above]{$\chi \looseid_{(I_A \times \transid)}$} (21);
\draw[doubleloose] (21) to node[above]{$\looseid_g l$} (31);
\draw[doubleloose] (00) to node[above]{$l \looseid_g $} (30);
\draw[doubletight] (02) to node[left]{$\tightid_{\tens(I\times \transid)}\beta$} (01);
\draw[doubletight] (12) to node[right]{$\tightid (\beta \times \beta) \tightid$} (11);
\draw[doubletight] (22) to node[left]{$\beta \tightid$} (21);
\draw[doubletight] (32) to node[left]{$\beta$} (31);
\draw[doubletight] (01) to node[left]{$\tightid_{\tens} (\tightid_I\times \beta)$} (00);
\draw[doubletight] (31) to node[left]{$\beta$} (30);
\node at (0.5,1.5){$\DDownarrow \tightid (N^{\beta} \times \tightid_{\looseid}) \tightid$};
\node at (1.5,1.5){$\DDownarrow \Sigma^{\beta} \tightid$};
\node at (2.5,1.5){$\DDownarrow \looseid_{\beta} \tightid_{l}$};
\node at (1.5,0.5){$\DDownarrow \gamma^g$};
\end{tikzpicture}
\end{aligned}\hspace{-2cm}
\end{equation*}\newpage

\begin{equation*}\hspace{-2cm}
\begin{aligned}
\begin{tikzpicture}[xscale=4,yscale=2]
\node (02) at (0,2){$\tens(f \times I_B ) $};
\node (12) at (1,2){$\tens(f \times fI_A) $};
\node (22) at (2,2){$f\tens(\transid \times I_A) $};
\node (32) at (3,2){$f $};
\node (01) at (0,1){$\tens(\transid \times I_B) f$};
\node (31) at (3,1){$f $};
\node (00) at (0,0){$\tens(\transid \times I_B ) g$};
\node (30) at (3,0){$g $};
\draw[doubleloose] (02) to node[above]{$\looseid_{\tens} (\looseid_f \times \iota_f) $} (12);
\draw[doubleloose] (12) to node[above]{$\chi \looseid_{(\transid \times I_A)}$} (22);
\draw[doubleloose] (22) to node[above]{$\looseid_f r$} (32);
\draw[doubleloose] (01) to node[above]{$r \looseid_f $} (31);
\draw[doubleloose] (00) to node[above]{$r \looseid_g $} (30);
\draw[=] (02) to node[left]{} (01);
\draw[=] (32) to node[left]{} (31);
\draw[doubletight] (01) to node[left]{$\tightid_{\tens (\transid \times I)} \beta$} (00);
\draw[doubletight] (31) to node[left]{$\beta$} (30);
\node at (1.5,1.5){$\DDownarrow \delta^f$};
\node at (1.5,0.5){$\DDownarrow \looseid_{r}\tightid_{\beta}$};
\end{tikzpicture}
\end{aligned}\hspace{-2cm}
\end{equation*}
\begin{equation}\label{eq:monicon2}
  =
\end{equation}
\begin{equation*}\hspace{-2cm}
\begin{aligned}
\begin{tikzpicture}[xscale=4,yscale=2]
\node (02) at (0,2){$\tens(f \times I_B ) $};
\node (12) at (1,2){$\tens(f \times fI_A) $};
\node (22) at (2,2){$f\tens(\transid \times I_A) $};
\node (32) at (3,2){$f$};
\node (01) at (0,1){$\tens(g \times I_B) $};
\node (11) at (1,1){$\tens(g\times gI_A ) $};
\node (21) at (2,1){$g\tens(\transid \times I_A)  $};
\node (31) at (3,1){$g $};
\node (00) at (0,0){$\tens(\transid \times I_B ) g$};
\node (30) at (3,0){$g $};
\draw[doubleloose] (02) to node[above]{$\looseid_{\tens} (\looseid \times\iota_f) $} (12);
\draw[doubleloose] (12) to node[above]{$\chi \looseid_{(\transid \times I_A)}$} (22);
\draw[doubleloose] (22) to node[above]{$\looseid_f r$} (32);
\draw[doubleloose] (01) to node[above]{$\looseid_{\tens} (\looseid_g \times \iota_g) $} (11);
\draw[doubleloose] (11) to node[above]{$\chi \looseid_{(\transid \times I_A)}$} (21);
\draw[doubleloose] (21) to node[above]{$\looseid_g r$} (31);
\draw[doubleloose] (00) to node[above]{$r \looseid_g $} (30);
\draw[doubletight] (02) to node[left]{$\tightid_{\tens(\transid\times I)}\beta$} (01);
\draw[doubletight] (12) to node[right]{$\tightid (\beta \times \beta) \tightid$} (11);
\draw[doubletight] (22) to node[left]{$\beta \tightid$} (21);
\draw[doubletight] (32) to node[left]{$\beta$} (31);
\draw[doubletight] (01) to node[left]{$\tightid_{\tens} (\beta\times \tightid_I)$} (00);
\draw[doubletight] (31) to node[left]{$\beta$} (30);
\node at (0.5,1.5){$\DDownarrow \tightid (\tightid \times N^{\beta}) \tightid$};
\node at (1.5,1.5){$\DDownarrow \Sigma^{\beta} \tightid$};
\node at (2.5,1.5){$\DDownarrow \looseid_{\beta} \tightid_{r}$};
\node at (1.5,0.5){$\DDownarrow \delta^g$};
\end{tikzpicture}
\end{aligned}\hspace{-2cm}
\end{equation*}\newpage

\begin{equation*}\hspace{-2cm}
\begin{aligned}
\begin{tikzpicture}[xscale=4,yscale=2]
\node (02) at (0,2){$\tens(\tens \times \transid)(f \times f\times f)$};
\node (12) at (1,2){$\tens(f\tens \times f) $};
\node (22) at (2,2){$f\tens(\tens \times \transid)$};
\node (32) at (3,2){$f\tens(\transid \times \tens)$};
\node (01) at (0,1){$\tens(\tens \times \transid)(f \times f\times f)$};
\node (11) at (1,1){$\tens(\transid \times \tens)(f \times f\times f)$};
\node (21) at (2,1){$\tens(f \times f\tens) $};
\node (31) at (3,1){$f \tens ( \transid \times \tens)$};
\node (00) at (0,0){$\tens(\tens \times \transid)(g \times g \times g \times g$};
\node (10) at (1,0){$\tens(\transid \times \tens)(g \times g\times g)$};
\node (20) at (2,0){$\tens(g \times g\tens) $};
\node (30) at (3,0){$g \tens ( \transid \times \tens)$};
\draw[doubleloose] (02) to node[above]{$\looseid_{\tens} (\chi \times \looseid_f)$} (12);
\draw[doubleloose] (12) to node[above]{$\chi \looseid_{(\tens \times \transid)}$} (22);
\draw[doubleloose] (22) to node[above]{$\looseid_{f} \alpha$} (32);
\draw[doubleloose] (01) to node[above]{$\alpha \looseid_{f\times f \times f}$} (11);
\draw[doubleloose] (11) to node[above]{$\looseid_{\tens} (\transid_f \times \chi)$} (21);
\draw[doubleloose] (21) to node[above]{$\chi \looseid_{\transid \times \tens}$} (31);
\draw[doubleloose] (00) to node[above]{$\alpha \looseid_{g\times g \times g} $} (10);
\draw[doubleloose] (10) to node[above]{$\looseid_{\tens} (\transid_g \times \chi)$} (20);
\draw[doubleloose] (20) to node[above]{$\chi \looseid_{\transid \times \tens}$} (30);
\draw[=] (02) to node[left]{} (01);
\draw[=] (32) to node[left]{} (31);
\draw[doubletight] (01) to node[left]{$\tightid (\beta\times \beta \times \beta)$} (00);
\draw[doubletight] (11) to node {$\tightid (\beta\times \beta \times \beta)$} (10);
\draw[doubletight] (21) to node {$\tightid (\beta\times \beta \tightid)$} (20);
\draw[doubletight] (31) to node[left]{$\beta \tightid$} (30);
\node at (1.5,1.5){$\DDownarrow \omega^f$};
\node at (0.5,0.5){$\DDownarrow \tightid_{\alpha} \looseid_{\beta \times \beta \times \beta}$};
\node at (1.5,0.5){$\DDownarrow \tightid (\tightid \times \Sigma^{\beta})$};
\node at (2.5,0.5){$\DDownarrow \Sigma^{\beta} \tightid$};
\end{tikzpicture}
\end{aligned}\hspace{-2cm}
\end{equation*}
\begin{equation}\label{eq:monicon3}
  =
\end{equation}
\begin{equation*}\hspace{-2cm}
\begin{aligned}
\begin{tikzpicture}[xscale=4,yscale=2]
\node (02) at (0,2){$\tens(\tens \times \transid)(f \times f\times f)$};
\node (12) at (1,2){$\tens(f\tens \times f) $};
\node (22) at (2,2){$f\tens(\tens \times \transid)$};
\node (32) at (3,2){$f\tens(\transid \times \tens)$};
\node (01) at (0,1){$\tens(\tens \times \transid)(g \times g\times g)$};
\node (11) at (1,1){$\tens(g\tens \times g)$};
\node (21) at (2,1){$g\tens(\tens \times \transid) $};
\node (31) at (3,1){$g \tens ( \transid \times \tens)$};
\node (00) at (0,0){$\tens(\tens \times \transid)(g \times g \times g \times g$};
\node (10) at (1,0){$\tens(\transid \times \tens)(g \times g\times g)$};
\node (20) at (2,0){$\tens(g \times g\tens) $};
\node (30) at (3,0){$g \tens ( \transid \times \tens)$};
\draw[doubleloose] (02) to node[above]{$\looseid_{\tens} (\chi \times \looseid_f)$} (12);
\draw[doubleloose] (12) to node[above]{$\chi \looseid_{(\tens \times \transid)}$} (22);
\draw[doubleloose] (22) to node[above]{$\looseid_{f} \alpha$} (32);
\draw[doubleloose] (01) to node[above]{$\looseid_{\tens} (\chi \times \looseid_g)$} (11);
\draw[doubleloose] (11) to node[above]{$\chi \looseid_{(\tens \times \transid)}$} (21);
\draw[doubleloose] (21) to node[above]{$\looseid_{g} \alpha$} (31);
\draw[doubleloose] (00) to node[above]{$\alpha \looseid_{g\times g \times g} $} (10);
\draw[doubleloose] (10) to node[above]{$\looseid_{\tens} (\transid_g \times \chi)$} (20);
\draw[doubleloose] (20) to node[above]{$\chi \looseid_{\transid \times \tens}$} (30);
\draw[=] (01) to node[left]{} (00);
\draw[doubletight] (12) to node {$\tightid (\beta \tightid \times \beta)$} (11);
\draw[doubletight] (22) to node[left] {$\beta \tightid $} (21);
\draw[=] (31) to node[left]{} (30);
\draw[doubletight] (02) to node[left]{$\tightid (\beta\times \beta \times \beta)$} (01);
\draw[doubletight] (32) to node[left]{$\beta \tightid$} (31);
\node at (1.5,0.5){$\DDownarrow \omega^g$};
\node at (0.5,1.5){$\DDownarrow \tightid (\Sigma^{\beta} \times \tightid)$};
\node at (1.5,1.5){$\DDownarrow \Sigma^{\beta} \tightid$};
\node at (2.5,1.5){$\DDownarrow \looseid_{\beta} \tightid_{\alpha}$};
\end{tikzpicture}
\end{aligned}\hspace{-2cm}
\end{equation*}\newpage

\subsubsection*{Braided Monoidal Icon}

\begin{equation*}\hspace{-2cm}
\begin{aligned}
\begin{tikzpicture}[xscale=4, yscale=2]
\node (02) at (0,2){$\tens(f \times f)$};
\node (12) at (1,2){$\tens \tau (f \times f)$};
\node (22) at (2,2){$f \tens \tau$};
\node (01) at (0,1){$\tens(f \times f)$};
\node (11) at (1,1){$f \tens $};
\node (21) at (2,1){$f \tens \tau$};
\node (00) at (0,0){$\tens (g \times g)$};
\node (10) at (1,0){$g\tens $};
\node (20) at (2,0){$g \tens \tau$};
\draw[doubleloose] (02) to node[above]{$\sigma \looseid_{f\times f}$} (12);
\draw[doubleloose] (12) to node[above]{$\chi \looseid_{\tau}$} (22);
\draw[doubleloose] (01) to node[above]{$\chi$} (11);
\draw[doubleloose] (11) to node[above]{$\looseid_f \sigma$} (21);
\draw[doubleloose] (00) to node[above]{$\chi$} (10);
\draw[doubleloose] (10) to node[above]{$\looseid_g \sigma$} (20);
\draw[=] (02) to (01);
\draw[doubletight] (01) to node[left]{$\looseid_{\tens} (\beta \times \beta)$}(00);
\draw[doubletight] (11) to node[left]{$\beta \looseid_{\tens}$}(10);
\draw[=] (22) to (21);
\draw[doubletight] (21) to node[left]{$\beta \looseid_{\tens \tau}$}(20);
\node at (1,1.5){$\DDownarrow u$};
\node at (.5,.5){$\DDownarrow\Pi^{\beta}$};
\node at (1.5,.5) {$\DDownarrow \looseid_{\beta} \tightid_{\sigma}$};
\end{tikzpicture}
\end{aligned}\hspace{-2cm}
\end{equation*}
\begin{equation}\label{eq:bricon}
=
\end{equation}
\begin{equation*}\hspace{-2cm}
\begin{aligned}
\begin{tikzpicture}[xscale=4, yscale=2]
\node (02) at (0,2){$\tens(f \times f)$};
\node (12) at (1,2){$\tens \tau (f \times f)$};
\node (22) at (2,2){$f \tens \tau$};
\node (01) at (0,1){$\tens(g \times g)$};
\node (11) at (1,1){$ \tens \tau (g\times g)$};
\node (21) at (2,1){$g \tens \tau$};
\node (00) at (0,0){$\tens (g \times g)$};
\node (10) at (1,0){$g\tens $};
\node (20) at (2,0){$g \tens \tau$};
\draw[doubleloose] (02) to node[above]{$\sigma \looseid_{f\times f}$} (12);
\draw[doubleloose] (12) to node[above]{$\chi \looseid_{\tau}$} (22);
\draw[doubleloose] (01) to node[above]{$\sigma \looseid_{g \times g}$} (11);
\draw[doubleloose] (11) to node[above]{$\chi \looseid_{\tau}$} (21);
\draw[doubleloose] (00) to node[above]{$\chi$} (10);
\draw[doubleloose] (10) to node[above]{$\looseid_g \sigma$} (20);
\draw[doubletight] (02) to node[left]{$\looseid_{\tens} (\beta \times \beta)$}(01);
\draw[=] (01) to (00);
\draw[doubletight] (12) to node[xshift=5pt]{$\looseid_{\tens \tau} \hspace{.1cm} (\beta \times \beta) $}(11);
\draw[doubletight] (22) to node[left]{$\beta \looseid_{\tens \tau}$}(21);
\draw[=] (21) to (20);
\node at (1,0.5){$\DDownarrow u$};
\node at (.5,1.5){$\DDownarrow \tightid_{\sigma} \looseid_{\beta \times \beta} $};
\node at (1.5,1.5) {$\DDownarrow \Pi^{\beta}$};
\end{tikzpicture}
\end{aligned}\hspace{-2cm}
\end{equation*}

\subsubsection*{Monoidal 3-cell}

\begin{equation}\label{eq:mon3cell1}
\begin{aligned}
 \begin{tikzpicture}[scale=1.8]
 \node (tm) at (0,1) {\scriptsize $f  I_A$};
 \node (tr) at (1,1) {\scriptsize$g  I_A$};
 \node (bm) at (0,0) {\scriptsize$f' I_A$};
 \node (br) at (01,0) {\scriptsize$g' I_A$}; 
 \draw[doubleloose] (tm)  to node[above]{\scriptsize$\alpha \looseid_I$} (tr);
 \draw[doubletight] (tm) to node[right, yshift=8] {\scriptsize$\gamma \tightid_I$} (bm);
 \draw[doubleloose] (bm) to node[above] {\scriptsize$\beta \looseid_I$}(br);
  \draw[-implies, double equal sign distance] (tr) to node[right] {\scriptsize $\delta \tightid_I$}(br);
 \node at (0.5,0.5) {\scriptsize $\DDownarrow \Gamma \tightid_{\looseid}$}; 
 \node (tl) at (-1,1) {\scriptsize$I_B$};
 \node (bl) at (-1,0) {\scriptsize$I_B$};
 \draw[doubleloose] (tl)  to node[above]{\scriptsize$\iota_f$} (tm);
 \draw[doubleeq] (tl) to (bl);
 \draw[doubleloose] (bl) to node[above]{\scriptsize$\iota_{f'}$}(bm);
 \node at (-0.5,.5) {\scriptsize $\DDownarrow N^{\gamma}$};
\node (bl1) at (-1,-.7){\scriptsize$I_B$};  
 \node (bm1) at (0,-.7) {\scriptsize$I_B$};
  \node (br1) at (1,-.7) {\scriptsize$g' I_A$}; 
 \draw[doubleloose] (bl1)  to node[above]{\scriptsize$\looseid_{I}$} (bm1);
 \draw[doubleloose] (bm1) to  node[above]{\scriptsize$\iota_{g'}$}(br1);
  \draw[doubleeq] (bl)  to (bl1);
    \draw[doubleeq] (br)  to (br1);
 \node at (0,-0.35) {\scriptsize $\DDownarrow M^{\beta}$}; 
 \end{tikzpicture}
\end{aligned}
 =
 \begin{aligned}
  \begin{tikzpicture}[scale=1.8]
 \node (ml) at (0,1) {\scriptsize$I_B$};
 \node (mm) at (1,1) {\scriptsize$I_B$};
 \node (bl) at (0,0) {\scriptsize$I_B$};
 \node (bm) at (01,0) {\scriptsize$I_B$}; 
 \draw[doubleloose] (ml)  to node[above]{\scriptsize$ \looseid_{I}$}(mm);
 \draw[doubleeq] (ml) to  (bl);
 \draw[doubleloose] (bl) to  node[above]{\scriptsize$ \looseid_{I}$}(bm);
 \draw[doubleeq] (mm) to (bm);
 \node at (0.5,0.5) {\scriptsize $=$}; 
 \node (tl) at (0,1.7) {\scriptsize$I_B$};
 \node (tm) at (1,1.7) {\scriptsize$f I_A$};
 \node (tr) at (2,1.7) {\scriptsize$g I_A$};
 \node (mr) at (2,1) {\scriptsize$g I_A$};
 \draw[doubleloose] (tl)  to node[above]{\scriptsize$\iota_f$} (tm);
 \draw[doubleloose] (tm) to node[above]{\scriptsize$\alpha \looseid_I$} (tr);
 \draw[doubleloose] (mm) to node[above]{\scriptsize$\iota_{g}$}(mr);
 \node at (1,1.35) {\scriptsize\footnotesize $\DDownarrow M^{\alpha}$};
  \node (br) at (2,0) {$g' I$};
 \draw[doubleloose] (bm)  to node[above]{\scriptsize$\iota_{g'}$} (br);
 \draw[doubletight] (mr) to  node[right]{\scriptsize$\delta \tightid_I$}(br);
 \draw[doubleeq] (tr) to (mr);
  \draw[doubleeq] (tl) to (ml);
 \node at (1.5,.5) {\scriptsize$\DDownarrow N^{\delta}$}; 
 \end{tikzpicture}
 \end{aligned}
\end{equation}

\begin{equation}\label{eq:mon3cell2}
\begin{aligned}
 \begin{tikzpicture}[yscale=2, xscale=2]
 \node (tm) at (0,1) {\scriptsize$f\ten$};
 \node (tr) at (1,1) {\scriptsize$g \ten$};
 \node (mm) at (0,0) {\scriptsize$f' \ten$};
 \node (mr) at (01,0) {\scriptsize$g' \ten$}; 
 \draw[doubleloose] (tm)  to node[above]{\scriptsize$\alpha  \looseid_{\ten}$} (tr);
 \draw[doubletight] (tm) to node[right, yshift=8]{\scriptsize$\gamma \tightid_{\ten}$} (mm);
 \draw[doubleloose] (mm) to node[above, xshift=1pt, yshift=-1pt] {\scriptsize$\beta \looseid_{\ten}$}(mr);
  \draw[doubletight] (tr) to node[right] {\scriptsize$\delta \tightid_{\ten}$}(mr);
 \node at (0.5,0.5) {\scriptsize$\DDownarrow \Gamma \tightid$}; 
 \node (tl) at (-1,1) {\scriptsize$\ten  (f\times f)$};
 \node (ml) at (-1,0) {\scriptsize$\ten  (f'\times f')$};
 \draw[doubleloose] (tl)  to node[above]{\scriptsize$\chi^f$} (tm);
 \draw[doubletight] (tl) to node[left]{\scriptsize$\tightid_{\ten} (\gamma \times \gamma)$} (ml);
 \draw[doubleloose] (ml) to node[above]{\scriptsize$\chi^{f'}$}(mm);
 \node at (-0.5,0.5) {\scriptsize $\DDownarrow \Sigma^{\gamma}$};
 \node (bl) at (-1,-.7) {\scriptsize$\ten (f'\times f')$};
  \node (bm) at (0,-.7) {\scriptsize$\ten (g'\times g')$};
  \node (br) at (1,-.7) {\scriptsize$g' \ten$};
  \draw[doubleeq] (ml) to (bl);
 \draw[doubleloose] (bl)  to node[above]{\scriptsize$\looseid_{\ten} (\beta \times \beta)$} (bm);
 \draw[doubleloose] (bm) to  node[above]{\scriptsize$\chi^{g'}$}(br);
   \draw[doubleeq] (mr) to (br);
 \node at (0,-0.35) {\scriptsize $\DDownarrow \Pi^{\beta}$}; 
 \end{tikzpicture}
\end{aligned}
 =
 \begin{aligned}
  \begin{tikzpicture}[yscale=2, xscale=2]
 \node (ml) at (0,1) {\scriptsize$\ten (f\times f)$};
 \node (mm) at (1,1) {\scriptsize$\ten (g\times g)$};
 \node (bl) at (0,0) {\scriptsize$\ten (f'\times f')$};
 \node (bm) at (01,0) {\scriptsize$\ten (g'\times g')$}; 
 \draw[doubleloose] (ml)  to node[above]{\scriptsize$\looseid_{\ten} (\alpha \times \alpha)$} (mm);
 \draw[doubletight] (ml) to node[left]{\scriptsize$\tightid_{\ten} (\gamma \times \gamma)$}  (bl);
 \draw[doubleloose] (bl) to node [below] {\scriptsize$\looseid_{\ten} (\beta \times \beta)$} (bm);
  \draw[doubletight] (mm) to node[above] {\scriptsize$\tightid_{\ten} (\delta \times \delta)$} (bm);
 \node at (0.5,0.5) {\scriptsize\footnotesize $\DDownarrow \tightid (\Gamma \times \Gamma)$}; 
 \node (tl) at (0,1.7) {\scriptsize$ \ten (f \times f$)};
 \node (tm) at (1,1.7) {\scriptsize$f \ten$};
 \node (tr) at (2,1.7) {\scriptsize$g \ten$};
   \node (mr) at (2,1) {\scriptsize$g \ten$};
   \node(br) at (2,0) {\scriptsize$g' \ten$};
 \draw[doubleloose] (tl)  to node[above]{\scriptsize$\chi^f$} (tm);
 \draw[doubleloose] (tm) to node[above]{\scriptsize$\alpha \looseid_{\ten}$} (tr);
 \draw[doubletick] (mm) to node[above]{\scriptsize$\chi^{g}$}(mr);
 \node at (1,1.35) {\scriptsize $\DDownarrow \Pi^{\alpha}$};
 \draw[doubleloose] (bm)  to node[below]{\scriptsize$\chi^{g'}$} (br);
 \draw[doubletight] (mr) to  node[right]{\scriptsize$\delta \tightid_{\ten}$}(br);
 \draw[doubleeq] (tr) to (mr);
  \draw[doubleeq] (tl) to (ml);
 \node at (1.5,.5) {\scriptsize $\DDownarrow \Sigma^{\delta}$}; 
 \end{tikzpicture}
 \end{aligned}
\end{equation}


\bibliographystyle{alpha}
\bibliography{smbicat}

\end{document}